\documentclass[12pt, 3p]{elsarticle}

\usepackage{braket}             % \set
\usepackage{siunitx}
\usepackage{stmaryrd} 			% double square brackets
\usepackage{amsmath, amssymb, amsthm, mathtools, bm}
\usepackage{caption}
\usepackage{subcaption}
\usepackage[english=american]{csquotes}
\usepackage[ISO]{diffcoeff}
\usepackage{hyperref}
\usepackage[english]{babel}
\usepackage[noabbrev]{cleveref}

\biboptions{sort&compress}

\crefname{figure}{figure}{figures}

\newcommand{\eps}{\ensuremath{\varepsilon}}

\DeclarePairedDelimiterX{\abs}[1]{\lvert}{\rvert}{#1}
\DeclarePairedDelimiterX{\norm}[1]{\lVert}{\rVert}{#1}
\DeclarePairedDelimiterX{\jump}[1]{\llbracket}{\rrbracket}{#1}
\DeclarePairedDelimiterX{\dotprod}[2]{\langle}{\rangle}{#1,#2}

\newcommand{\R}{\mathbb{R}}

\newcommand{\bvec}[1]{\bm{#1}}
\newcommand{\bmat}[1]{\bm{#1}}

\graphicspath{{figures/}}

\begin{document}

    \begin{frontmatter}

        \title{A Meshfree Point Collocation Method for Elliptic Interface Problems}

        \author[itwm,uk]{Heinrich Kraus \corref{hk}}
        \author[itwm]{Jörg Kuhnert}
        \author[uk]{Andreas Meister}
        \author[itwm,ulu]{Pratik Suchde}

        \cortext[hk]{Corresponding author: heinrich.kraus@itwm.fraunhofer.de}

        % \affiliation[itwm]{
        %     organization={Fraunhofer ITWM},
        %     addressline={Fraunhofer-Platz 1},
        %     city={67663 Kaiserslautern},
        %     country={Germany}
        % }

        % \affiliation[uk]{
        %     organization={Universität Kassel},
        %     addressline={Heinrich-Plett-Straße 40},
        %     postcode={34132 Kassel},
        %     country={Germany}
        % }

        % \affiliation[ulu]{
        %     organization={University of Luxembourg},
        %     addressline={2 Av. de l'Universite},
        %     city={4365 Esch-sur-Alzette},
        %     country={Luxembourg}
        % }

        \address[itwm]{Fraunhofer ITWM, Fraunhofer-Platz 1, 67663 Kaiserslautern, Germany}
        \address[uk]{Universität Kassel, Heinrich-Plett-Straße 40, 34132 Kassel, Germany}
        \address[ulu]{University of Luxembourg, 2 Av. de l'Universite, 4365 Esch-sur-Alzette, Luxembourg}

        \begin{abstract}
            We present a meshfree generalized finite difference method for solving Poisson's equation with a diffusion coefficient that contains jump discontinuities up to several orders of magnitude. To discretize the diffusion operator, we formulate a strong form method that uses a smearing of the discontinuity; and a conservative formulation based on locally computed Voronoi cells. Additionally, we propose a novel conservative formulation for enforcing Neumann boundary conditions that is compatible with the conservative formulation of the diffusion operator. Finally, we introduce a way to switch from the strong form to the conservative formulation to obtain a locally conservative and positivity preserving scheme. The presented numerical methods are benchmarked against four test cases of varying complexity and jump magnitude on point clouds with nodes that are not aligned to the discontinuity. Our results show that the new hybrid method that switches between the two formulations produces better results than the classical generalized finite difference approach for high jumps in diffusivity.
        \end{abstract}

        %%Research highlights
        \begin{highlights}
            \item Meshfree hybrid discretization of the diffusion operator with discontinuous coefficients
            \item Formulation of a conservative scheme to enforce Neumann boundary conditions
            \item Comparison of the introduced methods with respect to convergence and performance depending on jump magnitude
        \end{highlights}

        \begin{keyword}
            Diffusion Operator \sep Discontinuous Coefficients \sep Meshfree \sep Collocation \sep Generalized Finite Difference Method
        \end{keyword}

    \end{frontmatter}

    \section{Introduction} \label{sec:introduction}
    % Motivation and future work
The present work lays the foundation for simulating phase change processes with a meshfree generalized finite difference method (GFDM) \cite{Liszka_Orkisz_1980,Gavete_Urena_Benito_Garcia_Urena_Salete_2017}. GFDMs have been proven to be broadly applicable for different applications, from solving elliptic problems \cite{Milewski_2018}, transient heat conduction problems \cite{Qu_Gu_Zhang_Fan_Zhang_2019}, shallow water equations \cite{Li_Fan_2017} to Navier--Stokes equations in industrial applications \cite{Kuhnert_2014}. Our eventual goal, beyond the scope of the present work, is to use a meshfree GFDM in a Lagrangian framework \cite{Kuhnert_2003} for phase change processes in a one-fluid model, which describes a multiphase flow in a monolithic approach without explicitly distinguishing between the different phases or tracking the interfaces between them \cite{Suchde_Kraus_Bock-Marbach_Kuhnert_2022}.

% Present methods for simulating phase change
There exist numerous approaches for simulating phase change processes.
For problems with a known and usually fixed location of the interface, a common approach is a domain decomposition where the phases are treated as separate subdomains with interface boundary conditions connecting them \cite{Davydov_Safarpoor_2021,Shao_Song_Li_2021}. Since we do not know the position of the interface in practical applications, computing a domain decomposition based on physical properties such as the temperature is required. However, this is computationally heavy and unfeasible in Lagrangian methods where a domain decomposition would need to be computed once or even multiple times per time step.
Level-set methods are another common approach for simulating phase change processes \cite{Gibou_Chen_Nguyen_Banerjee_2007}. For these types of methods, a level-set function serves as an identifier to determine to which phase a point in the computational domain belongs, and is treated as an additional variable in the system of PDEs. Level-set methods are also used for smeared-out or diffuse interfaces, but treating the level-set function as an additional variable also means that an initial level-set function needs to be set at the beginning of the simulation together with appropriate boundary conditions. However, it is often desirable for the interface between the two phases to emerge naturally without the necessity to specify the time and location of its emergence.

% Motivate that our approach leads to jumps in the diffusion coefficient and elaborate on different present discretization techniques
The monolithic one-fluid approach described earlier leads to jumps in material coefficients potentially of several orders of magnitude. To illustrate this, let us consider the case of freezing water. Near the freezing temperature, the viscosity of liquid water is $\eta \approx \SI{1e-3}{\pascal\second}$ whereas the viscosity of ice cannot be determined since it exhibits no flow. To be able to simulate the freezing process, we might use the viscosity of pitch as a replacement for the viscosity of ice. Pitch is a material that appears to be solid but in reality, it is viscoelastic with a viscosity of $\eta\approx\SI{2e8}{\pascal\second}$ \cite{Edgeworth_Dalton_Parnell_1984}. If we were to model the viscosity of ice to be of a similar magnitude, we would obtain a jump of \num{1e11} in the phase change region. Viscosity and other material coefficients, such as thermal conductivity and density, appear as a diffusivity in diffusion operators of the form $\nabla\cdot(\eta\nabla u)$. This motivated us to study Poisson's equation in divergence form with a discontinuous diffusion parameter in more detail. For such elliptic interface problems, there exist several ways of discretizing the diffusion operator with jumping diffusion coefficients.
One common approach is the use of enrichment which has been used in meshfree methods \cite{Yoon_Song_2014,Suchde_Kuhnert_2019} and finite element methods \cite{Fries_Belytschko_2010,Bordas_Natarajan_Kerfriden_Augarde_Mahapatra_Rabczuk_Pont_2011} alike. The goal of enrichment is to extend the basis of the function space, usually by discontinuous functions, and to achieve more accurate results that way.
In present literature that deals with meshfree methods, the location of the interface is usually known \cite{Davydov_Safarpoor_2021} and the considered jumps are only of a few orders of magnitude \cite{Resendiz-Flores_Saucedo-Zendejo_2018,Saucedo-Zendejo_Resendiz-Flores_2019,Trask_Perego_Bochev_2017}.

In our formulation of the GFDM, we approximate strong form differential operators based on local neighborhoods. However, it fails to ensure flux conservation in its core formulation which can lead to numerical instabilities. Several approaches have been pursued to eradicate this issue but only offer partial solutions. \Citet{Suchde_Kuhnert_Schroeder_Klar_2017} provide a weaker notion of flux conservation. Other approaches have been presented \cite{Kwan-yu_Chiu_Wang_Hu_Jameson_2012,Trask_Perego_Bochev_2017} to establishing differential operators globally that violate the local nature of differential operators and thus are too expensive to be used for time-varying point clouds.

In this paper, we introduce a method that does not depend on any a priori knowledge of the interface. It strives to identify the interface based on algebraic properties from the strong formulation of the discrete diffusion operator and switches over to a conservative method that is based on locally computed control volumes. In the finite volume method, such control volumes are based on a globally defined mesh. In contrast, the control volumes here can be computed locally \cite{Suchde_Kuhnert_Schroeder_Klar_2017}.

This paper is organized as follows. In \cref{sec:method}, we formulate Poisson's equation in divergence form and give a brief introduction to the generalized finite difference method. After that, we establish properties that our numerical methods have to fulfill; and then develop the strong form weighted-least-squares and the conservative control volume-based discretization of the diffusion operator. Finally, we introduce a way of combining these approaches into a hybrid formulation. These methods are benchmarked in \cref{sec:results} against four test cases with increasing complexity. We conclude the paper with a brief discussion of the methods and numerical results in \cref{sec:discussion} and a conclusion in \cref{sec:conclusion}.

    \section{Method Formulation} \label{sec:method}
    We consider Poisson's equation in divergence form
\begin{equation}\label{eq:poisson}
    -\nabla\cdot(\eta\nabla u) = f
\end{equation}
on a domain $\Omega\subset\R^d$ with a diffusivity $\eta = \eta(\bvec{x}) > 0$. Additionally, we set Dirichlet and Neumann boundary conditions
\begin{subequations}
    \label{eq:poisson_boundary}
    \begin{align}
        \label{eq:poisson_boundary_dirichlet} u |_{\partial\Omega_D} &= g_D \\
        \label{eq:poisson_boundary_neumann}   \left. \diffp{u}{\bvec{n}} \right\vert_{\partial\Omega_N} &= g_N
    \end{align}
\end{subequations}
with a disjoint partition of the boundary $\partial\Omega$ into  $\partial\Omega_D$ and $\partial\Omega_N$. We also use the abbreviatory notation
\begin{equation}\label{eq:diffusion_operator}
    Lu := \nabla\cdot(\eta\nabla u).
\end{equation}

Since we lay the focus on diffusion coefficients with jumps on a $d-1$-dimensional manifold $\Gamma \subset \Omega$, we are looking for solutions $u \in H^1(\Omega)$ \cite{Chipot_2009} that satisfy the jump conditions
\begin{subequations} \label{eq:jump_conditions}
    \begin{align}
        \jump{u}_\Gamma &= 0, \\
        \jump{\eta\nabla u \cdot \bvec{n}}_\Gamma &= 0.
    \end{align}
\end{subequations}
The vector $\bvec{n}$ denotes the normal vector to $\Gamma$ and the notation $\jump{\cdot}_\Gamma$ describes the jump of a function on the set $\Gamma$
\begin{align*}
    \jump{f}_\Gamma \colon &\Gamma \to \R, \quad \bvec{x} \mapsto \lim\limits_{\bvec{y}\to\bvec{x}^+} f(\bvec{y}) - \lim\limits_{\bvec{y}\to\bvec{x}^-} f(\bvec{y})
\end{align*}
where the one-sided limit from each side of the set $\Gamma$ is represented by $ \bvec{y} \to \bvec{x}^\pm$.

\subsection{Generalized Finite Difference Method}
In our formulation of the GFDM, we use a finite set of points
\begin{equation*}
    \Omega_h = \set{\bvec{x}_1, \dots, \bvec{x}_N} \subset \overline{\Omega} = \Omega \cup \partial\Omega
\end{equation*}
together with an interaction radius $h \colon \Omega \to \R$ for spatially discretizing the domain $\Omega$. For an arbitrary function $u \colon \Omega \to \R$ we use the abbreviatory notation $u_i = u(\bvec{x}_i)$ and $\bvec{u} = \begin{pmatrix*} u_1, \dots, u_N \end{pmatrix*}^T $. The interaction radius is also sometimes called smoothing length and induces a distance function $d \colon \Omega \times \Omega \to [0,\infty)$ that measures the relative distance between points in the point cloud. Usually, a scaled Euclidean distance is used, for example
\begin{align*}
    d_1(\bvec{x}_j, \bvec{x}_i) &= \frac{\norm{\bvec{x}_j - \bvec{x}_i}}{h_i}
    \shortintertext{or}
    d_2(\bvec{x}_j, \bvec{x}_i) &= \frac{2\norm{\bvec{x}_j-\bvec{x}_i}}{h_j + h_i}.
\end{align*}
The distance function enables us to define neighborhoods
\begin{equation} \label{eq:pointcloud_neighborhood}
    B_i = \set{\bvec{x}_j \in \Omega_h | d(\bvec{x}_j, \bvec{x}_i) \le 1}
\end{equation}
and their respective index sets
\begin{equation} \label{eq:pointcloud_stencil}
    S_i = \set{j | \bvec{x}_j \in B_i}
\end{equation}
for each point $\bvec{x}_i$. Neighborhoods are called symmetric if \[ j \in S_i \implies i \in S_j \] is satisfied for each point. Hence, by using $d_1$ we obtain radial but not necessarily symmetric neighborhoods while $d_2$ leads to symmetric but not necessarily radial neighborhoods. Neighborhoods are used to establish a notion of connectivity between points, based on which numerical computations are executed.

In the strong form GFDM, a linear differential operator $D^*$ is discretized at a point $\bvec{x}_i$ by using a numerical representation $D_i^*$ that is given by coefficients $c_{ij}^*$ such that
\begin{equation*}
    D^* u(\bvec{x}_i) \approx D_i^* u = \sum_{j \in S_i} c_{ij}^* u_j.
\end{equation*}
The asterisk is a placeholder for any differential operator, such as the Laplacian $\Delta$ or the partial derivative $\partial_{x_k}$ with respect to variable $x_k$.
To compute the coefficients $c_{ij}^*$, we impose exact reproducibility of monomial test functions
\begin{equation} \label{eq:optimization_monomial_constraint}
    D^* (\bvec{x} - \bvec{x}_i)^{\bvec{\alpha}} = \sum_{j \in S_i} c_{ij}^* (\bvec{x}_j - \bvec{x}_i)^{\bvec{\alpha}}
\end{equation}
for all multi-indices $\bvec{\alpha} \in \mathbb{N}_0^d$ with $\abs{\bvec{\alpha}} \le K$. Here, the usual notations $\bvec{x}^{\bvec{\alpha}} = \prod_{n=1}^d x_n^{\alpha_n}$ and $\abs{\bvec{\alpha}} = \sum_{n=1}^d \alpha_n$ are used.
This approach generally leads to fewer conditions than there are points in the neighborhood resulting in an under-determined linear system. To be able to solve the linear system, we add a minimization constraint on the coefficients
\begin{equation} \label{eq:optimization_minimization}
    \min \sum_{j \in S_i} \left( \frac{c_{ij}^*}{w_{ij}} \right)^2,
\end{equation}
with weights
\begin{equation} \label{eq:optimization_weights}
    w_{ij} = \exp\left(-2 \frac{\norm{\bvec{x}_j - \bvec{x}_i}}{h_i + h_j}\right).
\end{equation}
Other types of weight functions have been discussed in the literature, see for example \citet{Jacquemin_Tomar_Agathos_Mohseni-Mofidi_Bordas_2020}. They show that the approximation error also depends on the weight function. However, in our case, the main source of error is not the choice of weight function but the discontinuity of the diffusion coefficient. We reformulate \cref{eq:optimization_monomial_constraint,eq:optimization_minimization} into an optimization problem in matrix form
\begin{subequations} \label{eq:formal_optimization_vectors}
    \begin{alignat}{2}
        \label{eq:formal_optimization_min}
        &\min & & \norm{\bmat{W}_i^{-1} \bvec{c}_i^*}  \\
        \label{eq:formal_optimization_consistency}
        &\text{s.t.} & & \bmat{K}_i\bvec{c}_i^* = \bvec{b}_i^*
    \end{alignat}
    where the solution is given by
    \begin{equation}
        \bvec{c}_i^* = \bmat{W}_i^2 \bmat{K}_i^T (\bmat{K}_i \bmat{W}_i^2 \bmat{K}^T)^{-1} \bvec{b}_i^*.
    \end{equation}
\end{subequations}
Here the weights $w_{ij}$ are written in the $\abs{S_i} \times \abs{S_i}$ diagonal matrix $\bmat{W}_i$. Each row of matrix $\bmat{K}_i$ represents one test function applied to the points in the neighborhood and $\bvec{b}_i^*$ are the exact evaluations of the differential operators to the respective monomial test function.

\subsection{Discretization of the Diffusion Operator}
For the diffusion operator \eqref{eq:diffusion_operator}, we seek coefficients $\gamma_{ij}$ such that we obtain an approximation of the form
\begin{equation} \label{eq:divetagrad_coefficients}
    Lu(\bvec{x}_i) \approx L_i u = \sum_{j\in S_i} \gamma_{ij} u_j.
\end{equation}
The coefficients $\gamma_{ij}$ form a matrix $\bmat{G}$ to obtain the discrete Poisson's equation
\begin{equation} \label{eq:poisson_discrete}
    \bmat{G}\bvec{u} = \bvec{f}
\end{equation}
where the boundary conditions are included in $\bmat{G}$ and $\bvec{f}$. The maximum principle
\begin{equation} \label{eq:maximum_principle}
    -Lu \ge 0 \implies u \ge 0
\end{equation}
holds even for low regularity of the diffusivity $\eta$ and source term $f$ \cite{Chipot_2009}, and often the numerical schemes satisfy it discretely \cite{Li_Ito_2001}. However, GFDM in its current strong form approach has no means to guarantee the discrete maximum principle
\begin{equation} \label{eq:maximum_principle_discrete}
    \bmat{G}\bvec{u} \ge \bvec{0} \implies \bvec{u} = \bmat{G}^{-1}\bvec{f} \ge \bvec{0}
\end{equation}
without solving a global linear system \cite{Suchde_2018}.

To guarantee the maximum principle discretely while maintaining the local nature of GFDM, one way is to enforce the sign condition
\begin{subequations} \label{eq:m_matrix_conditions}
    \begin{align}
        \gamma_{ii} &< 0, \\
        \intertext{and the diagonal dominance criterion}
        \sum_{\substack{j\in S_i \\ j \ne i}} \abs{\gamma_{ij}} &\le \abs{\gamma_{ii}}
    \end{align}
\end{subequations}
upon the coefficients in \cref{eq:divetagrad_coefficients}. The conditions in \cref{eq:m_matrix_conditions} lead to an M-matrix and thus inverse-positive matrix $\bmat{G}$ that guarantees the discrete maximum principle \eqref{eq:maximum_principle_discrete}, see \citet{Seibold_2006}.

\subsubsection{Classical Operator} \label{ssec:classical}
\begin{figure}
    \centering
    \includegraphics[width=0.49\linewidth]{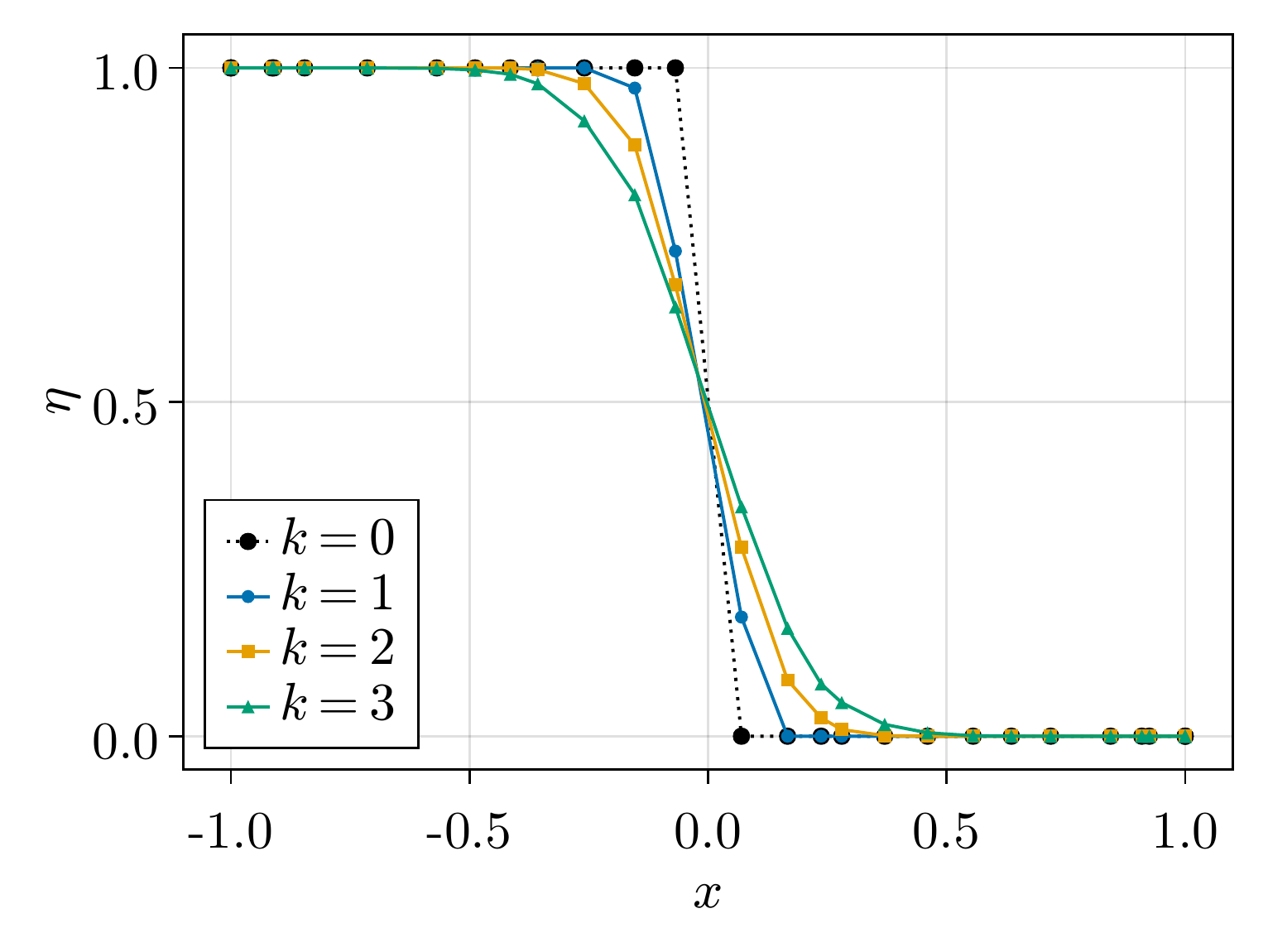}
    \caption{Qualitative smoothing comparison of a piecewise constant function and $k$ smoothing cycles.}
    \label{fig:comparison_smoothing_2d}
\end{figure}

\begin{figure}
    \begin{subfigure}{0.49\textwidth}
        \centering
        \includegraphics[width=\linewidth]{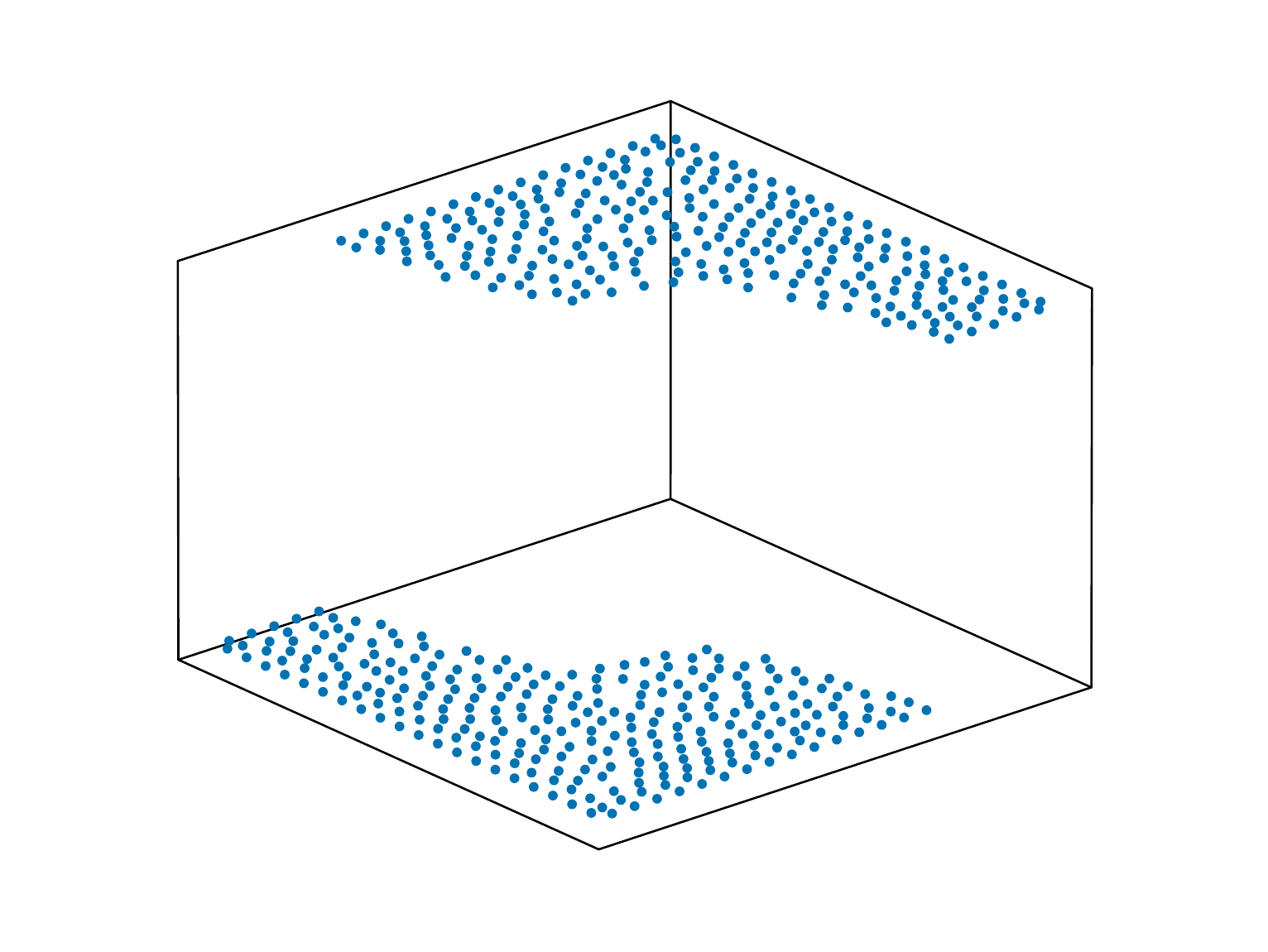}
    \end{subfigure}
    \begin{subfigure}{0.49\textwidth}
        \centering
        \includegraphics[width=\linewidth]{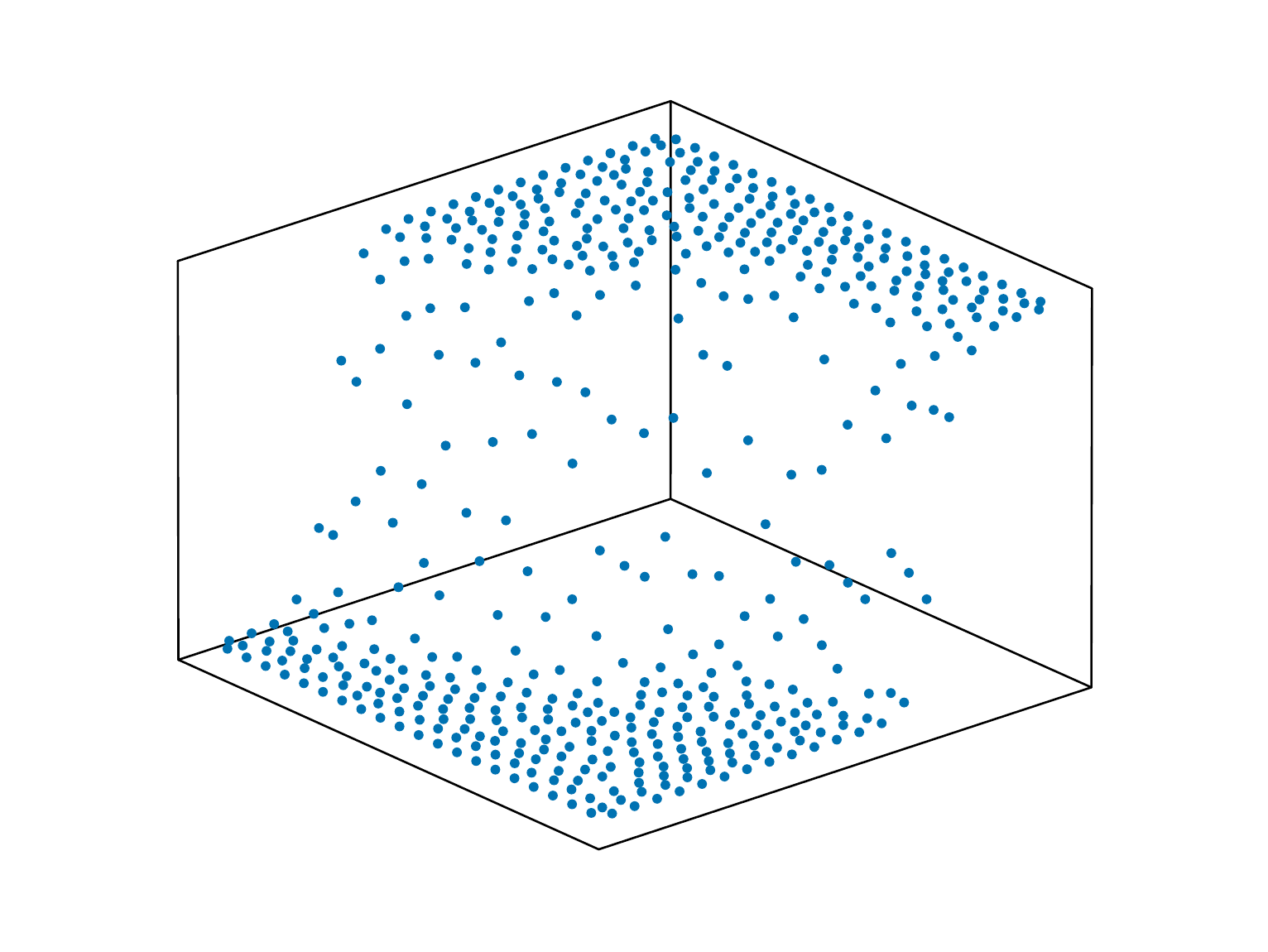}
    \end{subfigure}
    \caption{Comparison of discontinuous function (left) and its smoothed counterpart (right) with $2$ smoothing cycles.}
    \label{fig:comparison_smoothing}
\end{figure}

For the strong form GFDM, we assume that the diffusivity parameter $\eta$ is sufficiently smooth. To be able to use the strong form for discontinuous diffusivities, we smooth the diffusivity $\eta$ by using weights \[ s_{ij} = \exp\left(-3 \frac{\norm{\bvec{x}_j - \bvec{x}_i}^2}{h_i^2}\right). \] in a convex combination of function values over the neighborhood
\[ \tilde{\eta}_i = \frac{\sum_{j \in S_i} s_{ij} \eta_j}{\sum_{j \in S_i} s_{ij}}. \]
The smoothing process can be applied multiple times, leading to an incrementally increasing smearing of the jump. \Cref{fig:comparison_smoothing_2d} shows a comparison of a piecewise constant function with a different number of smoothing cycles and in \cref{fig:comparison_smoothing}, a discontinuous function with its smoothed equivalent is depicted.

For a smoothed diffusivity $\tilde{\eta}$, we apply the product rule on the diffusion operator \eqref{eq:diffusion_operator} to obtain
\begin{equation} \label{eq:divetagrad_classical}
    Lu \approx \nabla\cdot(\tilde{\eta}\nabla u) = \nabla\tilde{\eta}\cdot\nabla u + \tilde{\eta}\Delta u.
\end{equation}
Alternatively, it is sometimes advantageous to use the scaled diffusivity parameter $\mu = \log\eta$, especially when there are large jumps in $\eta$. Scaling the diffusivity that way yields, after potential smoothing steps, the diffusion operator
\begin{equation} \label{eq:divetagrad_classical_scaled}
    Lu = \nabla\cdot(\exp(\mu)\nabla u) \approx \nabla\cdot(\exp(\tilde{\mu})\nabla u) = \exp(\tilde{\mu}) \left( \nabla\tilde{\mu}\cdot\nabla u + \Delta u \right).
\end{equation}
By using monomial test functions of degree $1$ or higher, the gradient of either $\tilde{\eta}$ or $\tilde{\mu}$ has to be computed numerically.  The numerical computation of a smoothed and thus distorted diffusivity can be a source of inaccuracy for this method.

The coefficients $\gamma_{ij}$ are then computed by solving the optimization problem \eqref{eq:formal_optimization_vectors}. Since the optimization problem does not impose diagonal dominance, this can lead to non-diagonally dominant rows and, in some cases, instabilities in the numerical solution \cite{Suchde_2018}. To reduce numerical instabilities, a correction technique has been developed where a correction vector $\bvec{\xi}_i$ is computed that lies in the null space of the test function space
\[ \sum_{j\in S_i} \xi_{ij} (\bvec{x}_j - \bvec{x}_i)^{\bvec{\alpha}} = 0 \]
for all considered monomials. To obtain a non-trivial zero operator, an additional constraint $\xi_{ii} = 1$ is added and optimization \cref{eq:formal_optimization_vectors} is solved. With the correction vector $\bvec{\xi}_i$, we are able to define a new operator
\begin{equation} \label{eq:diagonal_dominance_correction}
    \bvec{a}_i = \bvec{\gamma}_i + \alpha \bvec{\xi}_i,
\end{equation}
that satisfies the same monomial reproducibility properties as $\bvec{\gamma}_i$ for each $\alpha\in\R$. In the next step, we minimize the function
\[ \phi(\alpha) = \sum_{j\in S_i} \frac{(\gamma_{ij} + \alpha \xi_{ij})^2}{(\gamma_{ii} + \alpha \xi_{ii})^2} \]
that attains its global minimum at
\[ \alpha_{\text{min}} = \frac{\gamma_{ii} \dotprod{\bvec{\gamma}_i}{\bvec{\xi}_i} - \xi_{ii} \dotprod{\bvec{\gamma}_i}{\bvec{\gamma}_i} }{\xi_{ii} \dotprod{\bvec{\gamma}_i}{\bvec{\xi}_i} - \gamma_{ii}\dotprod{\bvec{\xi}_i}{\bvec{\xi}_i} }. \]
This correction technique does not guarantee diagonally dominant operators and there are numerous examples of neighborhoods where this technique does not provide diagonally dominant rows, some of which are shown in \ref{app:node_selection}. Nevertheless, it has been shown that the stability of the discrete Laplace operator increases by using this technique \cite{Suchde_2018}.

\subsubsection{Conservative Operator} \label{ssec:conservative}
For the conservative approach, we omit optimization problem \eqref{eq:formal_optimization_vectors} and use Voronoi cells
\begin{equation} \label{eq:voronoi}
    \Omega_i = \Set{\bvec{x} \in \Omega | \norm{\bvec{x} - \bvec{x}_i} < \norm{\bvec{x} - \bvec{x}_j} \text{ for each } \bvec{x}_j \in \Omega_h \setminus \set{\bvec{x}_i}}
\end{equation}
instead. For a Voronoi cell $\Omega_i$, the area in 2D and volume in 3D is given by \[ \abs{\Omega_i} = \frac{1}{2d} \sum_{j\in S_i} \abs{\Gamma_{ij}} \lVert x_j - x_i \rVert \] where the boundary of the cell is represented as $\partial\Omega_i = \cup_j \Gamma_{ij}$ with $\Gamma_{ij} = \overline{\Omega}_i \cap \overline{\Omega}_j$. The size of $\Gamma_{ij}$, a line segment in 2D and a polygon in 3D, is denoted by $\abs{\Gamma_{ij}}$.

For an inner point $\bvec{x}_i \in \Omega$, we use a cell average \cite{Eymard_Gallouet_Herbin_2000} to compute the discrete diffusion operator
\begin{equation*}
    \abs{\Omega_i} L_i u = \int_{\Omega_i}\nabla\cdot(\eta\nabla u) \dl2\bvec{x}.
\end{equation*}
Applying the divergence theorem and exploiting the composition of the boundary of $\Omega_i$ into the faces $\Gamma_{ij}$ leads on to
\begin{equation*}
    \abs{\Omega_i} L_i u = \oint_{\partial\Omega_i} \eta\diffp{u}{\bvec{n}} \dl2S
                         = \sum_{j \in S_i} \int_{\Gamma_{ij}} \eta\diffp{u}{\bvec{n}} \dl2S.
\end{equation*}
This is a representation of the diffusion operator in terms of fluxes over the surfaces of the Voronoi cells. Note that for each point $\bvec{x}_j \in B_i$ that is not a Voronoi neighbor of $\bvec{x}_i$, we use $\abs{\Gamma_{ij}} = 0$. Applying Gaussian quadrature to each face leads to
\begin{equation*}
    \int_{\Gamma_{ij}} \eta \diffp{u}{\bvec{n}} \dl2S \approx \abs{\Gamma_{ij}} \eta_{ij} \frac{u_j - u_i}{d_{ij}}
\end{equation*}
where $\eta_{ij} \approx \eta(\bvec{x}_{ij})$ is an approximation of $\eta$ at the midpoint \[\bvec{x}_{ij} = \frac{\bvec{x}_i + \bvec{x}_j}{2}\] of the line segment \[L_{ij} = \set{\bvec{x}_i + \theta (\bvec{x}_j - \bvec{x}_i) \mid \theta \in [0, 1]}\] between $\bvec{x}_i$ and $\bvec{x}_j$ and $d_{ij} = \norm{\bvec{x}_j - \bvec{x}_i}$ is the length of that line segment. Note that if the line segment $L_{ij}$ and the surface $\Gamma_{ij}$ intersect, then they intersect at $\bvec{x}_{ij}$.

Finally, the coefficients $\gamma_{ij}$ for the discrete diffusion operator read
\begin{subequations}\label{eq:coeffs_divetagrad_conservative}
    \begin{align}
        \gamma_{ij} &= \frac{\abs{\Gamma_{ij}}}{\abs{\Omega_i}} \frac{\eta_{ij}}{d_{ij}}
    \shortintertext{where $i \ne j$ and for the diagonal entries}
        \gamma_{ii} &= -\sum_{\substack{j \in S_i \\ j \ne i}} \frac{\abs{\Gamma_{ij}}}{\abs{\Omega_i}} \frac{\eta_{ij}}{d_{ij}}.
    \end{align}
\end{subequations}
These coefficients trivially fulfill diagonal dominance and the sign restriction from \cref{eq:m_matrix_conditions} such that they automatically lead to diagonally dominant rows and thus fulfill the discrete maximum principle \eqref{eq:maximum_principle_discrete}.

In the case of meshfree GFDM, we have a fixed distribution of points $\bvec{x}_i \in \Omega$ based on which we can calculate the Voronoi cells. The quality requirements imposed on point clouds are rather loose \cite{Seifarth_2018} in comparison to the strong restrictions on the meshes used for mesh-based methods \cite{Knupp_2001}. Furthermore, the Voronoi cells used here are constructed locally from the neighborhoods $B_i$ and can thus easily be used for moving point clouds as well.

\paragraph{Calculation of $\eta_{ij}$}

To calculate the reconstructions $\eta_{ij}$, we use the harmonic mean
\[ \eta_{ij} = \frac{2}{\frac{1}{\eta_i} + \frac{1}{\eta_j}}. \]

Alternatively,  a discrete approximation operator can be computed with the optimization problem \eqref{eq:formal_optimization_vectors} that discretizes the identity mapping
\[ I_{ij} \eta = \eta(\bvec{x}_{ij}). \]
Similar operators are used for approximating function values at new points where the values were previously unknown. But for diffusion coefficients with large discontinuities, this can lead to oscillations and even negative reconstructions $\eta_{ij} < 0$.

\paragraph{Conservative Neumann Boundary Condition}

Applying the cell average idea to boundary cells $\Omega_i$ that are associated with Neumann boundary points $\bvec{x}_i \in \partial\Omega_N$, we obtain
\begin{equation*}
    \int_{\Omega_i} \nabla\cdot(\eta\nabla u) \dl2S = \int_{\Gamma_i} \eta\diffp{u}{\bvec{n}} \dl2S + \sum_{j \in S_i}\int_{\Gamma_{ij}} \eta\diffp{u}{\bvec{n}} \dl2S
\end{equation*}
where $\Gamma_i \subset \partial\Omega$ is a piece of the domain boundary associated with the boundary point $\bvec{x}_i$. Using \cref{eq:poisson,eq:poisson_boundary_neumann} finally leads to
\begin{equation*}
    -\int_{\Omega_i} f \dl2\bvec{x} - \int_{\Gamma_i} \eta g_N \dl2S = \sum_{j \in S_i}\int_{\Gamma_{ij}} \eta\diffp{u}{\bvec{n}} \dl2S.
\end{equation*}
Dividing by $\abs{\Omega_i}$ and using the coefficients in \cref{eq:coeffs_divetagrad_conservative}, we find
\begin{equation} \label{eq:coeffs_divetagrad_conservative_neumann}
    -f_i - \frac{\abs{\Gamma_i}}{\abs{\Omega_i}} \eta_i g_N(\bvec{x}_i) = \sum_{j \in S_i} \gamma_{ij} u_j.
\end{equation}

\paragraph{Discrete Divergence Theorem}
For the divergence theorem \[ \int_\Omega \nabla\cdot(\eta\nabla u) \dl2\bvec{x} = \oint_{\partial\Omega} \eta \diffp{u}{\bvec{n}} \dl2S \] to hold discretely, the coefficients $\gamma_{ij}$ from \cref{eq:divetagrad_coefficients} need to suffice to the column sum condition
\begin{equation} \label{eq:gauss_column_sum}
    \sum_{i = 1}^N \abs{\Omega_i} \gamma_{ij} = 0
\end{equation}
for interior points $\bvec{x}_j \in \Omega$ \cite{Suchde_2018}. The cell-based approach, given by \cref{eq:coeffs_divetagrad_conservative,eq:coeffs_divetagrad_conservative_neumann}, automatically satisfies \cref{eq:gauss_column_sum} and fulfills the discrete divergence theorem on the entire domain. However, enforcing \cref{eq:gauss_column_sum} in the row-based discrete operator calculation from \cref{ssec:classical}, requires us to solve one huge optimization problem instead of solving many small optimization problems as in \cref{eq:formal_optimization_vectors}.

\subsubsection{Hybrid} \label{ssec:hybrid}
The computations of the Voronoi cells \eqref{eq:voronoi}, even when done locally, can be costly and may not be necessary at all points. If the diffusion coefficient $\eta$ is smooth in a certain region, it is appropriate to use the classical approach from \cref{ssec:classical} to obtain satisfactory results and only use the conservative method in regions where the diffusivity is discontinuous. To identify regions where $\eta$ is discontinuous, we use the coefficients $\gamma_{ij}$ obtained from the strong form method and define the diagonal dominance error
\begin{equation} \label{eq:hybrid_condition}
    \sigma_i = \max\biggl( \sum_{\substack{j\in S_i \\ j \ne i}} \frac{\abs{\gamma_{ij}}}{\abs{\gamma_{ii}}} - 1, \gamma_{ii}, 0 \biggr).
\end{equation}
In the case that $\sigma_i \le 0$, we use the coefficients $\gamma_{ij}$ to discretize the diffusion operator, and for $\sigma_i > 0$ the coefficients are recomputed using the conservative formulation. This leads to a hybrid method that preserves the discrete maximum principle \eqref{eq:maximum_principle_discrete}.

Better results can be achieved by using the conservative formulation not only for those points where $\sigma_i > 0$ but also for each neighboring point $j \in S_i$. This results in a conservative scheme for symmetric neighborhoods. Note that condition \eqref{eq:hybrid_condition} checks only the diagonal dominance of a row. It can thus be triggered not only by high jumps in the diffusivity but also by local irregularities in the point cloud. Examples for such node configurations are presented in \ref{app:node_selection}.

Other jump identification techniques exist where the function is directly screened for discontinuities, see for example \citet{Allasia_Besenghi_Cavoretto_2009,Gutzmer_Iske_1997}. In their formulation, global properties of the examined function are used and thus these methods only work for one jump or multiple jumps of similar magnitudes.

    \section{Numerical Results} \label{sec:results}
    In the preceding sections, we presented numerical methods for approximating the diffusion operator. With that, we introduced several parameters for computing the discrete differential operator. Parameters of special interest are
\begin{itemize}
    \item number of smoothing cycles of the diffusivity,
    \item scaling of the diffusion coefficient
\end{itemize}
for the strong form method. For the hybrid method, we will take a closer look at
\begin{itemize}
    \item the possibility of extending the conservative formulation to the neighbors of points where a lack of diagonal dominance was discovered,
    \item conservative Neumann boundaries.
\end{itemize}
In this section, we compare the different methods for select test cases in terms of convergence and susceptibility to jumps in the diffusivity. A performance comparison between the strong form method and the hybrid method is given in \ref{app:performance}.

\begin{figure}
    \newcommand{\subfigurehspace}{.49\linewidth}
    \begin{subfigure}{\subfigurehspace}
        \centering
        \includegraphics[width=\linewidth]{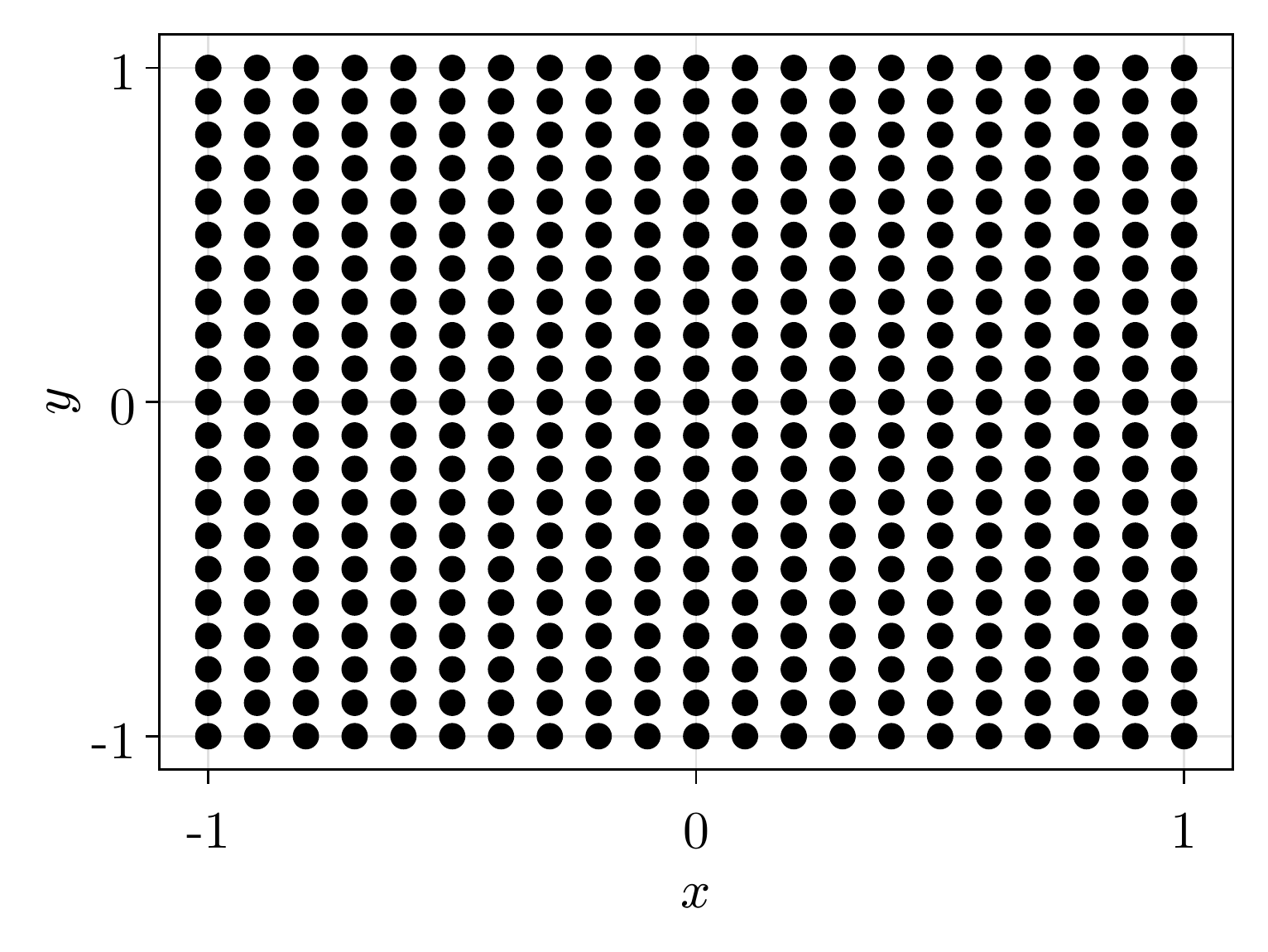}
        \caption{Uniform}
        \label{fig:comparison_pointclouds_uniform}
    \end{subfigure}
    \begin{subfigure}{\subfigurehspace}
        \centering
        \includegraphics[width=\linewidth]{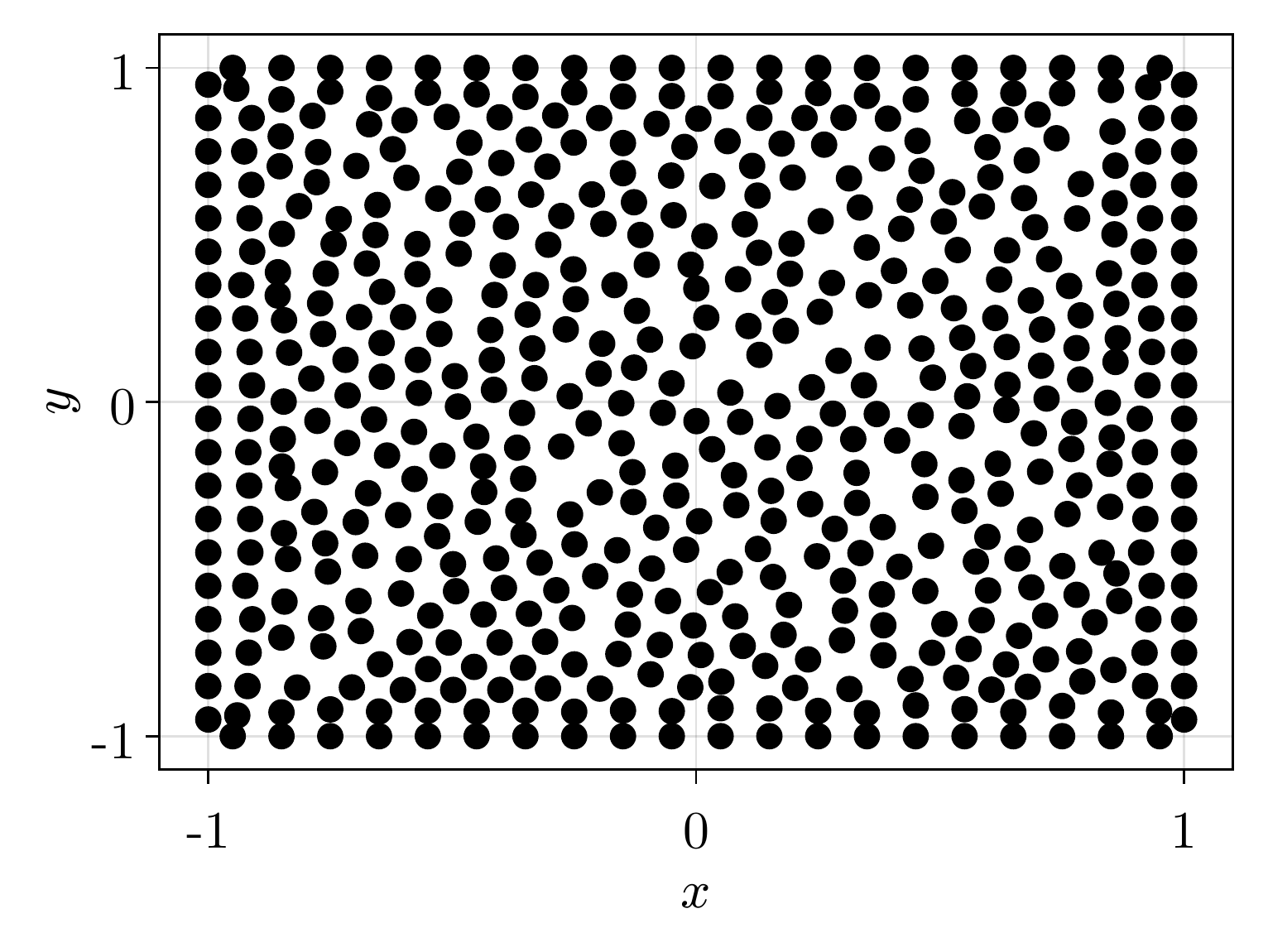}
        \caption{Irregular}
        \label{fig:comparison_pointclouds_distorted}
    \end{subfigure}
    \caption{Comparison of point cloud types}
    \label{fig:comparison_pointclouds}
\end{figure}

On the rectangle $\Omega = [-1, 1] \times [-1, 1]$ we use two different types of point clouds, uniform and irregular point clouds, with increasing refinement levels $h_k = \frac{2^{-k}}{5}$ for $k = 0, \dots, 5$. A comparison of the different point cloud types with similar refinements is shown in \cref{fig:comparison_pointclouds}.
The irregular point clouds are generated with an advancing front method \cite{Lohner_Onate_1998}, using the software MESHFREE developed by Fraunhofer ITWM \cite{MESHFREE}. We set coefficients $r_{\min} = 0.25$ and $r_{\max} = 0.45$ such that the relative distance to the closest neighbor is bounded below by
\[ r_{\min} \le \min_{\substack{j \in S_i \\ j \ne i}} \frac{\norm{\bvec{x}_j - \bvec{x}_i}}{h_i} \]
for each point $\bvec{x}_i$. Additionally, we require each point in the domain $\bvec{x}\in\Omega$ to fulfill the intersection condition
\[ \overline{B(\bvec{x}; r_{\max} h(\bvec{x}))} \cap \Omega_h \ne \emptyset \]
that defines the maximum size of holes that appear in the point cloud. By $B(\bvec{x}; r)$, we denote a ball around $\bvec{x}$ with radius $r$ with respect to the Euclidean metric.
We used the irregular point clouds to test if and to what extent they influence the numerical results compared to uniform point clouds. In our experiments, we compared the numerical results from the different point cloud types and could observe a small difference in the discretization error. Since the convergence behavior was similar, we will only show results from irregular point clouds as seen in \cref{fig:comparison_pointclouds_distorted} and specifically state where a uniform point cloud was used. Because the point clouds are generated independently of the diffusivity, we are able to show that our method works well on point clouds that do not conform to the interface.

To simplify notation, we decompose the boundary of $\Omega$ \[ \partial\Omega = \partial\Omega_L \cup \partial\Omega_R \cup \partial\Omega_T \cup \partial\Omega_B \] into the left, right, top and bottom part respectively.

To analyze the convergence behavior of the presented methods, we compute the \emph{relative} $L^2$ error of the numerical solution $u_h$ to the analytical solution $u$ \[ \norm{u - u_h} = \frac{\norm{u-u_h}_{L^2}}{\norm{u}_{L^2}}, \] where the discrete $L^2$ norm for a function $f \colon \Omega \to \R$ is defined by \[ \norm{f}_{L^2} = \left( \sum_{i = 1}^N v_i \abs{f(\bvec{x}_i)}^2 \right)^{1/2} \] and $v_i > 0$ are areas in 2D or volumes in 3D associated with each point $\bvec{x}_i$.

\subsection{Two-Strip Problem} \label{ssec:two_strip}
\begin{figure}
    \newcommand{\subfigurehspace}{0.49\linewidth}
    \begin{subfigure}{\subfigurehspace}
        \centering
        \includegraphics[width=\linewidth]{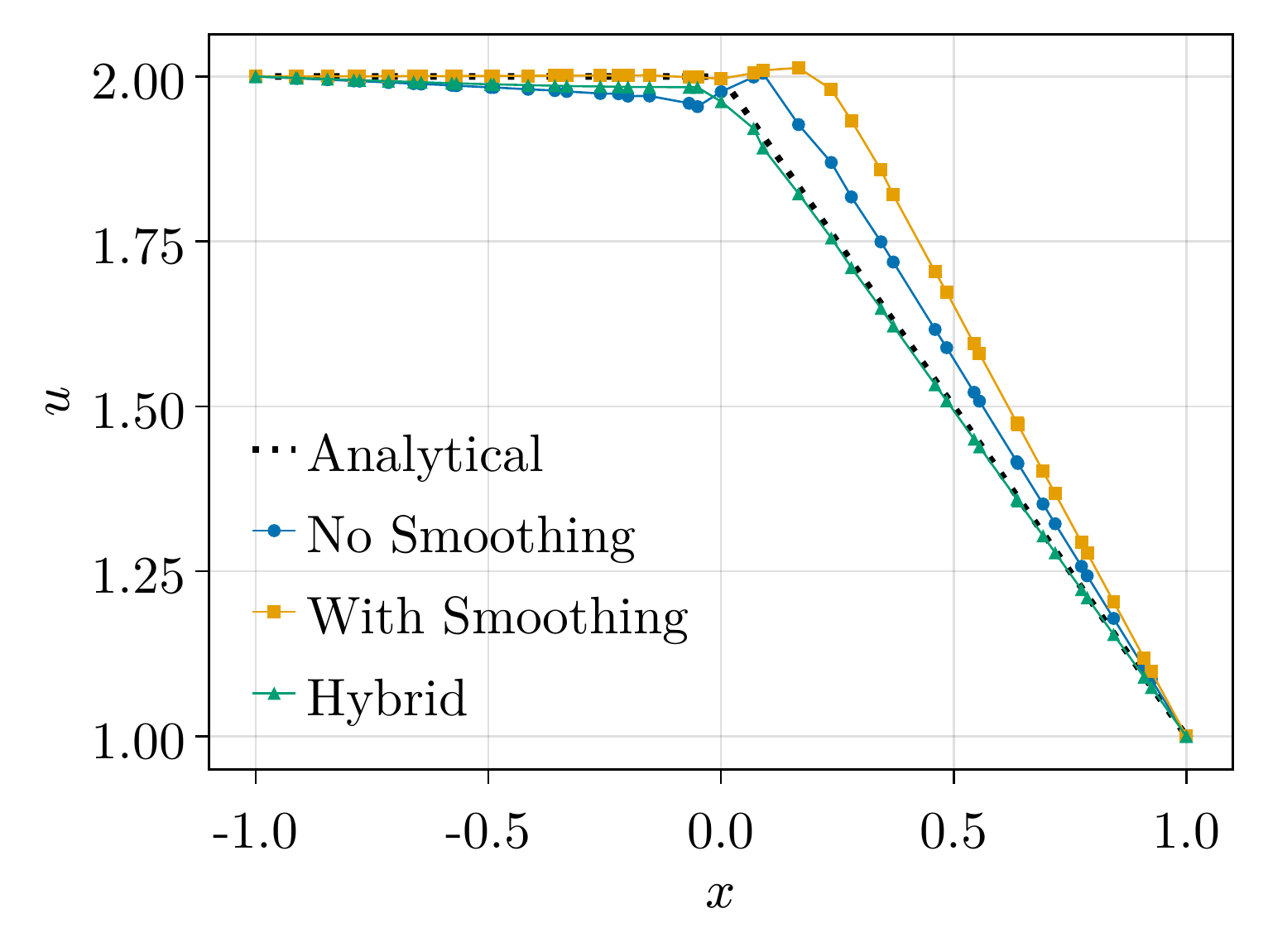}
        \caption{$\delta\eta = \num{1e-8}$}
        \label{sfig:two_strip_solutions_1e-8}
    \end{subfigure}
    \begin{subfigure}{\subfigurehspace}
        \centering
        \includegraphics[width=\linewidth]{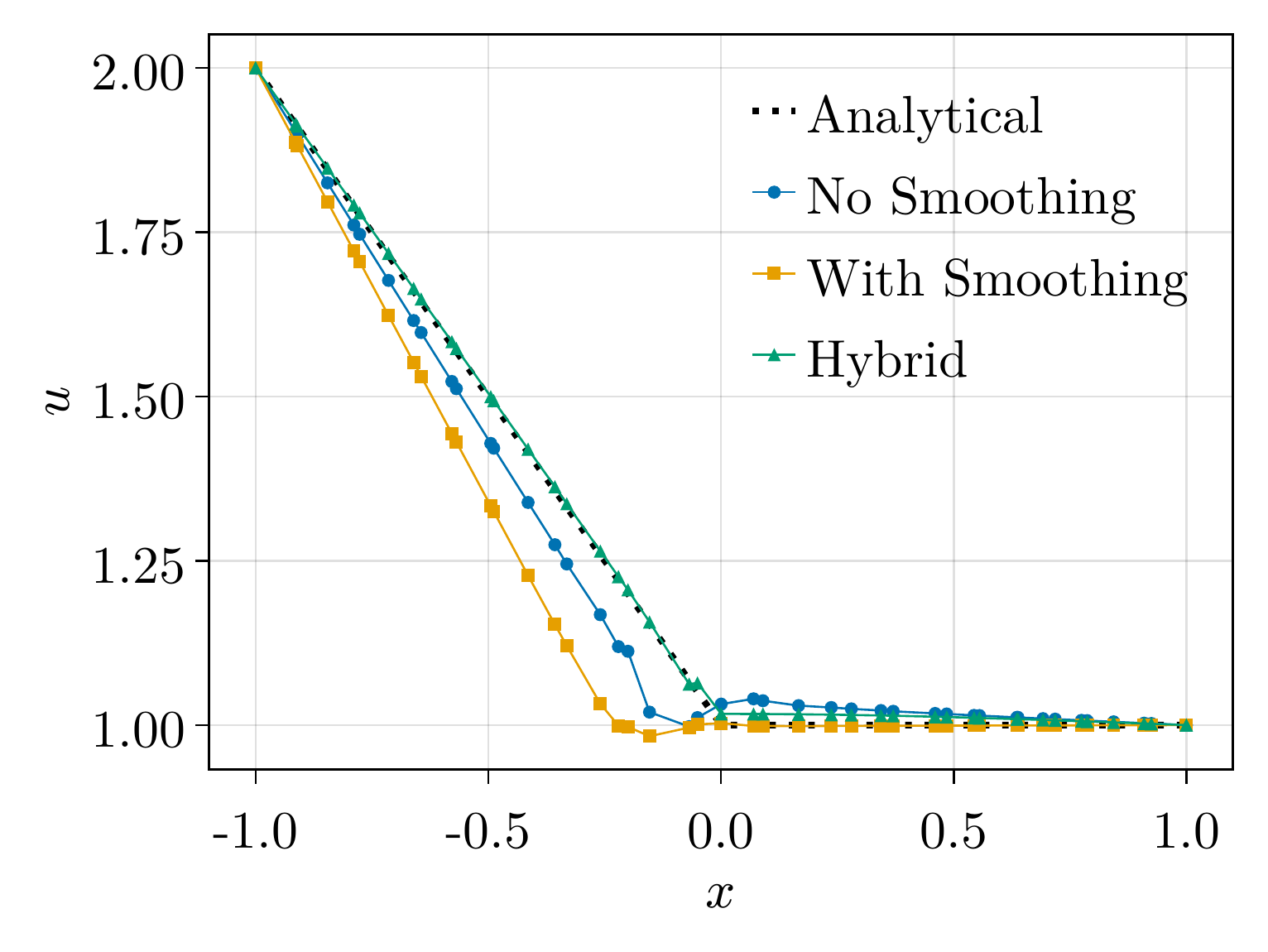}
        \caption{$\delta\eta = \num{1e8}$}
        \label{sfig:two_strip_solutions_1e+8}
    \end{subfigure}
    \caption{Solution profile at $y \approx 0$ to the two-strip test case with a jump of $\delta\eta$.}
    \label{fig:two_strip_solutions}
\end{figure}

For the two-strip problem, we solve the homogeneous Poisson's equation and define the diffusivity as a piecewise constant function
\[ \eta(x,y) = \begin{cases}
    \eta_L, & x < 0, \\
    \eta_R, & x > 0
\end{cases} \]
where $\eta_L$, $\eta_R > 0$.
The Dirichlet and Neumann boundary conditions
\begin{align*}
    u|_{\partial\Omega_L} &= 2 \\
    u|_{\partial\Omega_R} &= 1 \\
    \diffp{u}{\bvec{n}} \bigg\vert_{\partial\Omega_T} &= 0\\
    \diffp{u}{\bvec{n}} \bigg\vert_{\partial\Omega_B} &= 0
\end{align*}
lead to an analytical solution that is monotonously decreasing and piecewise linear in $x$-direction and constant in $y$-direction. We characterize a two-strip problem uniquely by the jump
\[ \delta\eta = \frac{\eta_R}{\eta_L}. \] The solution profile of the analytical solution can be seen in \cref{fig:two_strip_solutions} for a jump of $\delta\eta = \num{1e-8}$ and $\delta\eta = \num{1e8}$ respectively.

\begin{figure}
    \newcommand{\subfigurehspace}{.49\linewidth}
    \setlength{\belowcaptionskip}{1.5\baselineskip}
    \begin{subfigure}{\subfigurehspace}
        \centering
        \includegraphics[width=\linewidth]{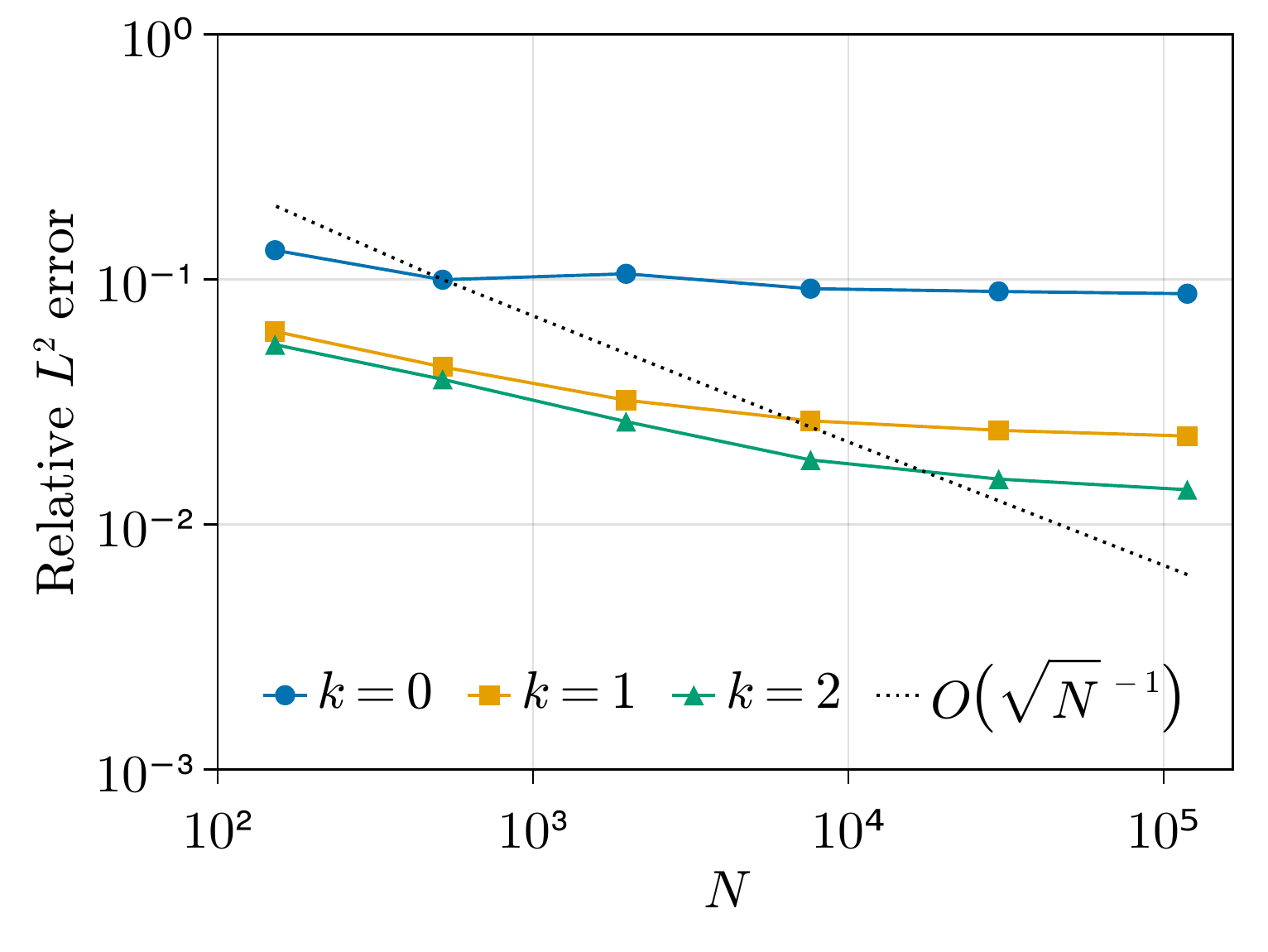}
        \caption{No scaling, $\delta\eta = \num{1e1}$}
        \label{sfig:two_strip_strong_form_smoothing_noscaling_1e1}
    \end{subfigure}
    \begin{subfigure}{\subfigurehspace}
        \centering
        \includegraphics[width=\linewidth]{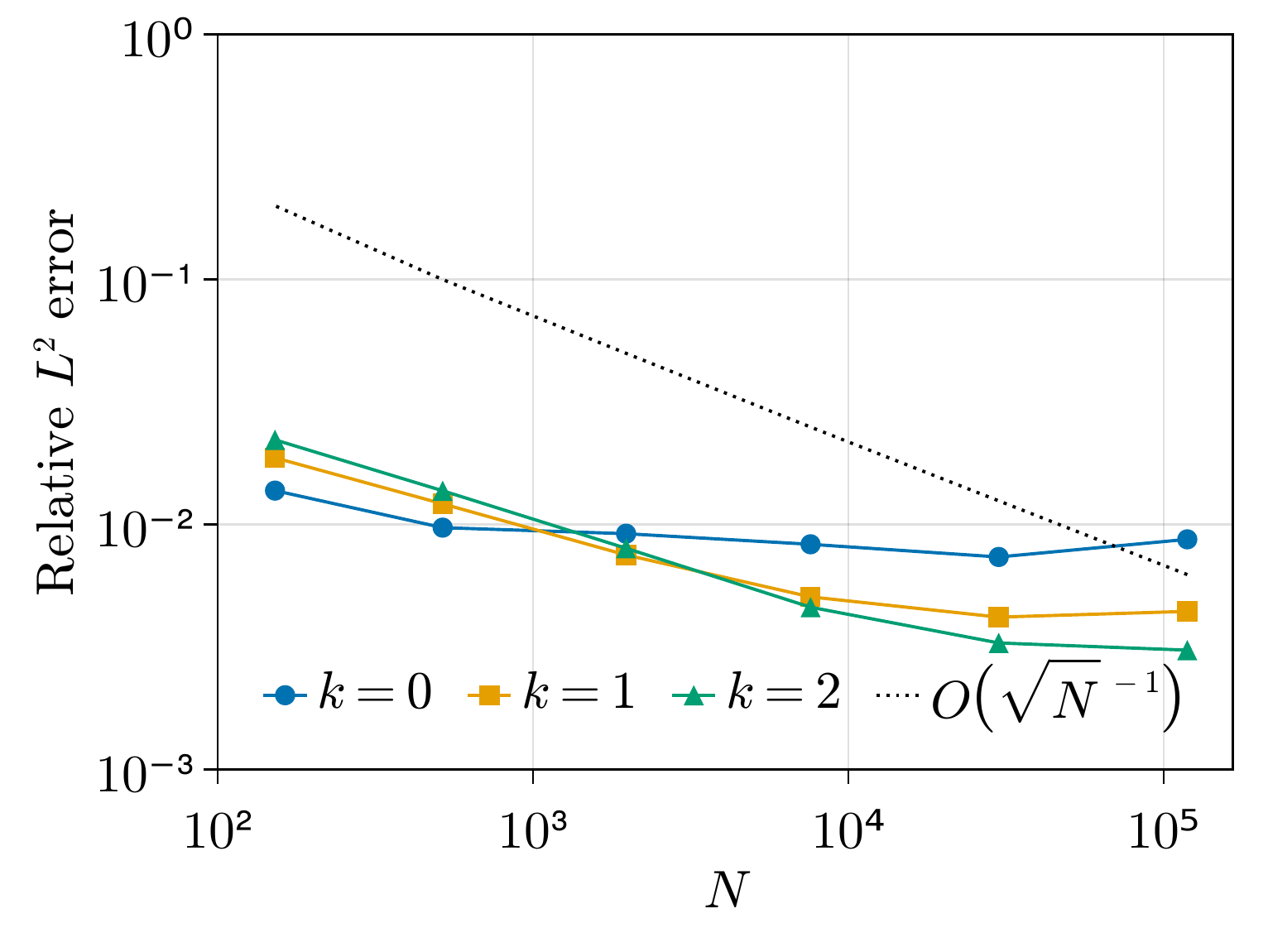}
        \caption{Logarithmic scaling, $\delta\eta = \num{1e1}$}
        \label{sfig:two_strip_strong_form_smoothing_scaling_1e1}
    \end{subfigure}

    \par

    \begin{subfigure}{\subfigurehspace}
        \centering
        \includegraphics[width=\linewidth]{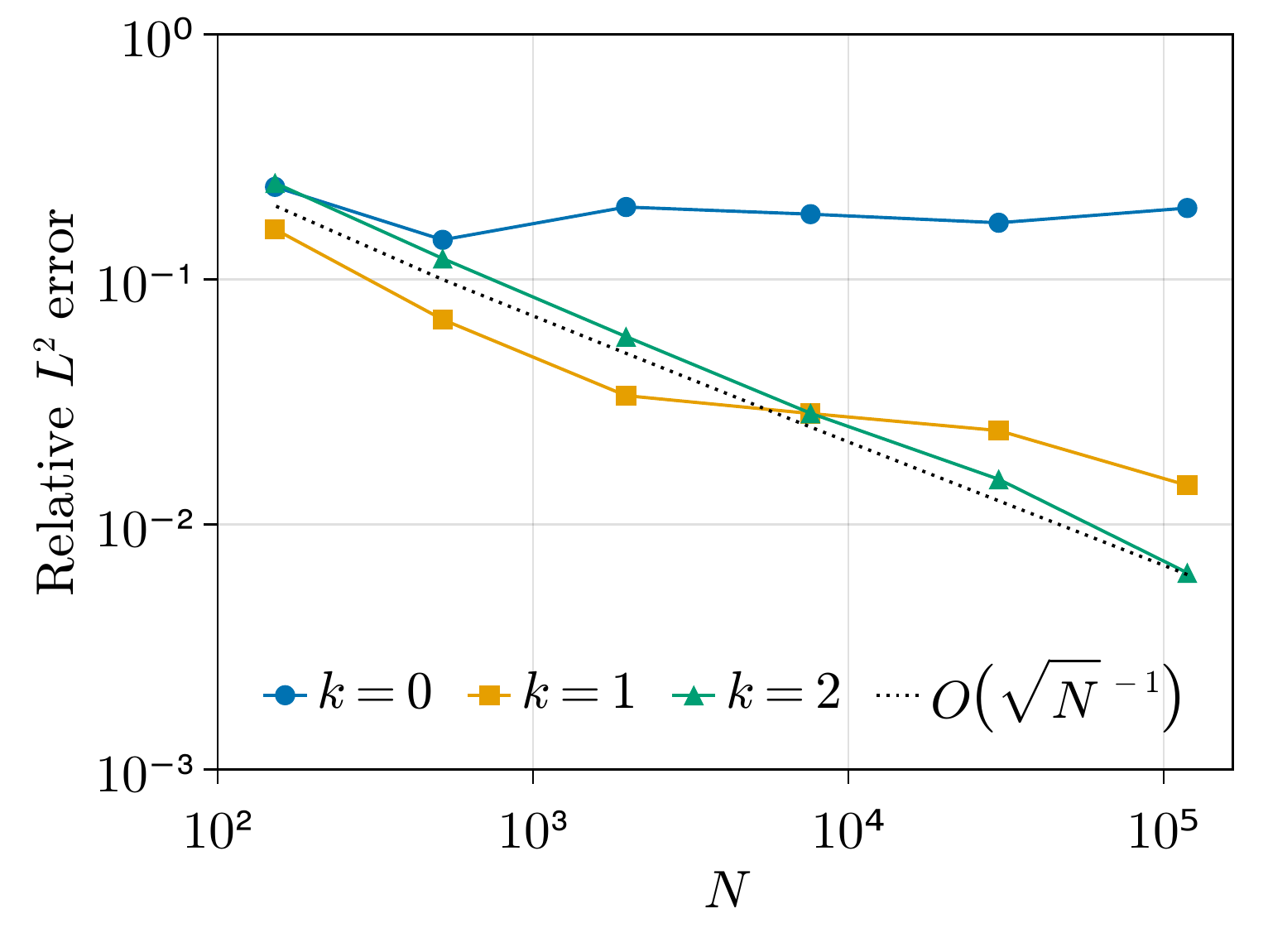}
        \caption{No scaling, $\delta\eta = \num{1e10}$}
        \label{sfig:two_strip_strong_form_smoothing_noscaling_1e10}
    \end{subfigure}
    \begin{subfigure}{\subfigurehspace}
        \centering
        \includegraphics[width=\linewidth]{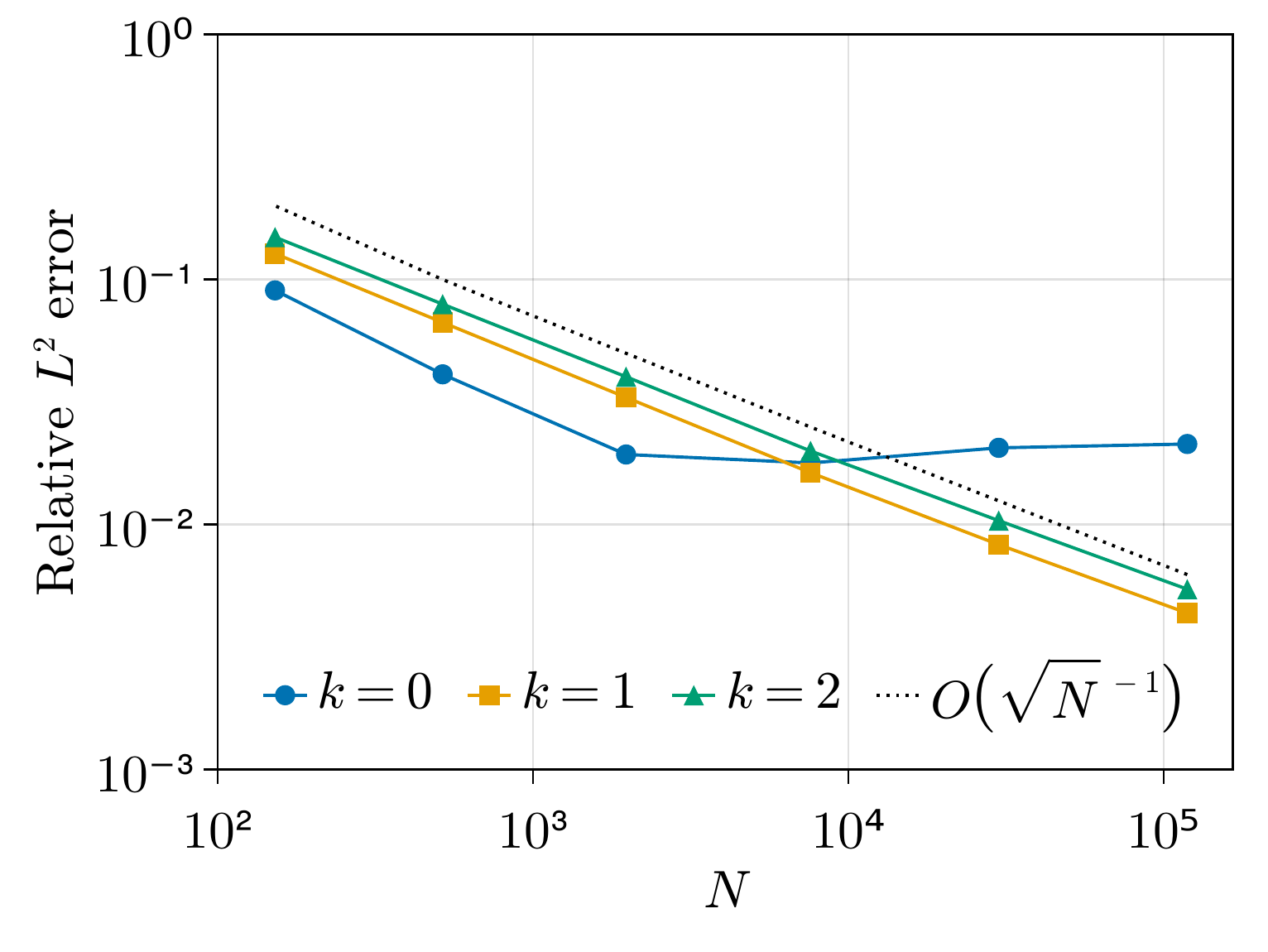}
        \caption{Logarithmic scaling, $\delta\eta = \num{1e10}$}
        \label{sfig:two_strip_strong_form_smoothing_scaling_1e10}
    \end{subfigure}
    \caption{Convergence plots of the classical strong form method for the two-strip test case with and without logarithmic scaling with different jumps $\delta\eta$ and number of smoothing cycles $k$.}
    \label{fig:two_strip_strong_form_smoothing}
\end{figure}

\begin{figure}
    \newcommand{\subfigurehspace}{.49\linewidth}

    \begin{subfigure}{\subfigurehspace}
        \centering
        \includegraphics[width=\linewidth]{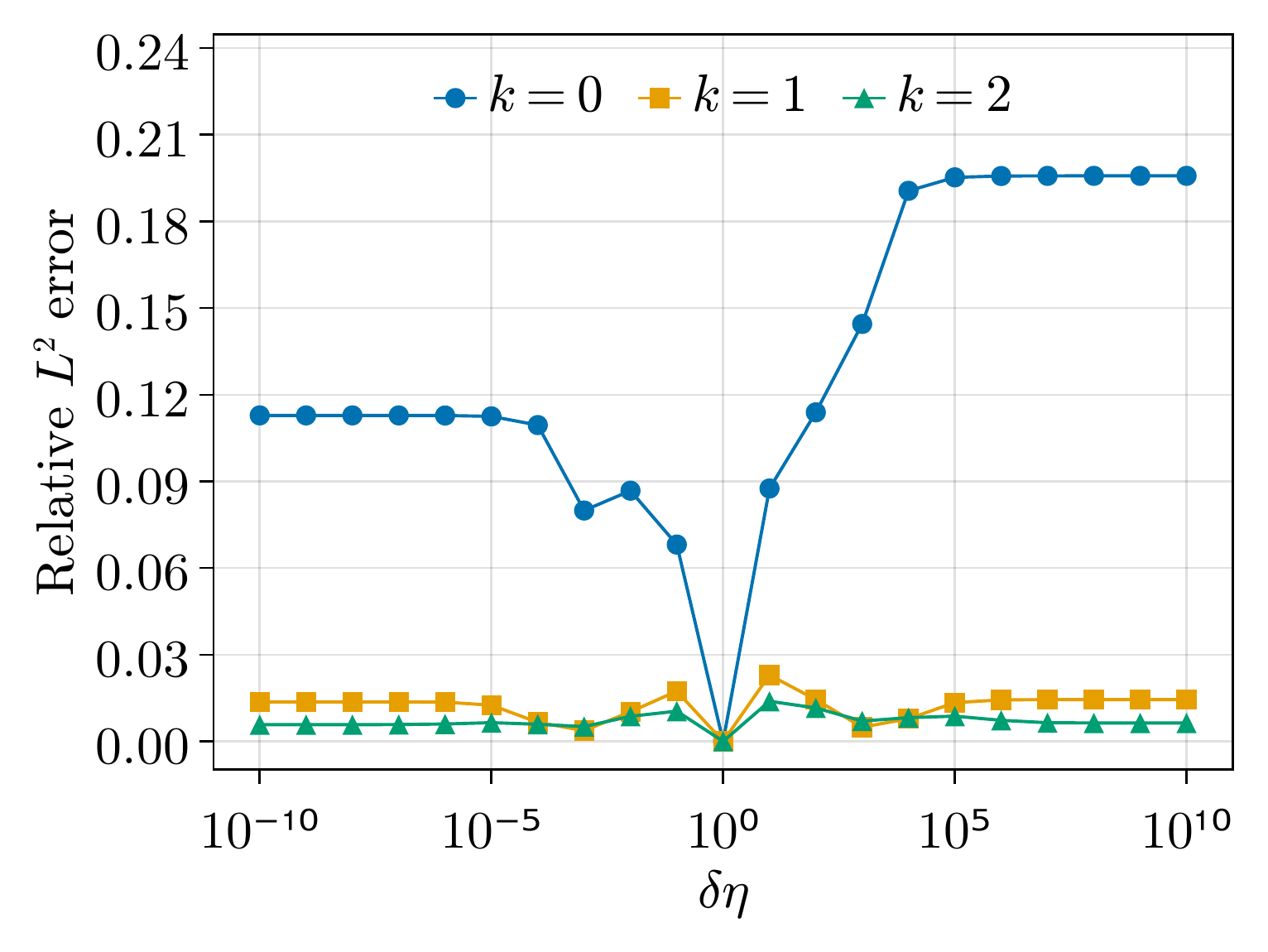}
        \caption{No scaling}
    \end{subfigure}
    \begin{subfigure}{\subfigurehspace}
        \centering
        \includegraphics[width=\linewidth]{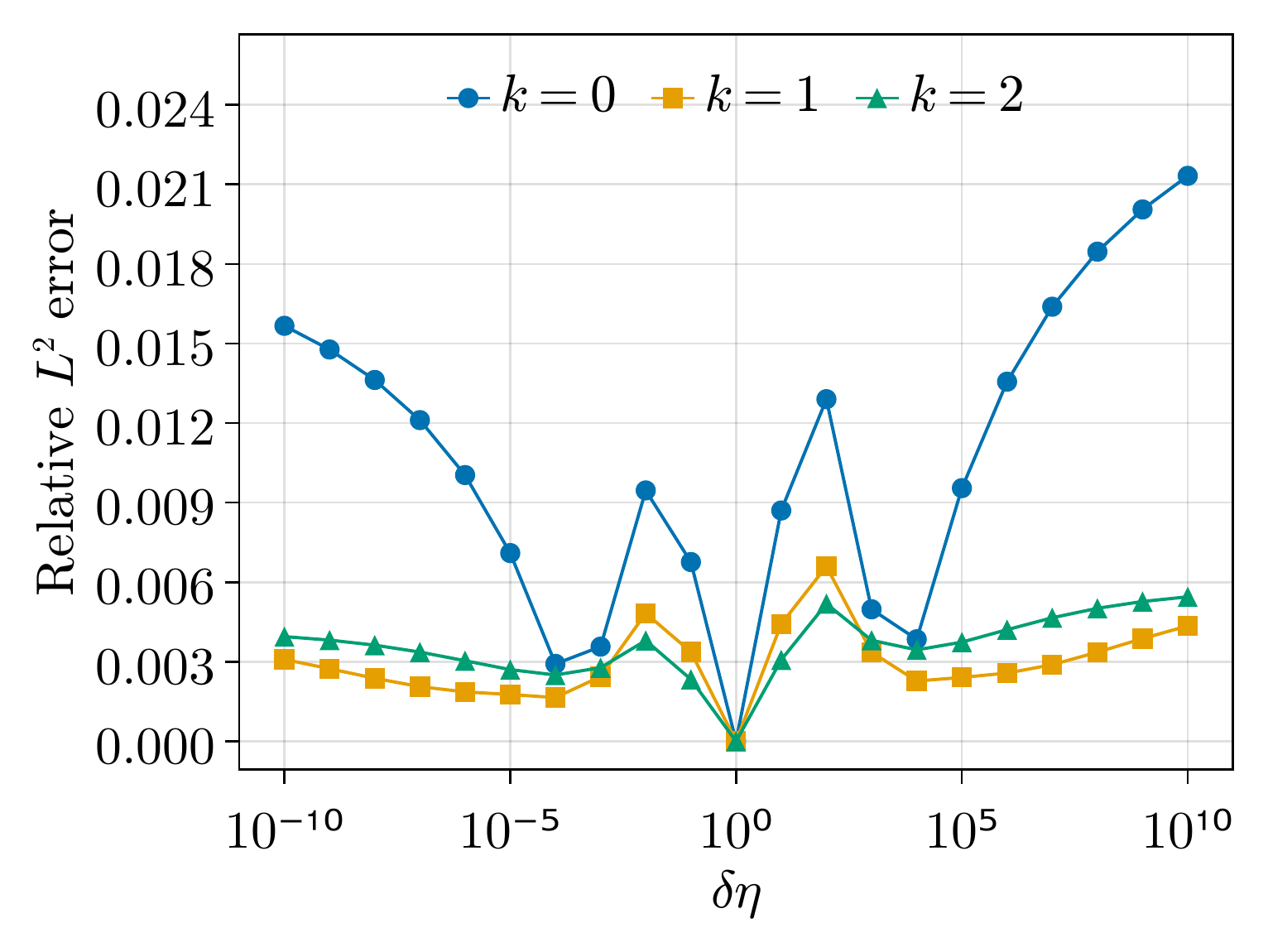}
        \caption{Logarithmic scaling}
    \end{subfigure}
    \caption{$L^2$ errors for the two-strip test case depending on the jump $\delta\eta$.}
    \label{fig:two_strip_jump_to_error}
\end{figure}

The purpose of this test case is to examine if smoothing and the logarithmic scaling of the diffusivity improve the numerical results. We produced results with zero, one, and two smoothing cycles. \Cref{fig:two_strip_strong_form_smoothing} shows that the classical method without smoothing or scaling does not converge. Applying only the scaling, we observe that the method starts to converge for point clouds with a small point density but with an increasing number of points, the error stagnates. By smoothing the diffusivity, we observe numerous effects. First, the method has first-order convergence even without the scaling of the diffusivity. A reason for that is that by smoothing the diffusivity, the discontinuity is smeared out and the discrete gradient operates better. Secondly, if we additionally scale the diffusivity, the error is reduced for all methods. We also find that increasing the number of smoothing cycles does not necessarily lead to better results, see \cref{sfig:two_strip_strong_form_smoothing_scaling_1e10}, due to an increasing distortion of the model with an increasing number of smoothing cycles. Finally, we observe that for some point clouds, the non-smoothed diffusivity leads to better results than the smoothed diffusivity, as can be seen in \cref{sfig:two_strip_strong_form_smoothing_scaling_1e1,sub@sfig:two_strip_strong_form_smoothing_scaling_1e10}. The relative $L^2$ errors depending on the jump magnitude in \cref{fig:two_strip_jump_to_error} show that the methods with a smoothed diffusivity are less susceptible to high jump magnitudes. However, it also shows that the error can be reduced by a factor of $\num{6}$ to $\num{10}$ for the non-smoothed version and by a factor of approximately $\num{4}$ for the smoothed versions by adding logarithmic scaling. The asymmetry of the error is due to the irregularity of the point cloud and because we compute the \emph{relative} error. As seen in \cref{fig:two_strip_solutions}, the $L^2$ norm of the analytical solution is different for different jump magnitudes and thus the relative error differs.

With this information, we establish the \enquote{best} strong form method as the method that includes the scaling of the diffusivity and, depending on the magnitude of the jump, some smoothing of it. From the results, it can be doubted if the smoothing improves the numerical solution. But it stabilizes the convergence and additionally, the non-scaled version tends to have overshoots in the region of the interface. The smoothing causes the overshoot to smear out but not disappear, see \cref{fig:two_strip_solutions}. This behavior is characteristic of strong form methods and is also present in the smoothing method presented by \citet{Li_Ito_2006}. In the following, we use the strong form method with a scaled diffusivity and $\num{2}$ smoothing cycles for comparisons.

\begin{figure}
    \newcommand{\subfigurehspace}{0.49\linewidth}
    \begin{subfigure}{\subfigurehspace}
        \centering
        \includegraphics[width=\linewidth]{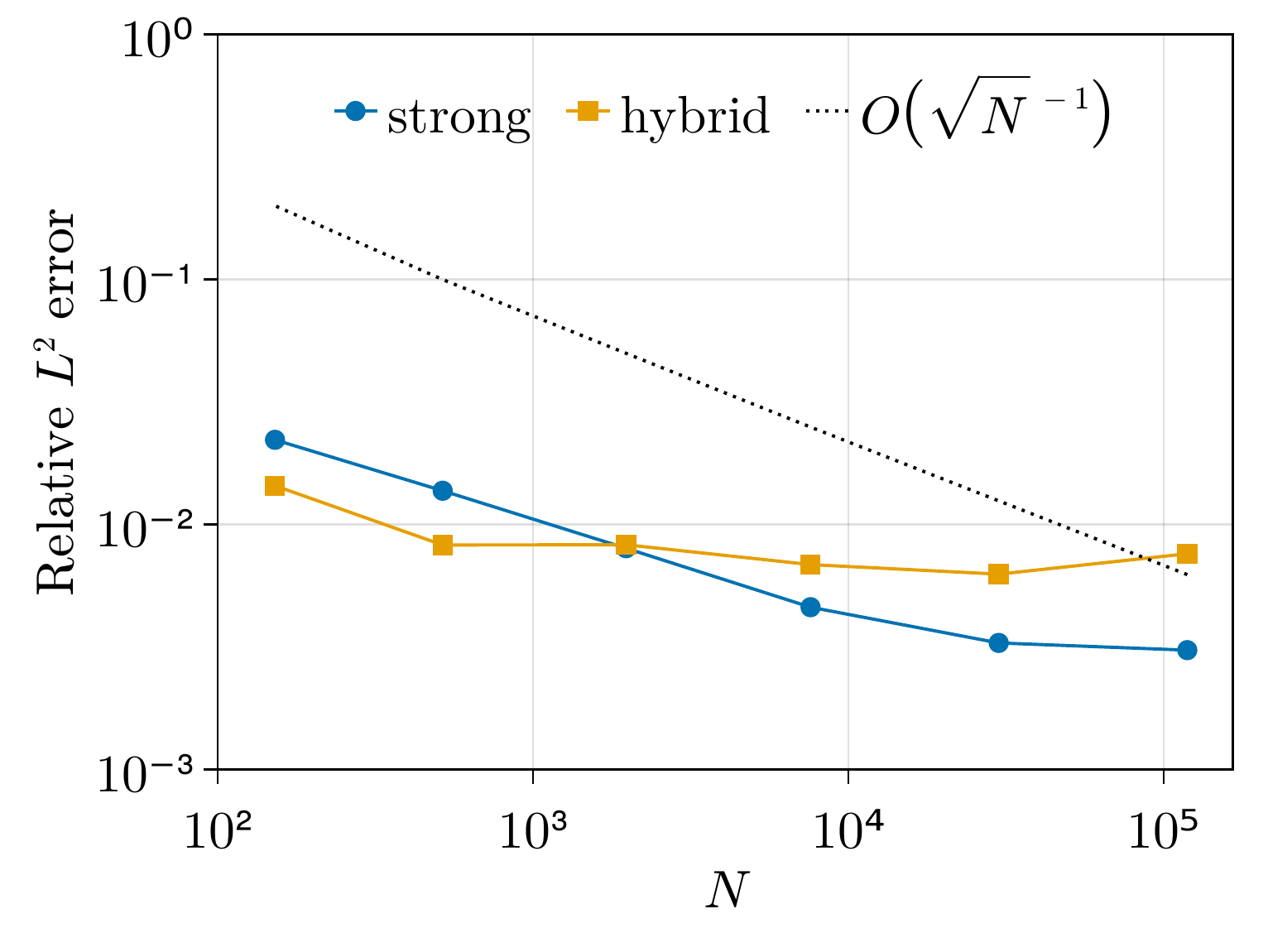}
        \caption{$\delta\eta = \num{1e1}$}
        \label{fig:two_strip_strong_vs_hybrid_low}
    \end{subfigure}
    \begin{subfigure}{\subfigurehspace}
        \centering
        \includegraphics[width=\linewidth]{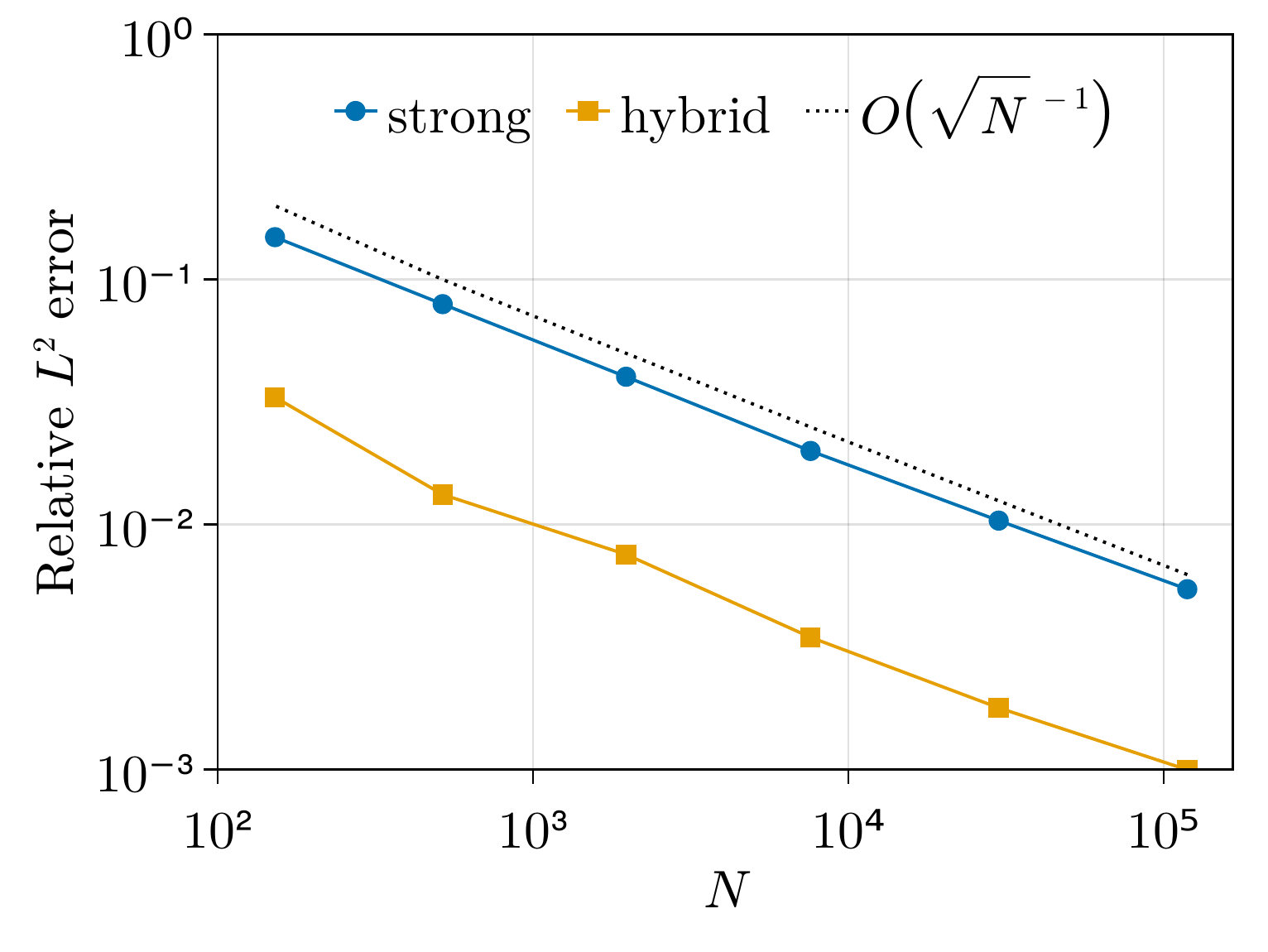}
        \caption{$\delta\eta = \num{1e10}$}
        \label{fig:two_strip_strong_vs_hybrid_high}
    \end{subfigure}
    \caption{$L^2$ errors of strong form method and hybrid method for the two-strip test case with different jump magnitudes.}
    \label{fig:two_strip_strong_vs_hybrid}
\end{figure}

Next, we compare the strong form method with the hybrid method. We distinguish between the hybrid method where we switch to the conservative formulation only for those points where diagonal dominance could not be achieved by the strong form method; and the hybrid method where we switch the formulation for the entire neighborhood of such points. For now, let us only consider the first hybrid method; the second method will be introduced in later examples. The numerical solution of the hybrid method can be seen in \cref{fig:two_strip_solutions}. \Cref{fig:two_strip_strong_vs_hybrid} illustrates the difference between the strong form and the hybrid method for a small jump of only $\delta\eta = \num{1e1}$ (\cref{fig:two_strip_strong_vs_hybrid_low}) and a high jump of $\delta\eta = \num{1e10}$ (\cref{fig:two_strip_strong_vs_hybrid_high}). In \cref{fig:two_strip_strong_vs_hybrid_low} we can clearly see that the hybrid method does not show any convergence. The method gives almost identical results to the scaled strong form method without smoothing as in \cref{sfig:two_strip_strong_form_smoothing_scaling_1e1}. This shows that for small jumps, the classical formulation already yields diagonally dominant rows, and only a few points that have distorted neighborhoods will use the conservative scheme. In fact, we obtain the same results if we compare both methods on uniform point clouds as seen in \cref{fig:comparison_pointclouds_uniform} because on uniform point clouds we observed diagonally dominant rows for each point. Increasing the magnitude of the jump, the error of the hybrid method reduces in comparison to the strong form method, see \cref{fig:two_strip_strong_vs_hybrid_high}. This shows that diagonal dominance depends on both, the quality of the point cloud and sufficient regularity of the diffusivity.

\subsection{Curved Interface} \label{ssec:curved_interface}
We now consider a test case with a curved interface \cite{Davydov_Safarpoor_2021}. The interface of the diffusion coefficient
\[ \eta(x,y) = \begin{cases}
    \eta_L, & y > 2x^3, \\
    \eta_R, & y < 2x^3,
\end{cases} \]
is given by the curve
\[ \Gamma = \set{(x,y) \in \Omega \mid y = 2x^3}. \]

\begin{figure}
    \newcommand{\subfigurehspace}{.49\textwidth}
    \begin{subfigure}{\subfigurehspace}
        \centering
        \includegraphics[width=\linewidth]{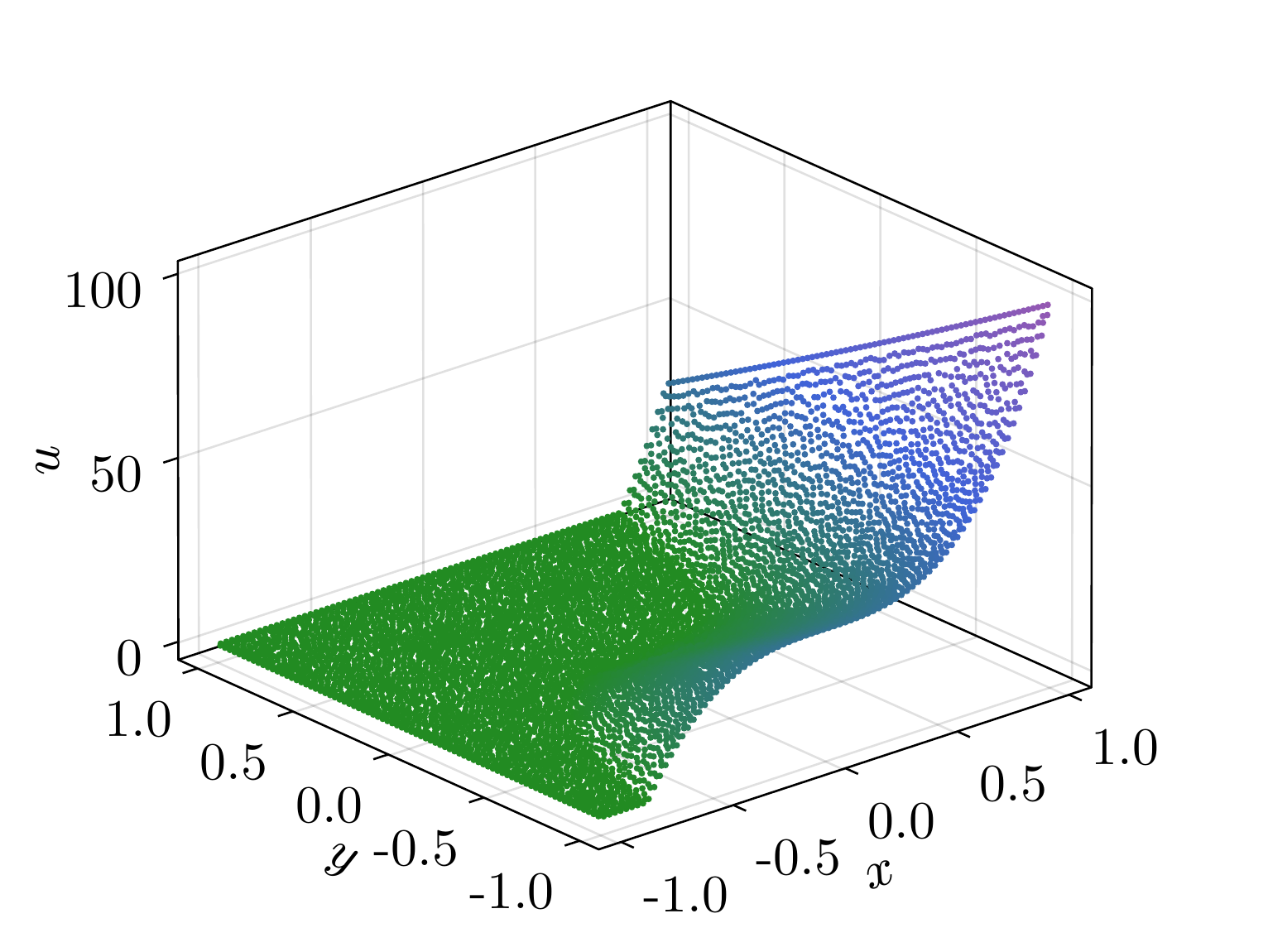}
        \caption{$\delta\eta = \num{1e-8}$}
    \end{subfigure}
    \begin{subfigure}{\subfigurehspace}
        \centering
        \includegraphics[width=\linewidth]{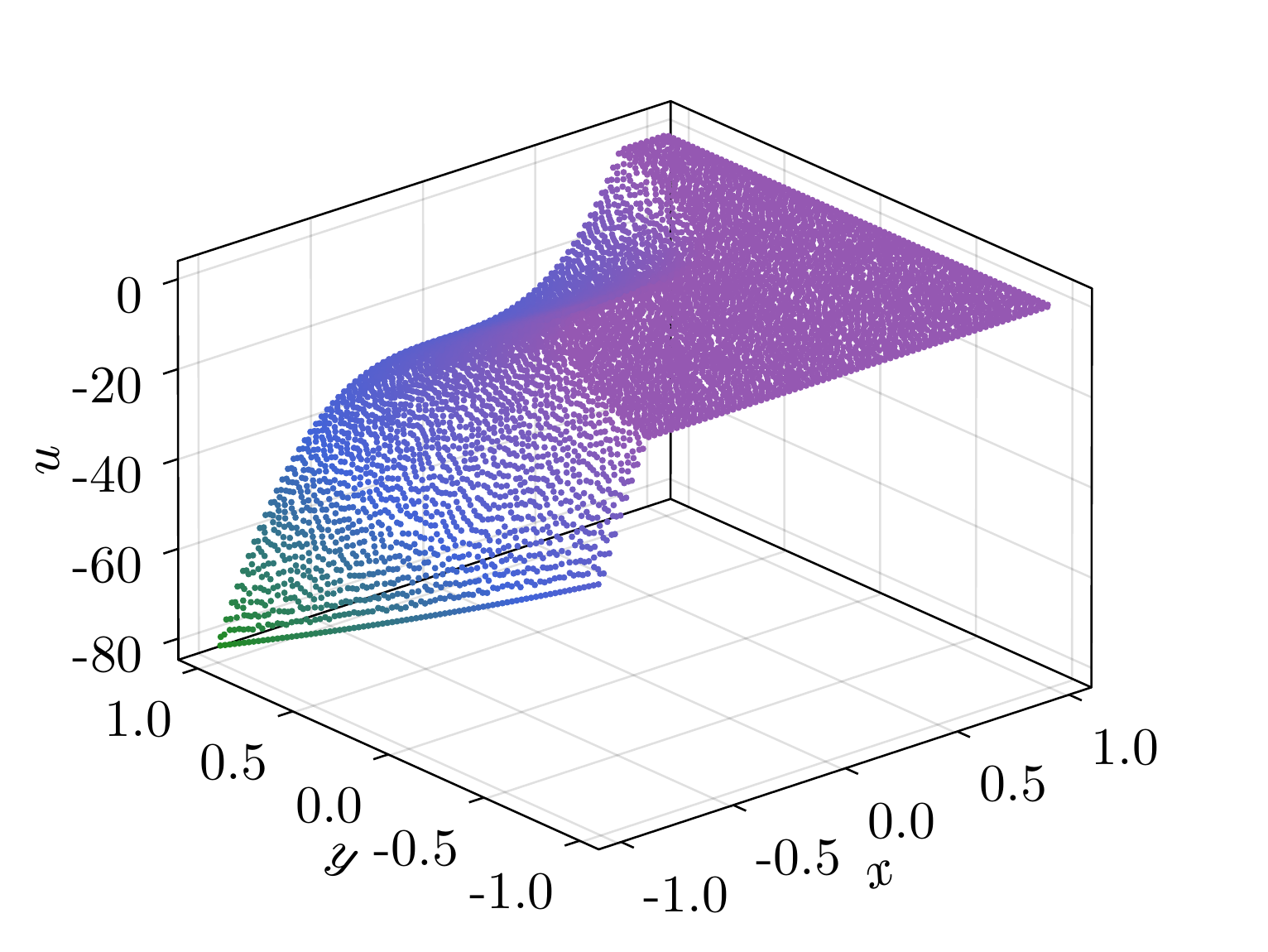}
        \caption{$\delta\eta = \num{1e+8}$}
    \end{subfigure}
    \caption{Analytical solution $u$ of the curved interface problem with different jumps $\delta\eta$.}
    \label{fig:curved_interface_solution}
\end{figure}

For the inhomogeneous Poisson's equation, we define the source term
\[ f(x,y) =  -120x^4 + 24x(y - 15) - 2 \]
and set Dirichlet boundary conditions on the complete boundary such that the analytical solution
\[ u(x,y) = \frac{(y-2x^3)^2 - 30(y-2x^3)}{\eta(x,y)} \]
is obtained.
Unlike the previous test case, the solution in this test case depends not only on the jump $\delta\eta$ but also on the values $\eta_L$ and $\eta_R$. In order to avoid high gradients for $u$, we use
\begin{align*}
    \eta_R &= 10^m
    \intertext{for a fixed $\eta_L = 1$ and}
    \eta_L &= 10^m
\end{align*}
for a fixed $\eta_R = 1$ where $m = 0,\dots,10$ and define the jump magnitude as in the two-strip test case. Examples for analytical solutions can be viewed in \cref{fig:curved_interface_solution} for two different jump magnitudes. In these figures, we observe a weak discontinuity on the curve $\Gamma$.

\begin{figure}
    \newcommand{\subfigurehspace}{.49\linewidth}
    \setlength{\belowcaptionskip}{1.5\baselineskip}
    \begin{subfigure}{\subfigurehspace}
        \centering
        \includegraphics[width=\linewidth]{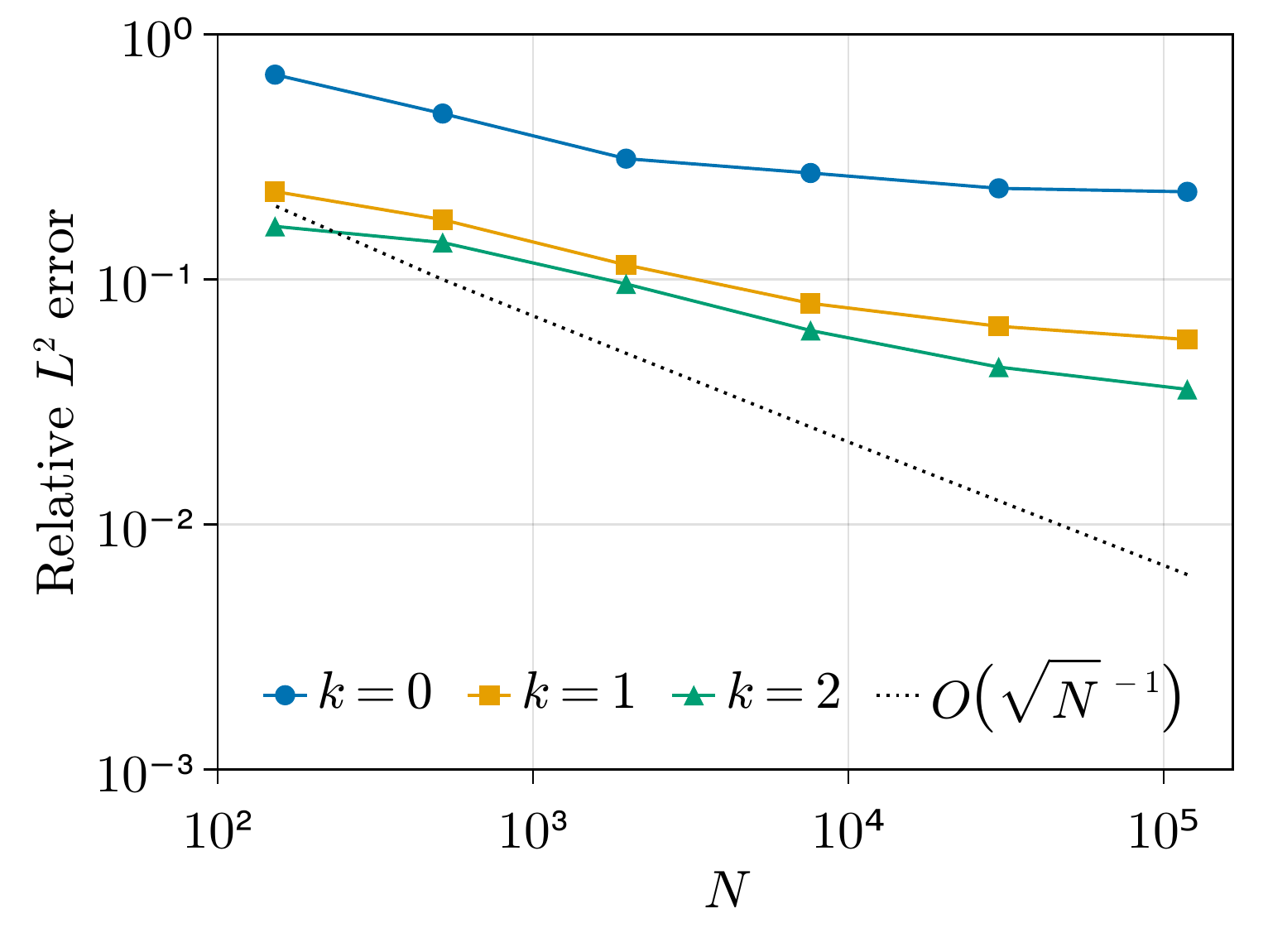}
        \caption{No scaling, $\delta\eta = \num{1e1}$}
        \label{sfig:curved_interface_strong_form_smoothing_noscaling_1e1}
    \end{subfigure}
    \begin{subfigure}{\subfigurehspace}
        \centering
        \includegraphics[width=\linewidth]{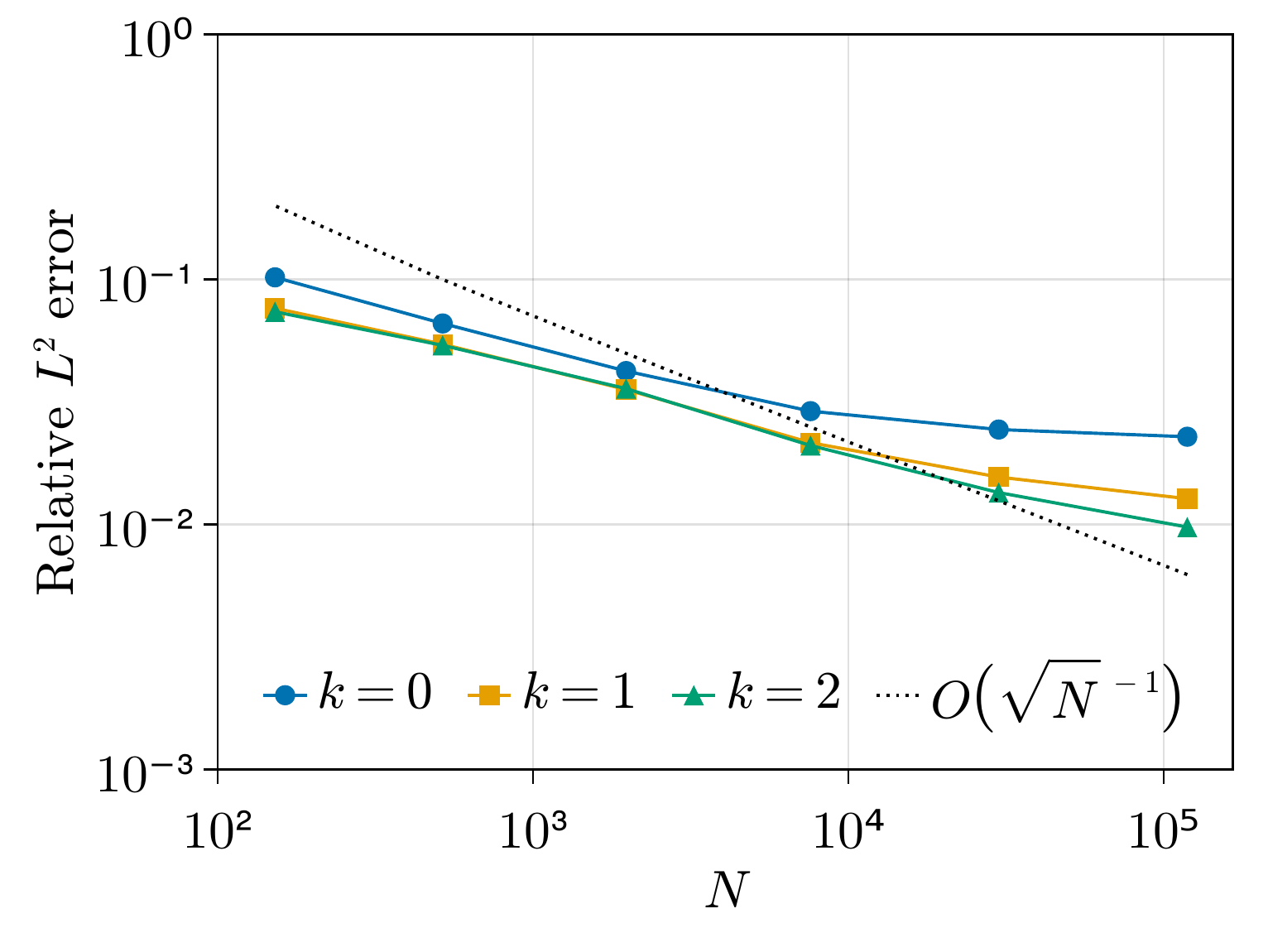}
        \caption{Logarithmic scaling, $\delta\eta = \num{1e1}$}
        \label{sfig:curved_interface_strong_form_smoothing_scaling_1e1}
    \end{subfigure}
    \par
    \begin{subfigure}{\subfigurehspace}
        \centering
        \includegraphics[width=\linewidth]{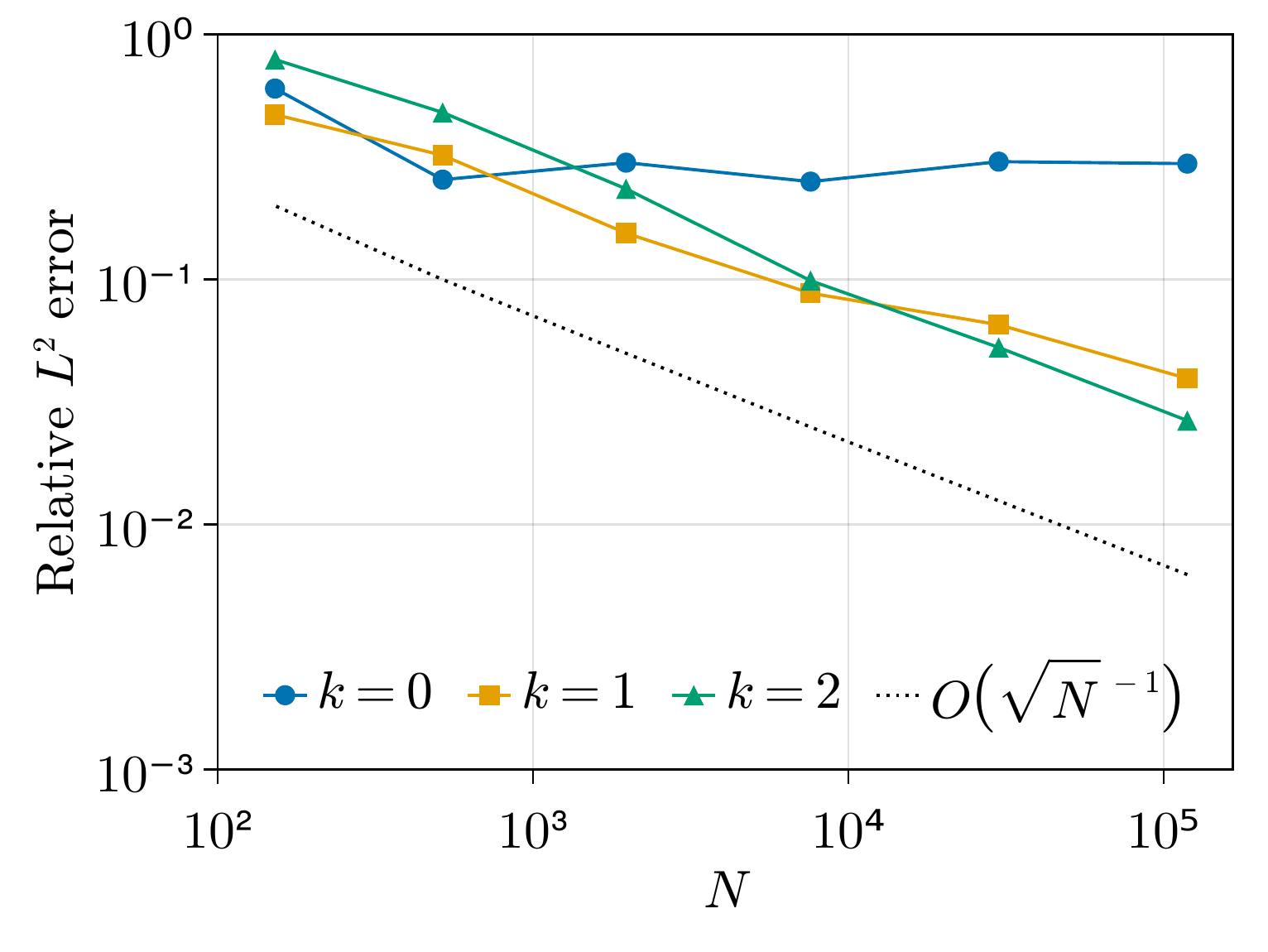}
        \caption{No scaling, $\delta\eta = \num{1e10}$}
        \label{sfig:curved_interface_strong_form_smoothing_noscaling_1e10}
    \end{subfigure}
    \begin{subfigure}{\subfigurehspace}
        \centering
        \includegraphics[width=\linewidth]{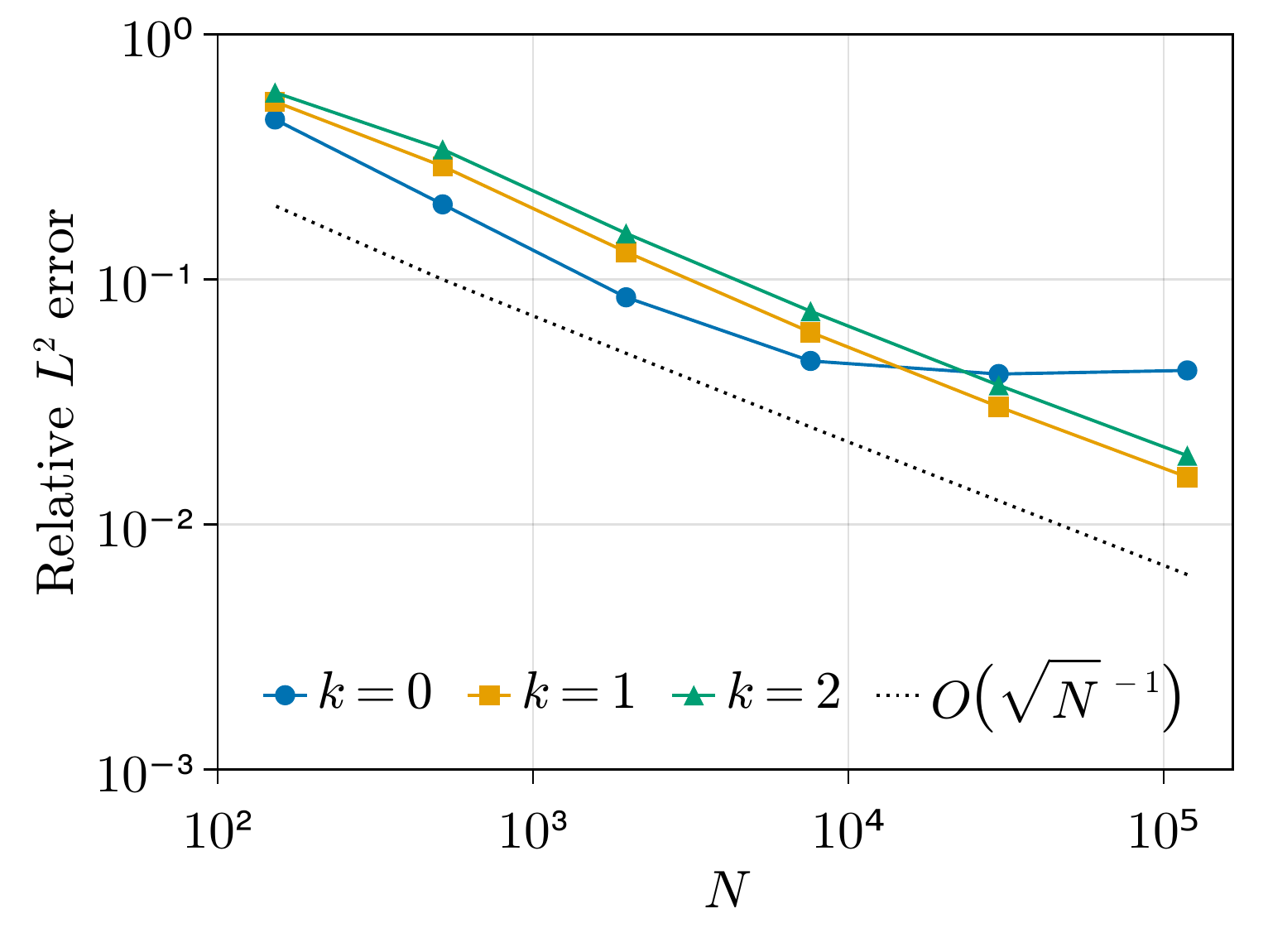}
        \caption{Logarithmic scaling, $\delta\eta = \num{1e10}$}
        \label{sfig:curved_interface_strong_form_smoothing_scaling_1e10}
    \end{subfigure}
    \caption{Convergence plots of the classical strong form method for the curved interface test case with and without logarithmic scaling with different jumps $\delta\eta$ and number of smoothing cycles $k$.}
    \label{fig:curved_interface_strong_form_smoothing}
\end{figure}

\begin{figure}
    \newcommand{\subfigurehspace}{.49\linewidth}
    \begin{subfigure}{\subfigurehspace}
        \centering
        \includegraphics[width=\linewidth]{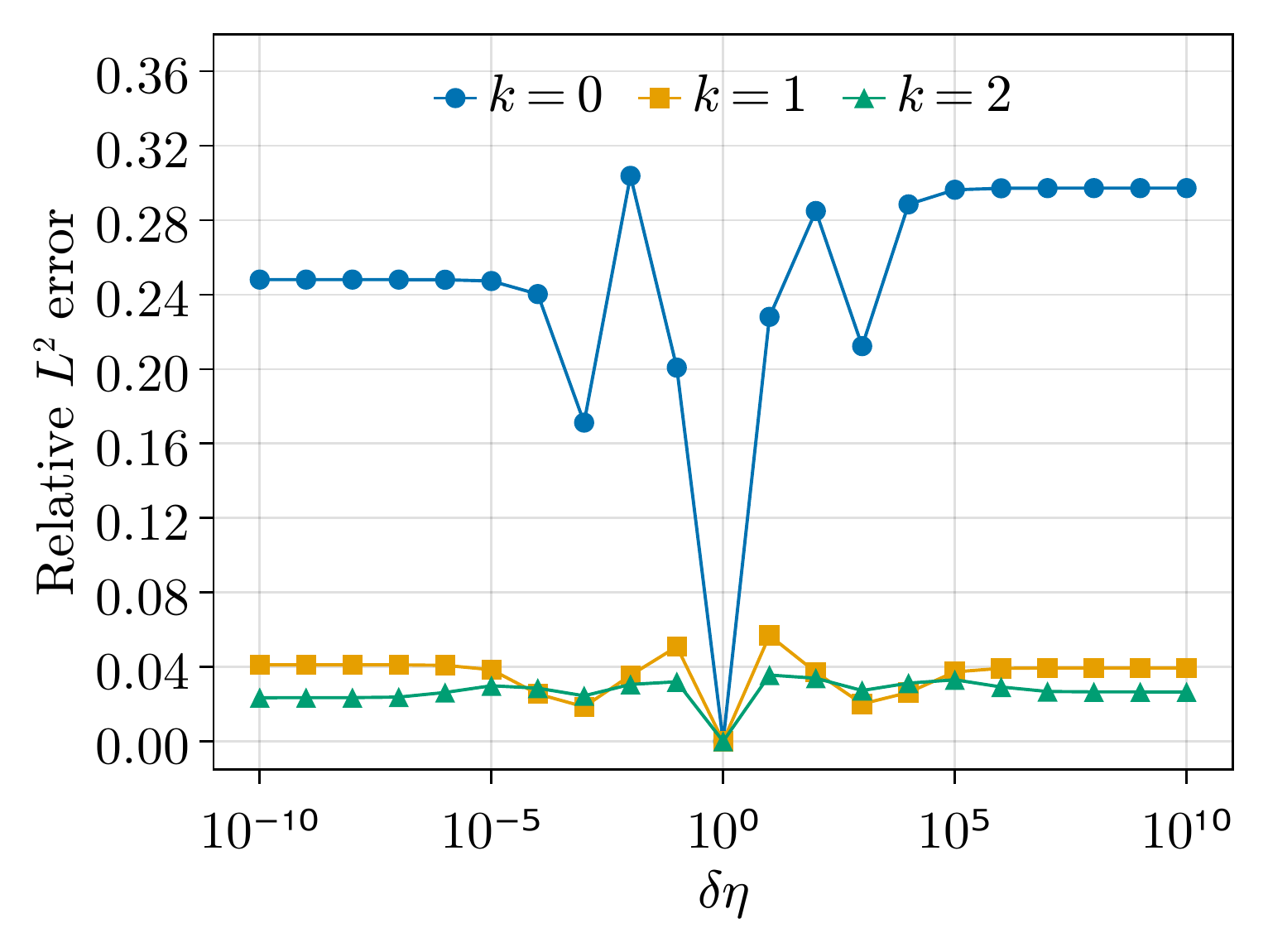}
        \caption{No scaling}
    \end{subfigure}
    \begin{subfigure}{\subfigurehspace}
        \centering
        \includegraphics[width=\linewidth]{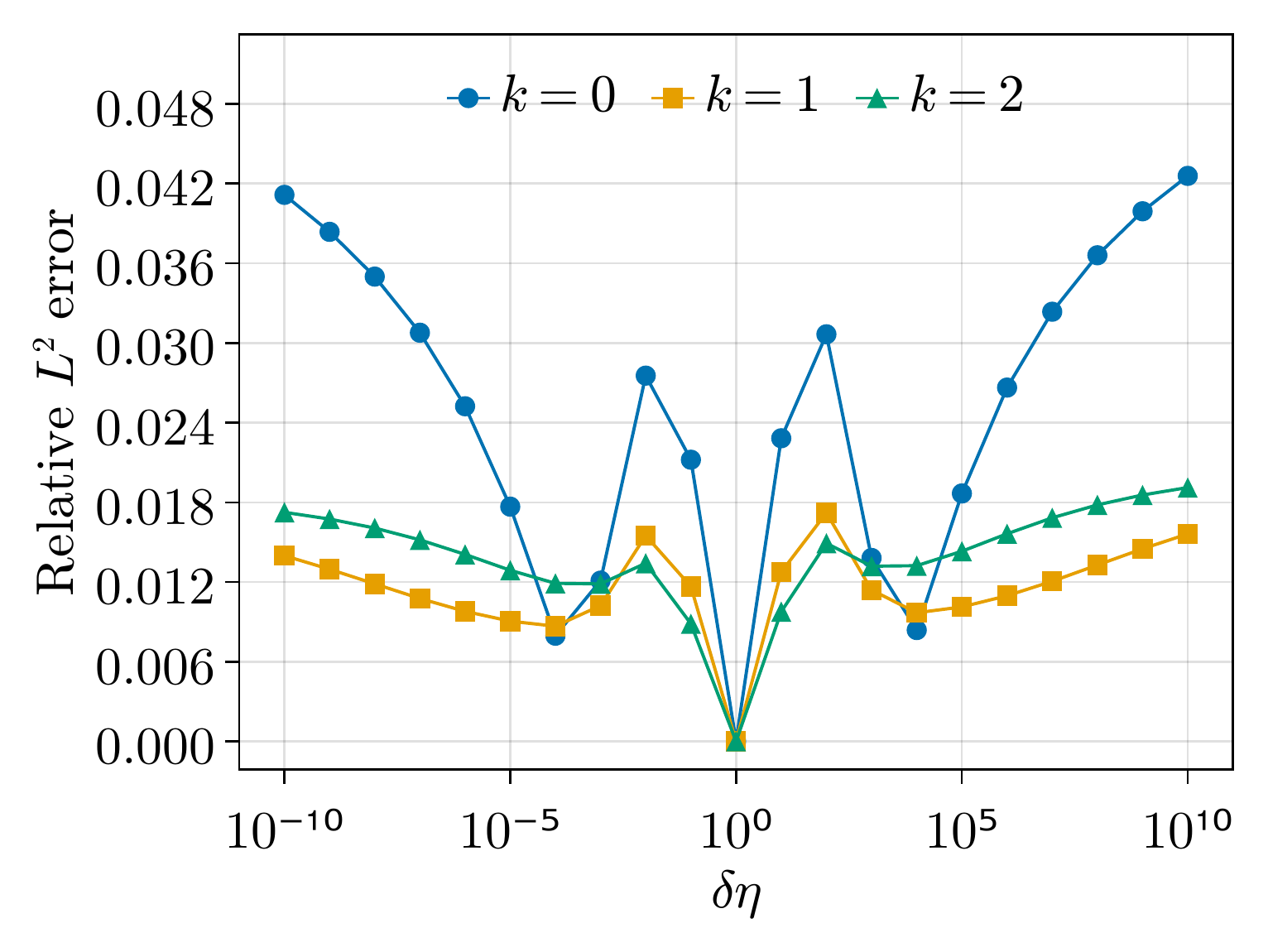}
        \caption{Logarithmic scaling}
    \end{subfigure}
    \caption{$L^2$ errors for the curved interface test case depending on the jump $\delta\eta$.}
    \label{fig:curved_interface_jump_to_error}
\end{figure}

\begin{figure}
    \newcommand{\subfigurehspace}{0.49\linewidth}
    \begin{subfigure}{\subfigurehspace}
        \centering
        \includegraphics[width=\linewidth]{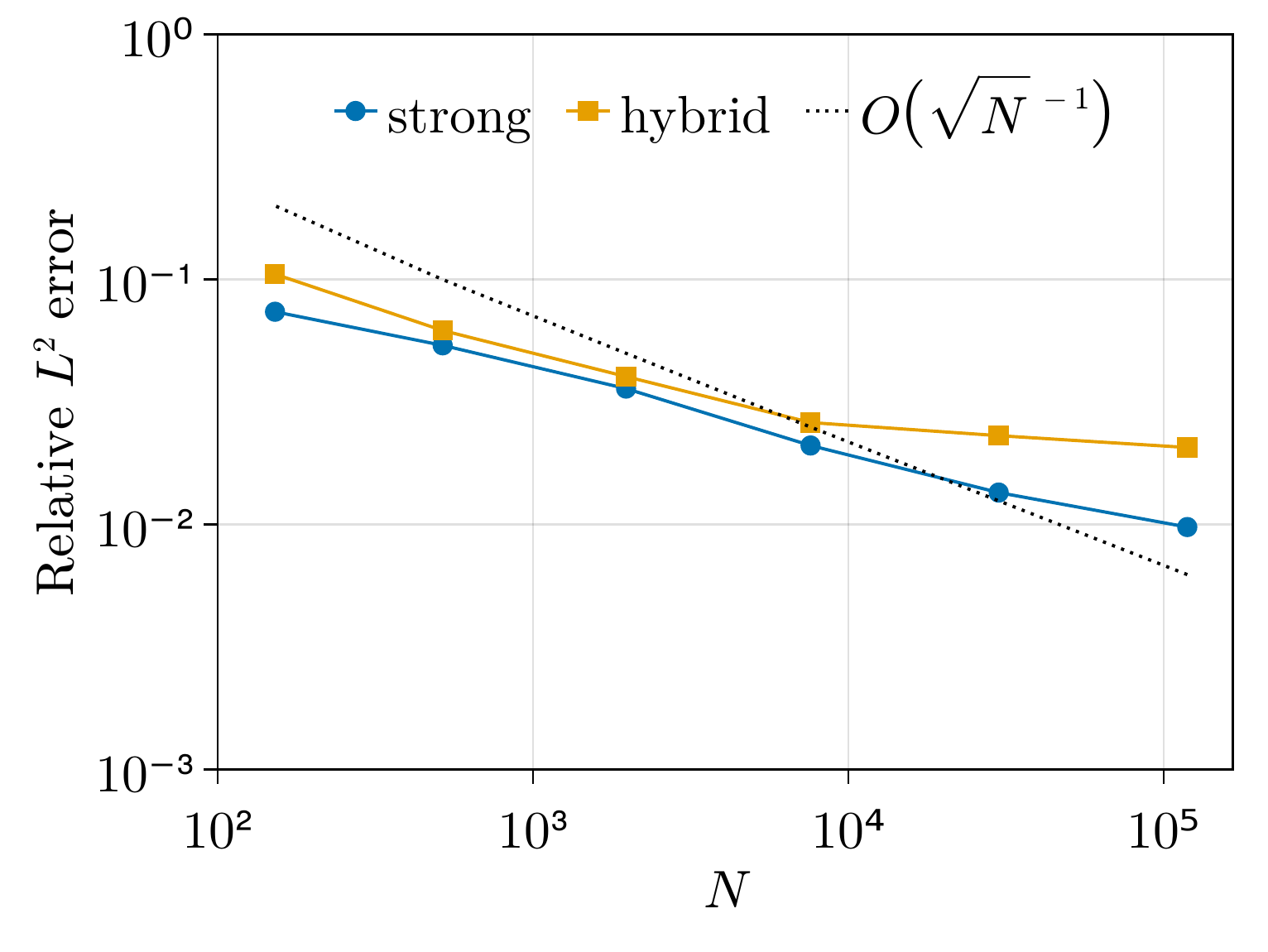}
        \caption{$\delta\eta = \num{1e1}$}
        \label{fig:curved_interface_strong_vs_hybrid_low}
    \end{subfigure}
    \begin{subfigure}{\subfigurehspace}
        \centering
        \includegraphics[width=\linewidth]{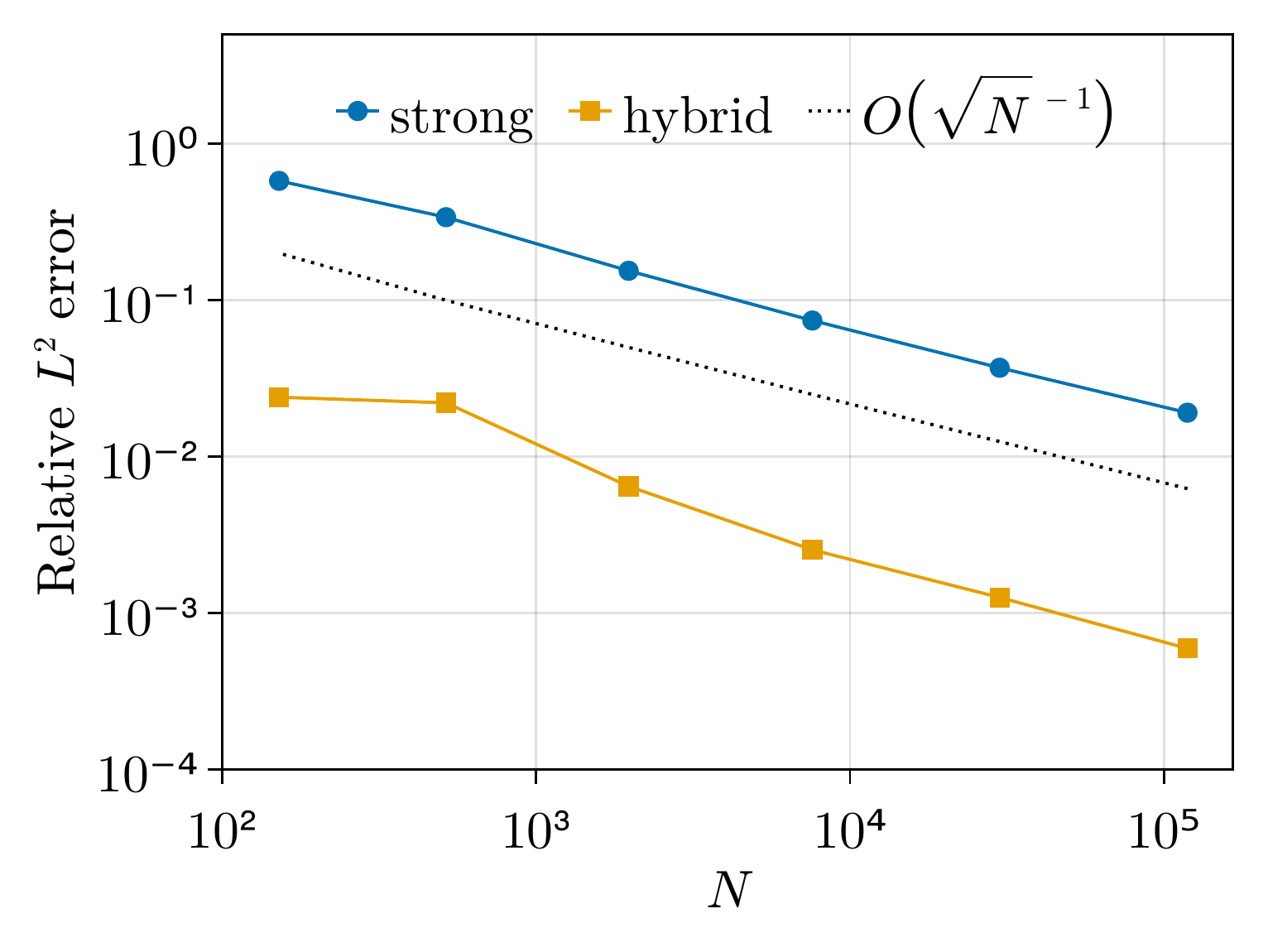}
        \caption{$\delta\eta = \num{1e10}$}
        \label{fig:curved_interface_strong_vs_hybrid_high}
    \end{subfigure}
    \caption{$L^2$ errors of strong form method and hybrid method for the curved interface test case with different jump magnitudes.}
    \label{fig:curved_interface_strong_vs_hybrid}
\end{figure}

For the curved interface test case, we qualitatively observe a similar impact of smoothing and scaling in comparison to the two-strip test case, see \cref{fig:curved_interface_jump_to_error,fig:curved_interface_strong_form_smoothing}. Comparing the errors of both test cases yields that the curved interface test case has higher errors. Similarly, higher errors can be observed in the comparison of the hybrid and the strong form method in \cref{fig:curved_interface_strong_vs_hybrid}. In summary, it can be said that the curvature of the interface has no significant impact on the numerical methods and we draw the same conclusions as for the two-strip test case.

\subsection{Interior Interface} \label{ssec:interior_interface}
Unlike the above two test cases, in this test case, we consider a diffusivity that has a discontinuity along a curve that does not intersect the boundary of the domain. For this test case, we define the function
\[ \phi(x, y) = \cos(\pi x / 2) \cos(\pi y / 2)\] and for the inhomogeneous Poisson's equation, we set the right-hand side
\[ f(x, y) = -\Delta\phi(x, y) = \frac{\pi^2}{2} \phi(x ,y). \]
Additionally, the Dirichlet boundary condition
\[ u|_{\partial\Omega} = 0 \]
is set. Given a height $H \in (0, 1)$, we define the diffusivity
\[ \eta(x, y) = \begin{cases}
    \eta_{\text{out}}, & \phi(x, y) < H, \\
    \eta_{\text{in}},  & \phi(x, y) > H
\end{cases} \]
with $\eta_{\text{out}} = 1$ and $\eta_{\text{in}} = 10^k$ for $k = 0,\dots,10$. Thus, the function
\[ u = \frac{\phi - H}{\eta} + H \]
solves the inhomogeneous Poisson's equation with jump conditions \eqref{eq:jump_conditions}. In our case $H = \frac{3}{4}$ is set which leads to a fully interior interface $\Gamma\subset\Omega$ that separates the domain into $\Omega_{\text{in}}$ and $\Omega_{\text{out}}$. This test case has the additional challenge that the interior domain $\Omega_{\text{in}}$ is entirely contained within $\Omega$ and not directly connected to a Dirichlet boundary. This means that the jump conditions have to be propagated across the interface by the discrete diffusion operator. Since in the previous test cases small jumps of only one or two orders of magnitude could not be correctly identified by the diagonal dominance criterion, we focus on larger jumps of several orders of magnitude instead.

We distinguish between three methods. The strong form method is based on the monomial reproducibility enforced with optimization problem \eqref{eq:formal_optimization_vectors}. The positivity preserving hybrid method is the hybrid approach where we switch to the conservative formulation only for points where diagonal dominance could not be achieved by the strong form method. This results in an M-matrix in the linear system \eqref{eq:poisson_discrete}; hence the method is positivity preserving in the sense of \eqref{eq:maximum_principle_discrete}. Finally, the conservative hybrid method switches to the conservative formulation in the entire neighborhoods of such points. This method is locally conservative additionally to positivity preserving due to the row sum condition \eqref{eq:gauss_column_sum}.

\begin{figure}
    \newcommand{\subfigurehspace}{0.49\linewidth}
    \setlength{\belowcaptionskip}{1.5\baselineskip}

    \begin{subfigure}{\subfigurehspace}
        \centering
        \includegraphics[width=\linewidth]{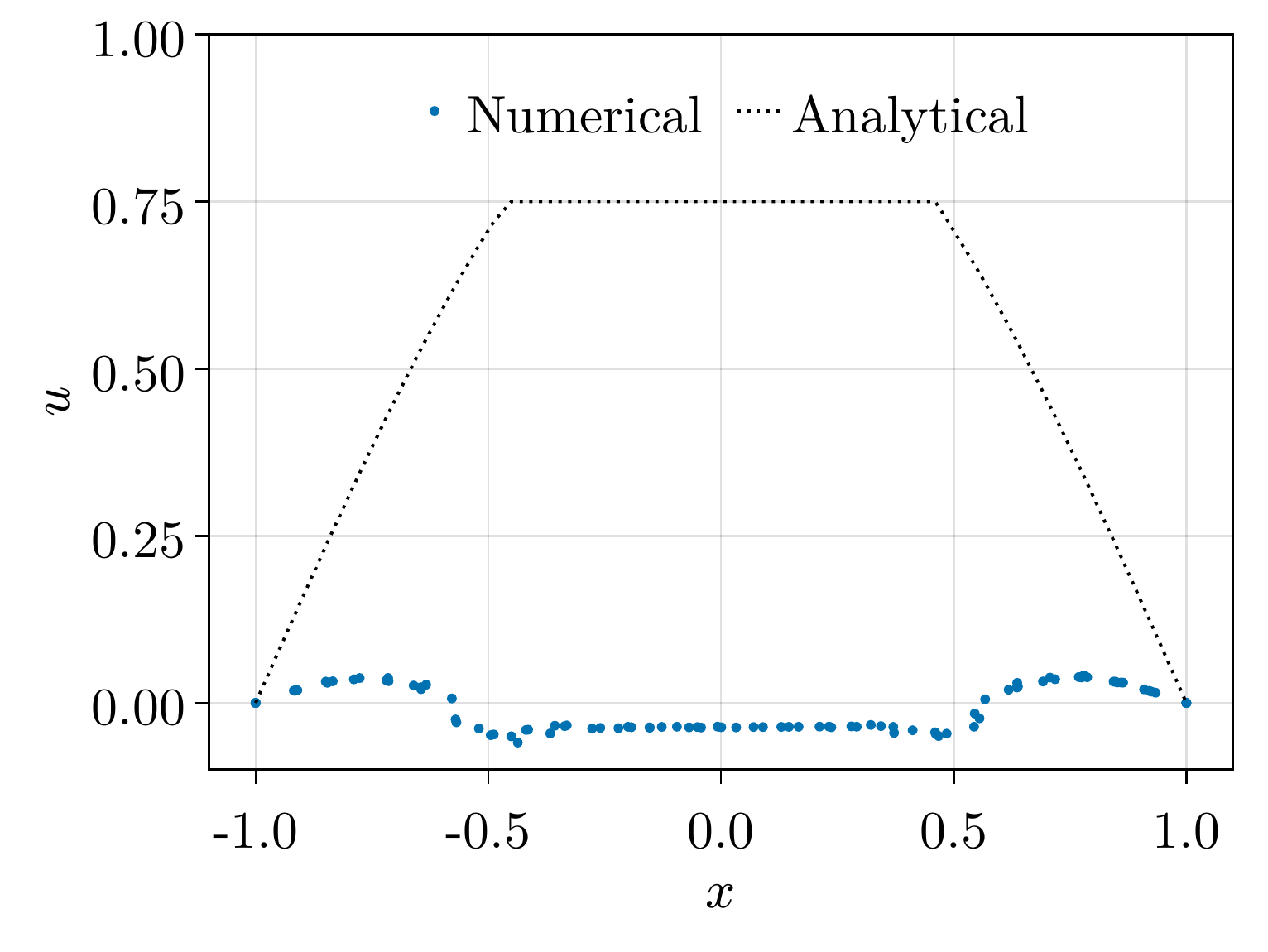}
        \caption{Strong form method without smoothing}
        \label{sfig:interior_interface_solutions_strong_no_smoothing}
    \end{subfigure}
    \begin{subfigure}{\subfigurehspace}
        \centering
        \includegraphics[width=\linewidth]{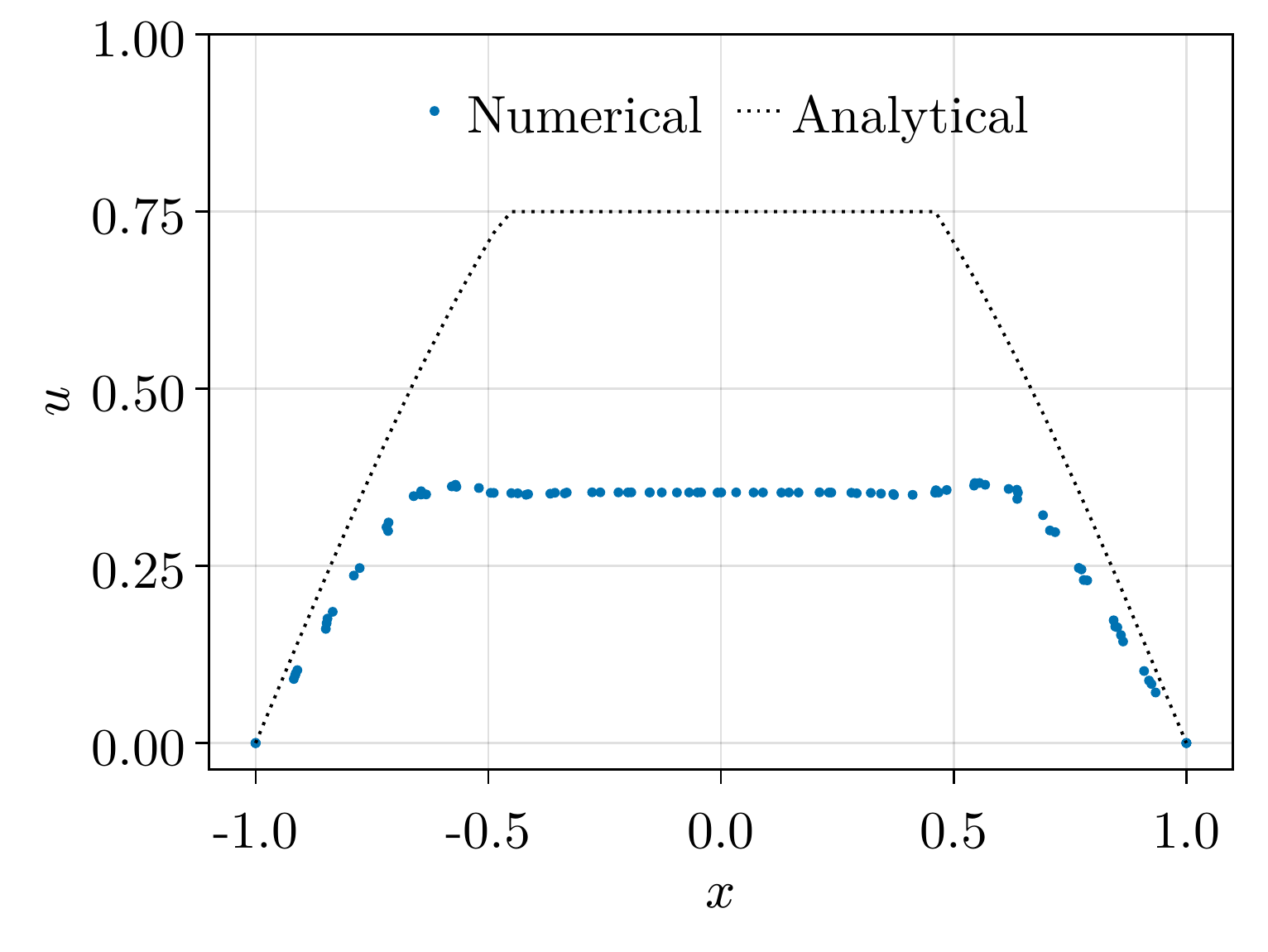}
        \caption{Strong form method with smoothing}
        \label{sfig:interior_interface_solutions_strong_with_smoothing}
    \end{subfigure}
    \par
    \begin{subfigure}{\subfigurehspace}
        \centering
        \includegraphics[width=\linewidth]{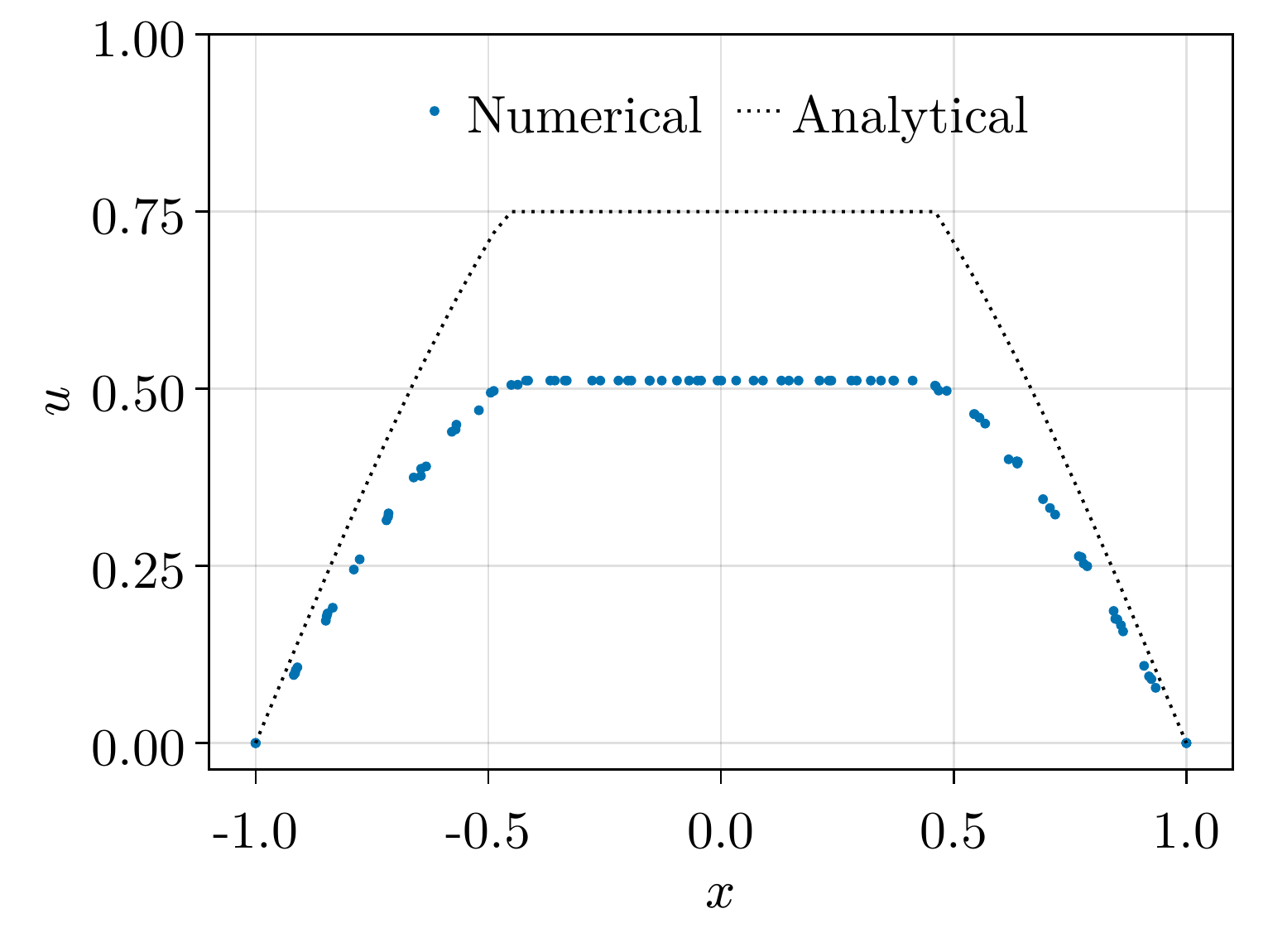}
        \caption{Positivity preserving hybrid method}
        \label{sfig:interior_interface_solutions_positive}
    \end{subfigure}
    \begin{subfigure}{\subfigurehspace}
        \centering
        \includegraphics[width=\linewidth]{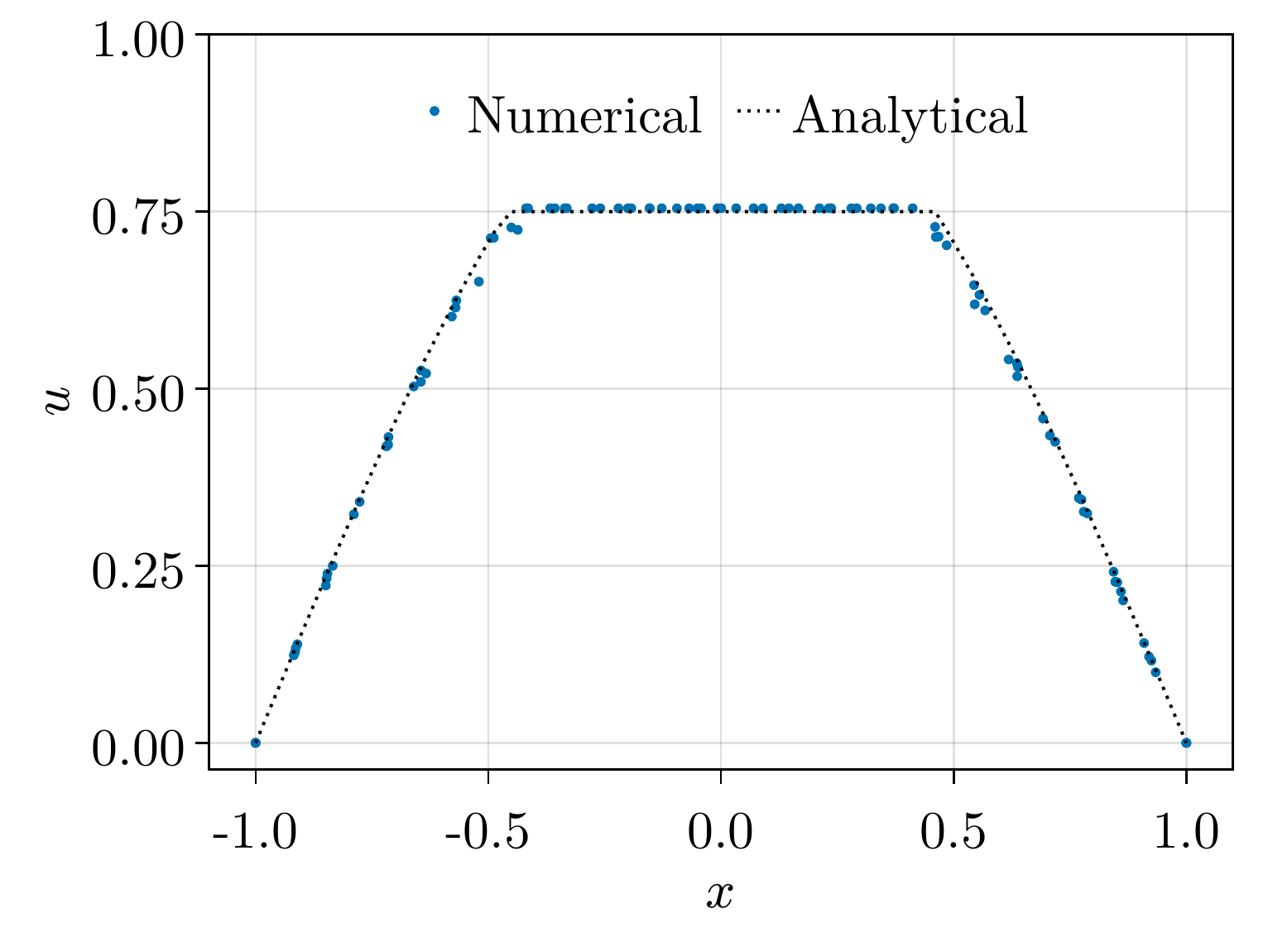}
        \caption{Conservative hybrid method}
        \label{sfig:interior_interface_solutions_conservative}
    \end{subfigure}
    \caption{Solution profile at $y \approx 0$ to the interior interface test case with a jump of $\delta\eta = \num{1e6}$.}
    \label{fig:interior_interface_solutions}
\end{figure}

\begin{figure}
    \setlength{\belowcaptionskip}{1.5\baselineskip}
    \newcommand{\subfigurehspace}{.49\linewidth}
    \begin{subfigure}{\subfigurehspace}
        \centering
        \includegraphics[width=\linewidth]{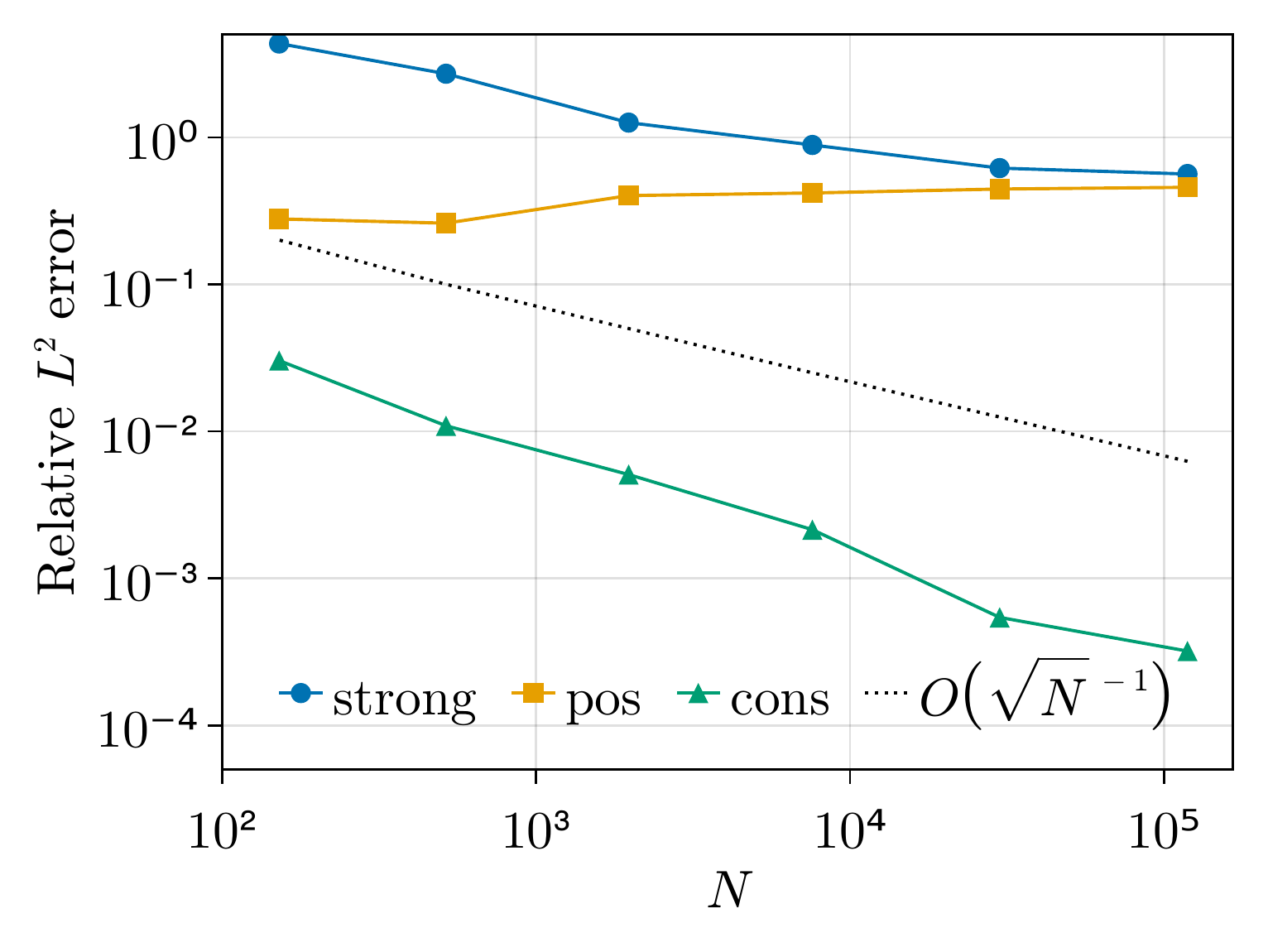}
        \caption{$\delta\eta = \num{1e4}$}
    \end{subfigure}
    \begin{subfigure}{\subfigurehspace}
        \centering
        \includegraphics[width=\linewidth]{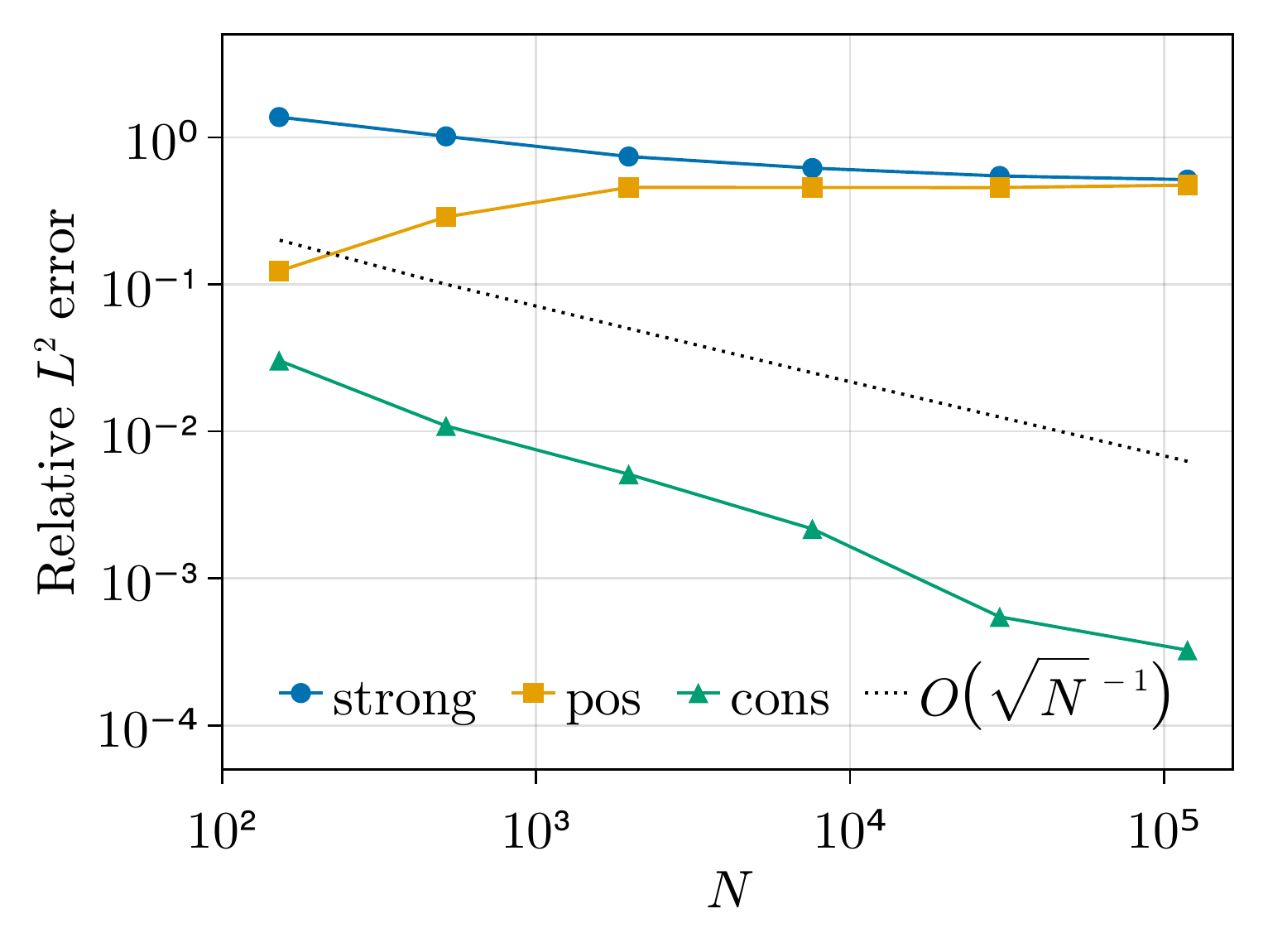}
        \caption{$\delta\eta = \num{1e6}$}
    \end{subfigure}
    \par
    \begin{subfigure}{\subfigurehspace}
        \centering
        \includegraphics[width=\linewidth]{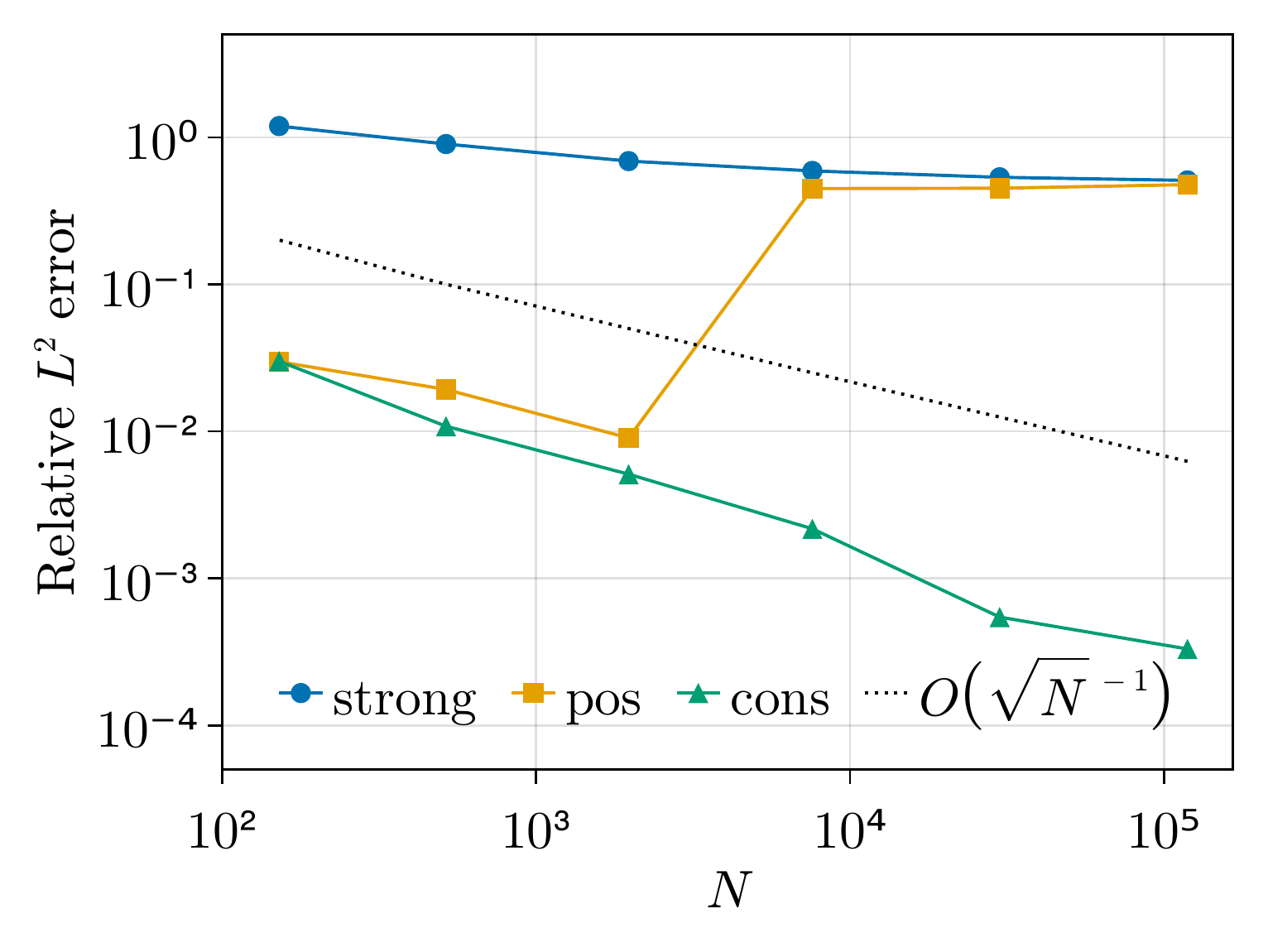}
        \caption{$\delta\eta = \num{1e8}$}
    \end{subfigure}
    \begin{subfigure}{\subfigurehspace}
        \centering
        \includegraphics[width=\linewidth]{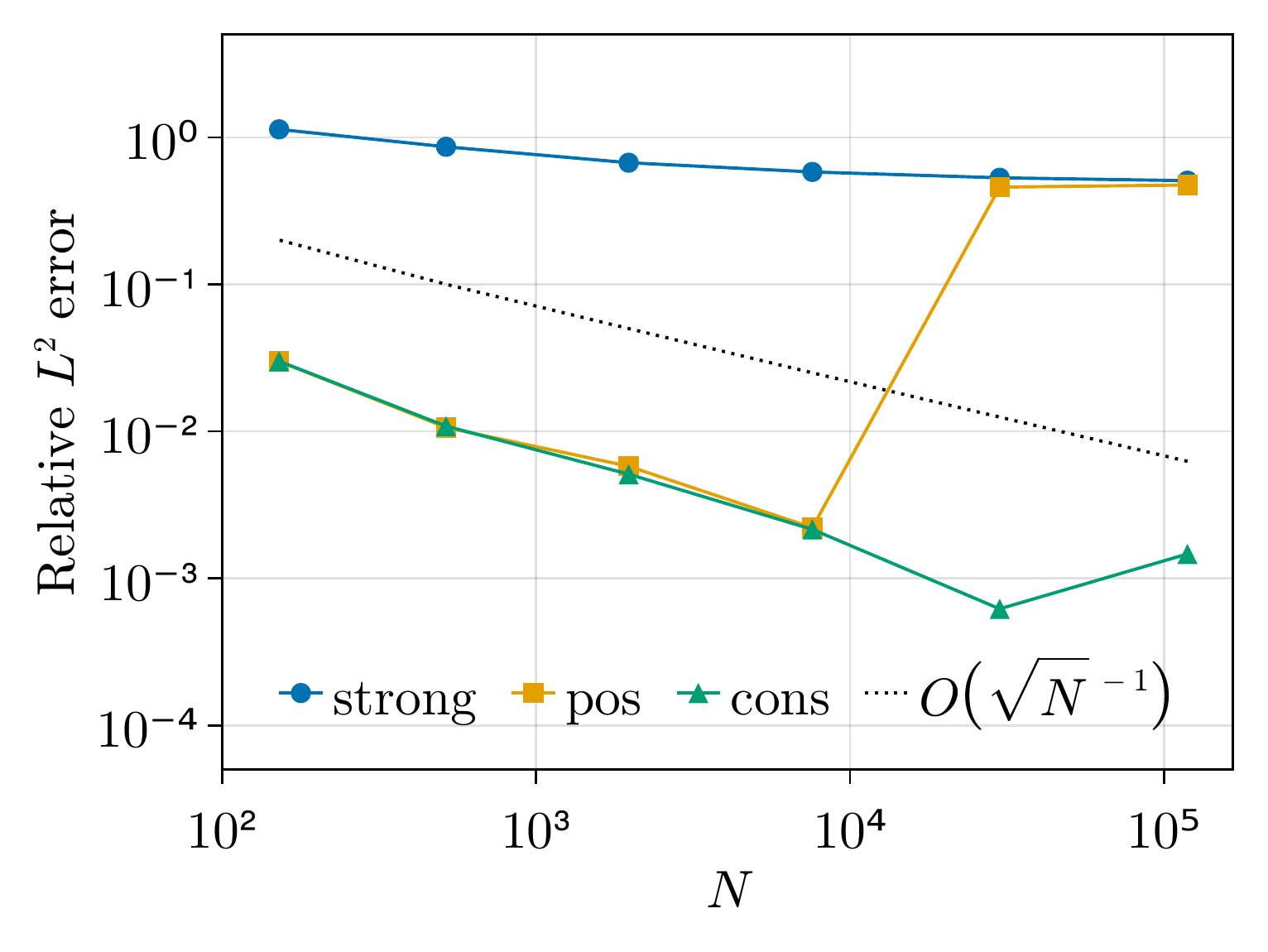}
        \caption{$\delta\eta = \num{1e10}$}
    \end{subfigure}
    \caption{Convergence plots of the classical strong form method, the positivity-preserving method (pos) and the conservative method (cons) for the interior interface test case with different jumps $\delta\eta$.}
    \label{fig:interior_cos_compare_extension}
\end{figure}

\begin{figure}
    \centering
    \includegraphics[width=0.49\linewidth]{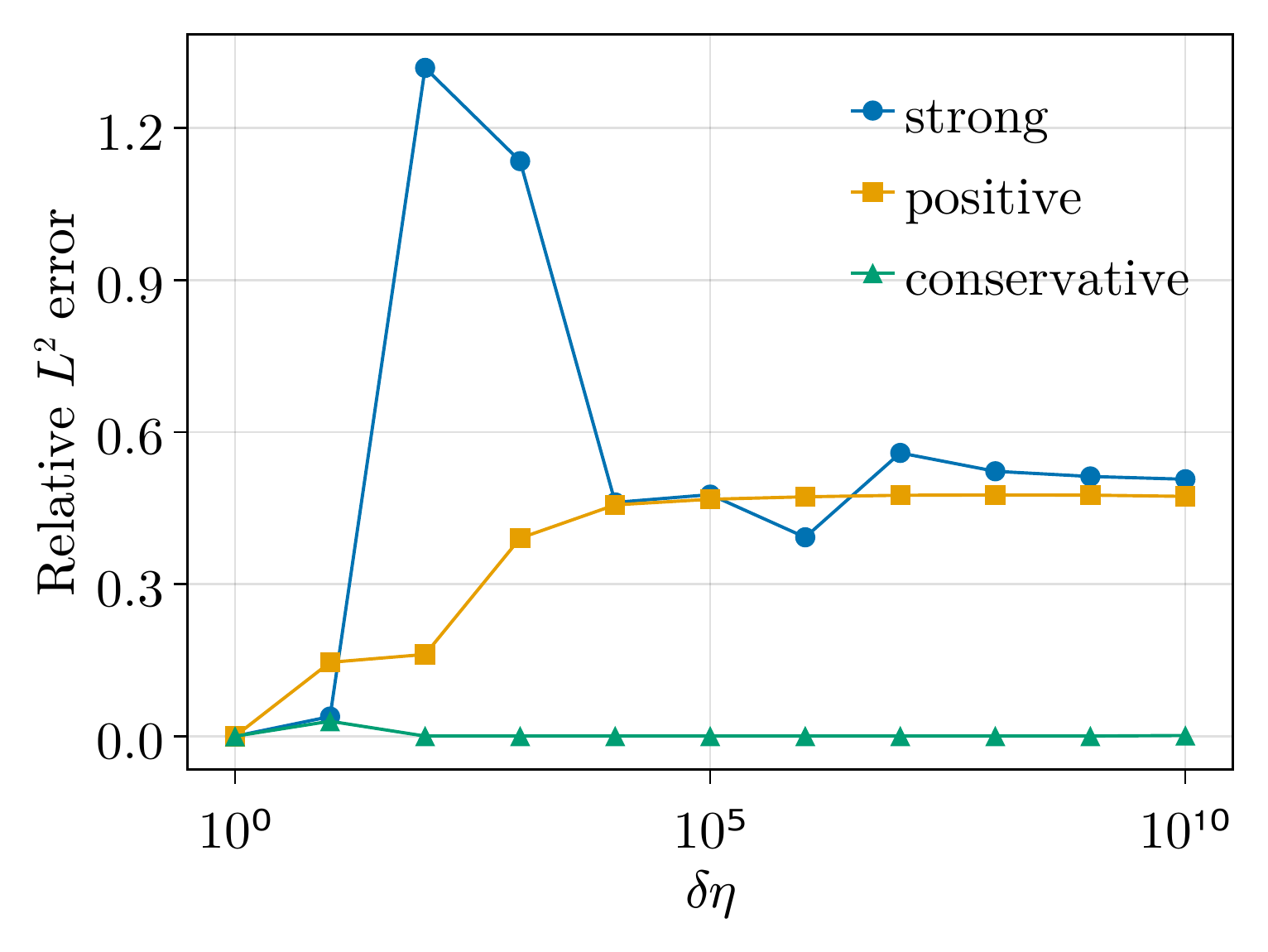}
    \caption{$L^2$ errors for the interior interface test case depending on the jump $\delta\eta$.}
    \label{fig:interior_cos_error_over_jump}
\end{figure}

Note that since $f \ge 0$ on $\Omega$, according to the maximum principle \eqref{eq:maximum_principle} we expect the numerical solution to fulfill $u_i \ge 0$ for each $\bvec{x}_i \in \Omega_h$. \Cref{sfig:interior_interface_solutions_strong_no_smoothing} shows that the strong form method without smoothing fails to provide positive solutions. For the strong form method with smoothing, we observe the characteristic smeared-out overshoots that could already be seen for the two-strip problem, see \cref{sfig:interior_interface_solutions_strong_with_smoothing}.
Switching to the positivity preserving scheme, a positive numerical solution is guaranteed by the method. However, the flux jump condition from \cref{eq:jump_conditions} appears to be violated by the numerical solution shown in \cref{sfig:interior_interface_solutions_positive}. This is caused by points $\bvec{x}_i$ with, for example, $\eta_i = \eta_{\text{in}}$ that have a neighbor $j \in S_i$ across the interface with $\eta_j = \eta_{\text{out}}$ where a diagonally dominant strong form operator could be established and thus flux conservation is not necessarily satisfied. To enforce flux conservation for these points, we need to extend the conservative scheme to these points as well. This results in a good alignment of the numerical and the analytical solution in \cref{sfig:interior_interface_solutions_conservative}. This behavior can be confirmed for different jump magnitudes in \cref{fig:interior_cos_compare_extension} where the positivity preserving method fails to provide first-order convergence. However, the conservative method has first-order convergence and this shows that the strong form scheme produces diagonally dominant discrete operators for some points that have their neighborhoods across the interface, as described earlier. Additionally, the conservative method is insensitive to the jump magnitude as shown in \cref{fig:interior_cos_error_over_jump}. In summary, the conservative scheme shows the best results for high jumps and thus should be the preferred method.

\subsection{Three-Strip Problem} \label{ssec:three_strip}
Similarly to the two-strip test case, we define the diffusivity with three positive values $\eta_L$, $\eta_M$ and $\eta_R$
\[  \eta(x,y) =
\begin{cases}
    \eta_L, & x < -\frac{1}{3}, \\
    \eta_M, & -\frac{1}{3} < x < \frac{1}{3}, \\
    \eta_R, & x > \frac{1}{3}
\end{cases}\]
that leads to three equally sized strips
\begin{align*}
    \Omega_L &= \Set{(x,y) \in \Omega \mid x < - \frac{1}{3}}, \\
    \Omega_M &= \Set{(x,y) \in \Omega \mid - \frac{1}{3} < x < \frac{1}{3}}, \\
    \Omega_R &= \Set{(x,y) \in \Omega \mid x > -\frac{1}{3}}.
\end{align*}
Dirichlet and Neumann boundary conditions are set exactly as for the two-strip test case such that the analytical solution is again monotonously decreasing and piecewise linear in $x$-direction. A three-strip problem can be characterized by the two jumps at the interfaces
\begin{align*}
    \delta\eta_L &= \frac{\eta_M}{\eta_L}, \\
    \delta\eta_R &= \frac{\eta_R}{\eta_M}
\end{align*}
and with $\eta_L = 1$ we have
\begin{align*}
    \eta_M &= \delta\eta_L, \\
    \eta_R &= \delta\eta_L \delta\eta_R.
\end{align*}

In this test case, we added additional complexity to the problem by partitioning the domain into the three strips where the interior strip $\Omega_M$ is not directly connected to a Dirichlet boundary. The only boundary condition that appears for this subdomain is a Neumann boundary condition which introduces additional challenges to the numerical scheme.
\begin{enumerate}
    \item We need the numerical scheme to propagate the jump conditions \eqref{eq:jump_conditions} from the surrounding domains $\Omega_L$ and $\Omega_R$ correctly.
    \item Boundary points with a Neumann boundary condition that lie in the proximity of the interface have to propagate the correct jump conditions across the interface boundary.
\end{enumerate}

\begin{figure}
    \newcommand{\subfigurehspace}{0.49\linewidth}
    \centering
    \includegraphics[width=\subfigurehspace]{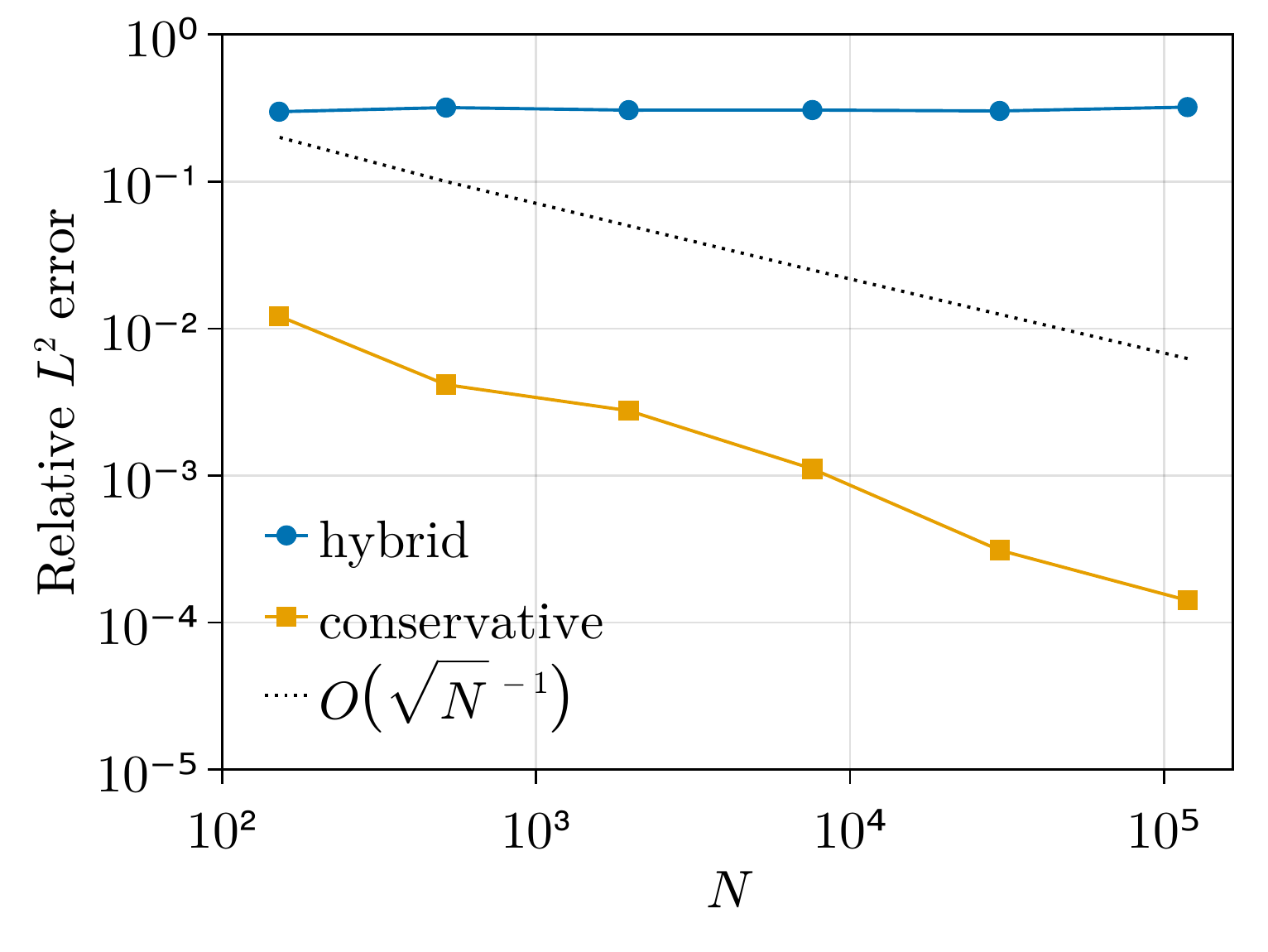}
    \caption{$L^2$ errors for the three-strip test case of hybrid method with and without conservative Neumann boundary conditions for a jump $\delta\eta_L = \num{1e6}$ and $\delta\eta_R = \num{1e-4}$.}
    \label{fig:convergence_three_strip}
\end{figure}

\begin{figure}
    \newcommand{\subfigurehspace}{0.49\linewidth}
    \begin{subfigure}{\subfigurehspace}
        \centering
        \includegraphics[width=\linewidth]{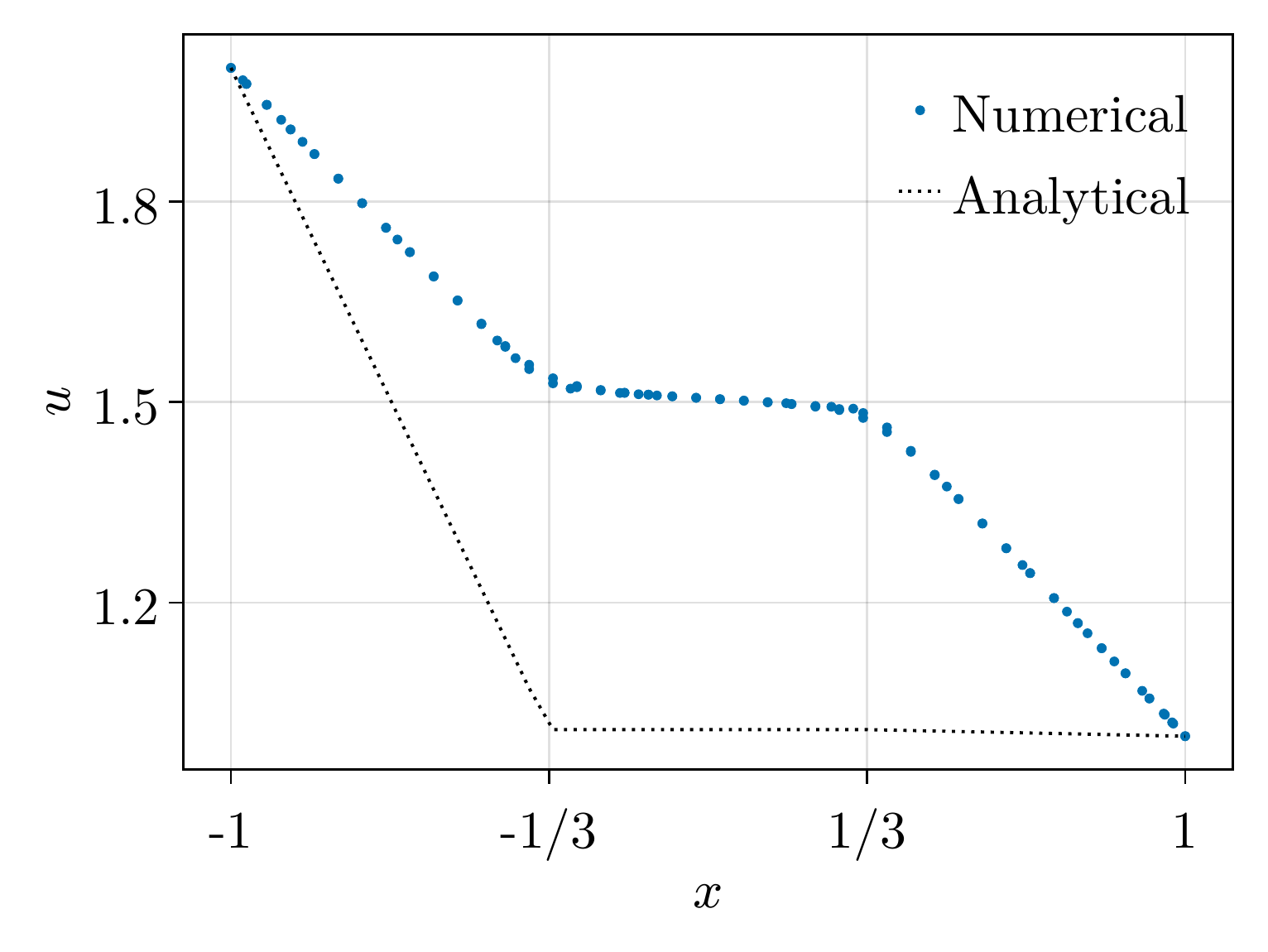}
        \caption{Regular Neumann boundaries}
        \label{fig:solutions_three_strip_regular}
    \end{subfigure}
    \begin{subfigure}{\subfigurehspace}
        \centering
        \includegraphics[width=\linewidth]{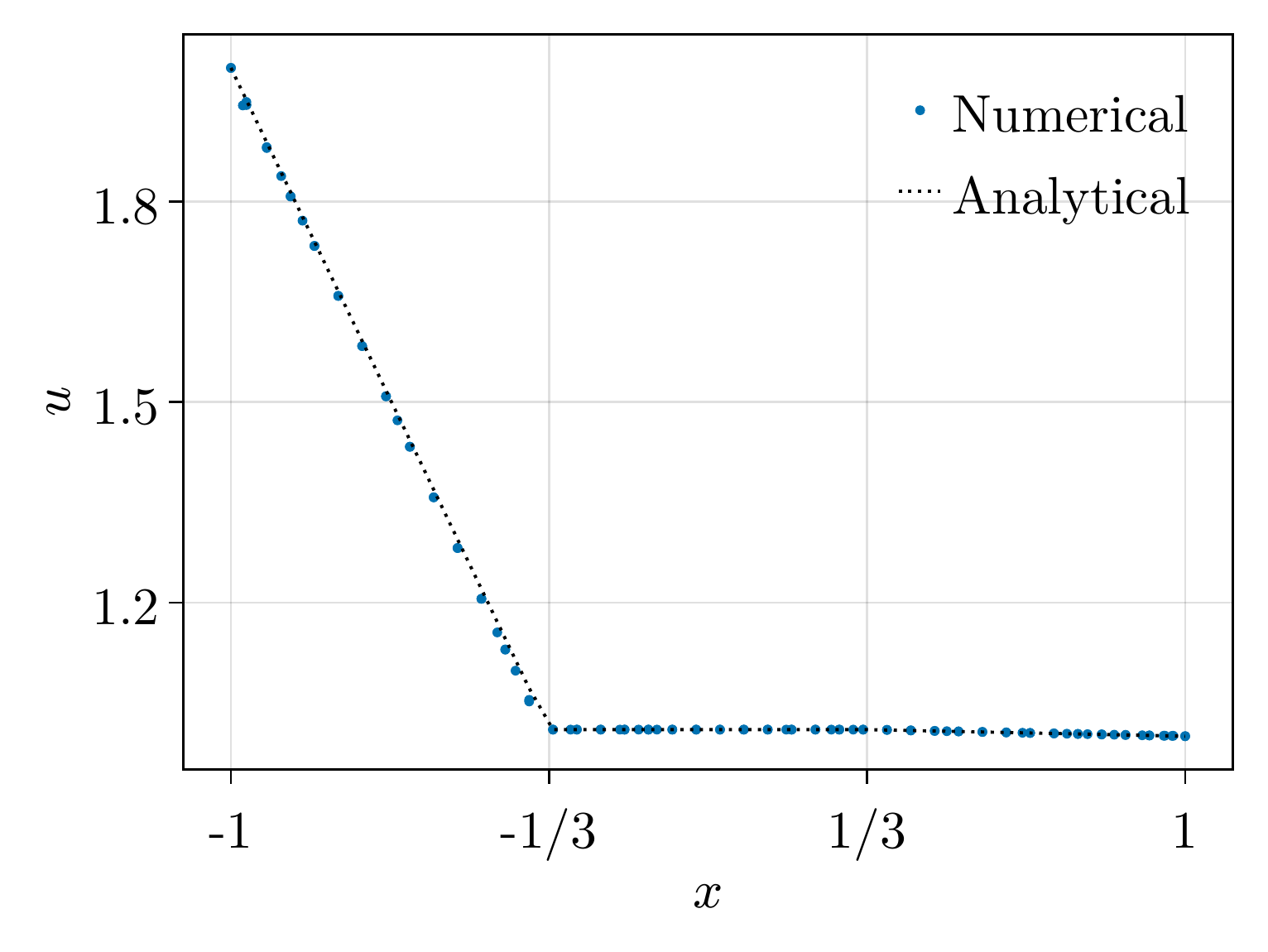}
        \caption{Conservative Neumann boundaries}
        \label{fig:solutions_three_strip_conservative}
    \end{subfigure}
    \caption{Solution profiles comparison at $y \approx 0$ for the three-strip test case of hybrid method with and without conservative Neumann boundary conditions for a jump $\delta\eta_L = \num{1e6}$ and $\delta\eta_R = \num{1e-4}$.}
    \label{fig:solutions_three_strip}
\end{figure}

\begin{figure}
    \newcommand{\subfigurehspace}{0.49\linewidth}
    \begin{subfigure}{\subfigurehspace}
        \centering
        \includegraphics[width=\linewidth]{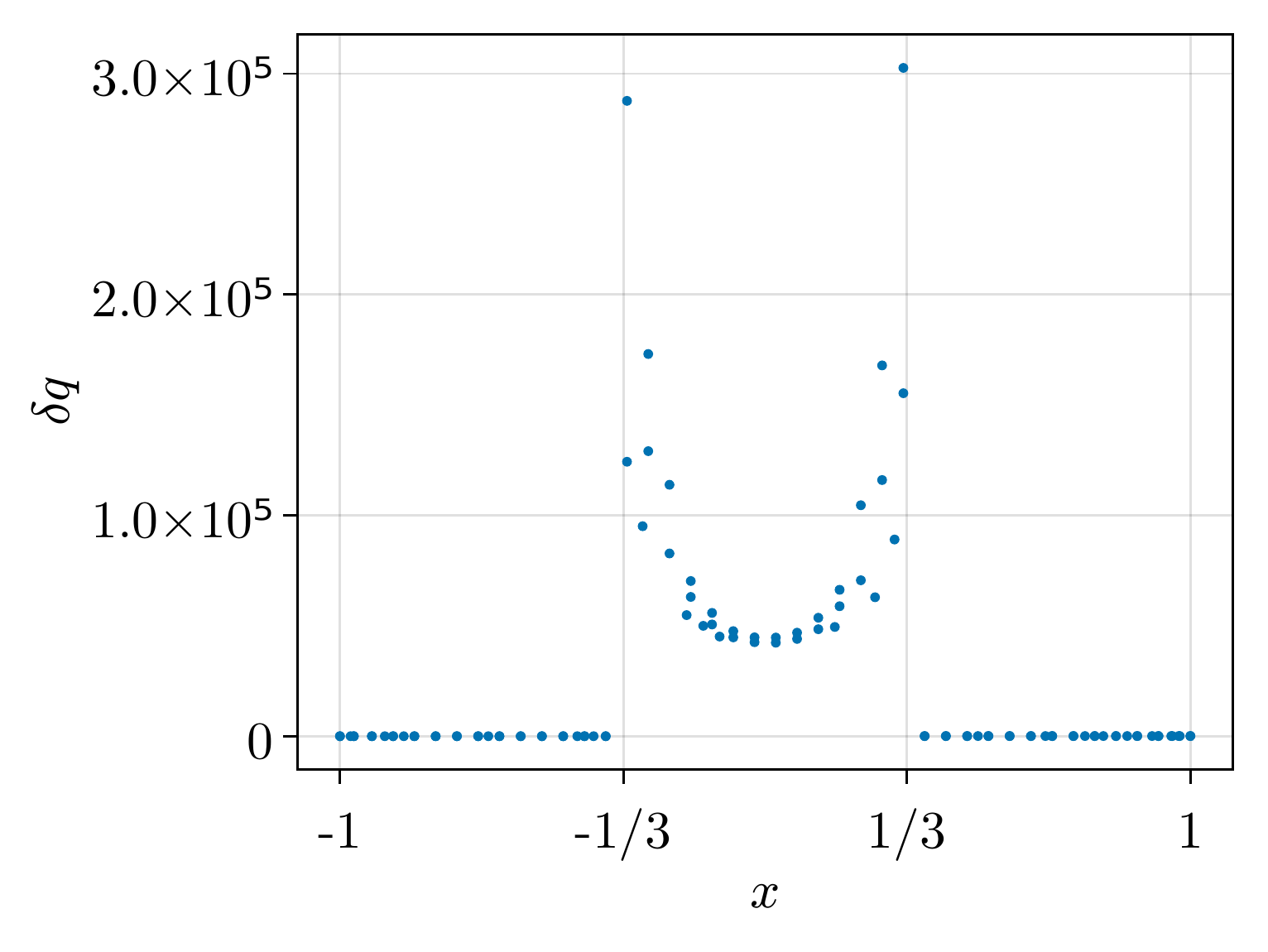}
        \caption{Regular Neumann boundaries}
        \label{fig:etagrad_three_strip_regular}
    \end{subfigure}
    \begin{subfigure}{\subfigurehspace}
        \centering
        \includegraphics[width=\linewidth]{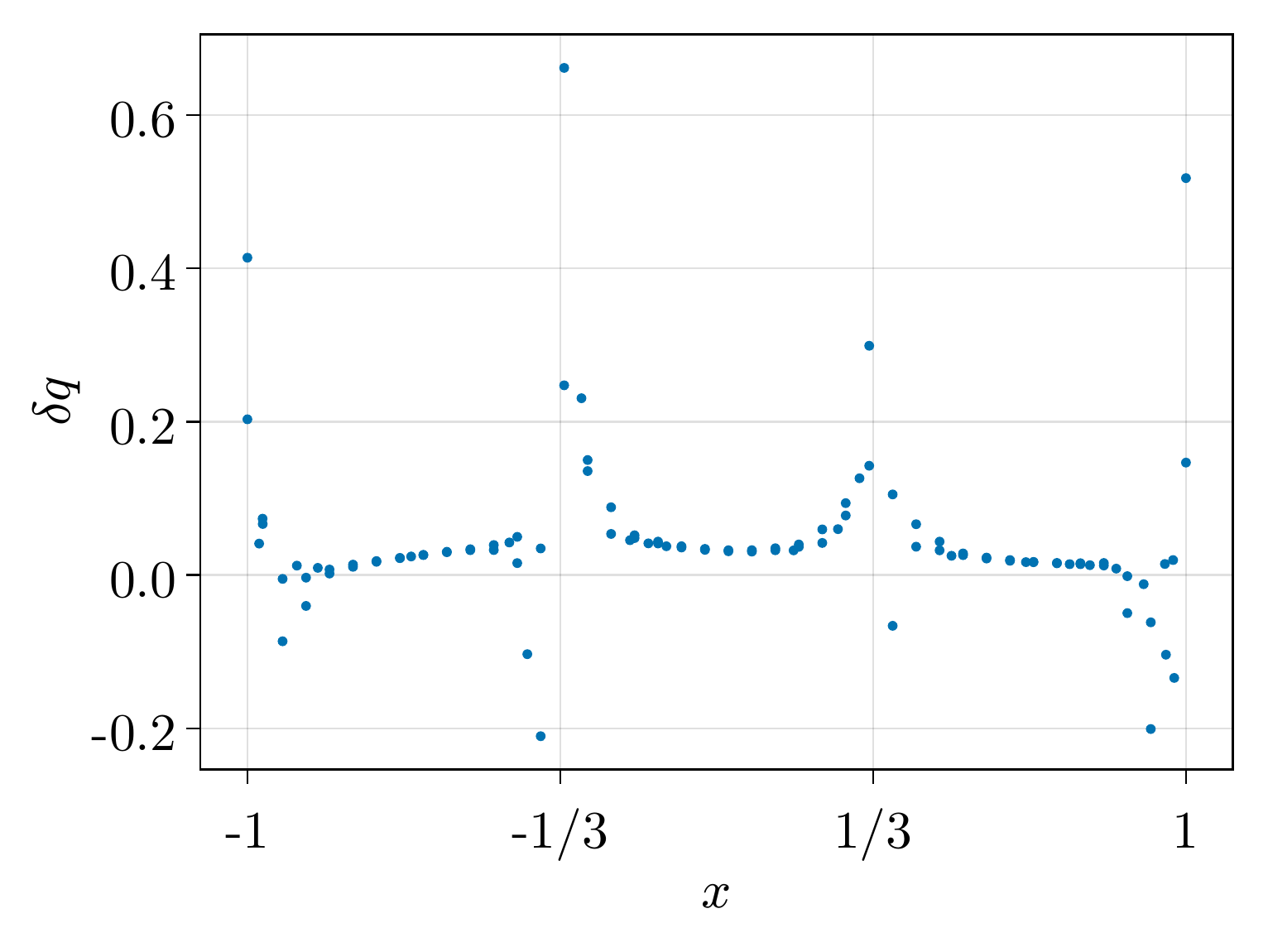}
        \caption{Conservative Neumann boundaries}
        \label{fig:etagrad_three_strip_conservative}
    \end{subfigure}
    \caption{Profiles of $\delta q = \eta\partial_x (u-u_h)$ at $y \approx 0$ for the three-strip test case of hybrid method with and without conservative Neumann boundary conditions for a jump $\delta\eta_L = \num{1e6}$ and $\delta\eta_R = \num{1e-4}$ on an irregular point cloud.}
    \label{fig:etagrad_three_strip}
\end{figure}

In the previous experiment with the interior interface test case, we observed that the jump propagation across interfaces works best with the conservative method. However, problems arise from the Neumann boundary where false jump conditions are propagated. The goal of this test case is to examine if the conservative Neumann boundary conditions can be used to improve numerical results. We use conservative Neumann boundary conditions in the neighborhood of points where a switch to the conservative formulation is performed. We observe that the standard hybrid method does not converge for this test case as can be seen in \cref{fig:convergence_three_strip}. A view on \cref{sub@fig:solutions_three_strip_regular} reveals that the analytical solution is not reached due to the lack of flux conservation \cite{Suchde_Kuhnert_Schroeder_Klar_2017}. For a numerical solution $u_h$ let us define the error in flux by \[ \delta q = \eta \partial_x (u - u_h). \] According to \cref{eq:jump_conditions}, the flux error $\delta q$ has to be piecewise constant and continuous which leads to it being constant throughout the domain. In \cref{fig:etagrad_three_strip_regular} we see that the flux error $\delta q$ for the numerical solution without conservative Neumann boundary conditions is not constant. This observation can be made especially at the boundaries $\partial\Omega_T$ and $\partial\Omega_B$ which leads to the assumption that flux conservation has to be met at the boundary as well. Finally, switching to conservative Neumann boundary conditions, we see in \cref{fig:solutions_three_strip_conservative,fig:etagrad_three_strip_conservative} that the analytical solution matches the numerical solution visually and $\delta q$ has lower oscillations. \Cref{fig:convergence_three_strip} confirms that the method is first-order accurate with conservative Neumann boundary conditions.

\begin{figure}
    \centering
    \includegraphics[width=0.49\linewidth]{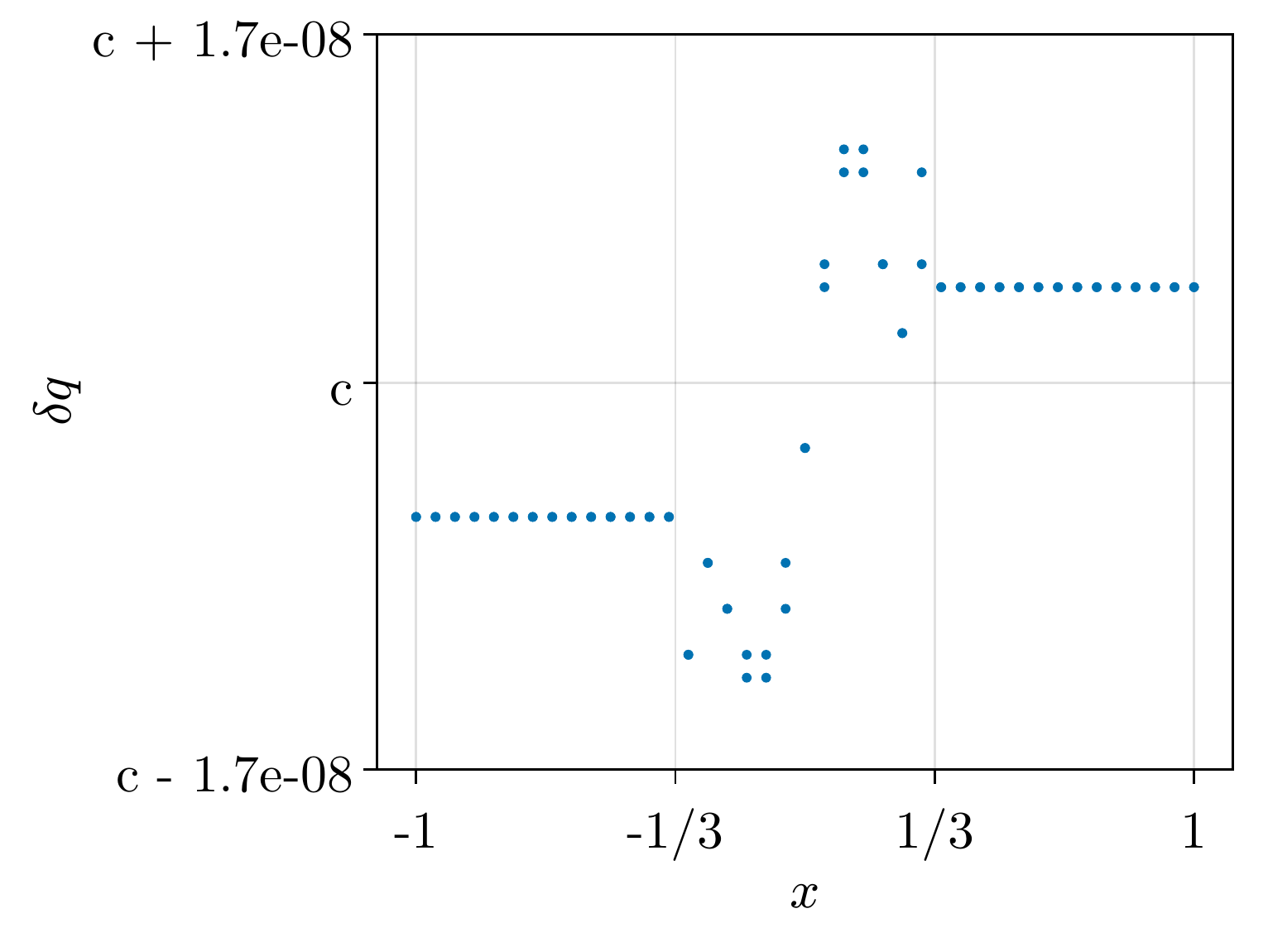}
    \caption{Profile of $\delta q = \eta\partial_x (u-u_h)$ at $y \approx 0$ for the three-strip test case of hybrid method with conservative Neumann boundary conditions for a jump $\delta\eta_L = \num{1e6}$ and $\delta\eta_R = \num{1e-4}$ on a uniform point cloud. The constant $c \approx \num{-1.834e-02}$ was numerically determined.}
    \label{fig:etagrad_three_strip_uniform}
\end{figure}

To better illustrate $\delta q$ in \cref{fig:etagrad_three_strip}, a domain decomposition based on $\eta$ was executed to be able to calculate the gradient $\nabla u|_{\Omega_*}$ on each subdomain $\Omega_L$, $\Omega_M$ and $\Omega_R$ without the influence of the neighboring subdomains. The numerical results were computed without a domain decomposition. Note that the outliers on \cref{fig:etagrad_three_strip_conservative} are due to the irregularity of the point cloud and the circumstance that we calculate the gradient near the interface which serves as a boundary to the subdomains. Computing $\delta q$ on a uniform point cloud, we obtain \cref{fig:etagrad_three_strip_uniform} and thus the flux error at each point is given by \[ \delta q_i = c + \eps_i \]
for a numerically computed constant $c \approx \num{-1.834e-02}$ and an error \[\max_{i = 1,\dots,N}\abs{\eps_i} = O(\num{1e-8})\] which is negligible.

    \section{Discussion} \label{sec:discussion}
    In this section, we give a summary of the results from the previous section and discuss the advantages and disadvantages of the new hybrid method.

Our first observation was that artificial smoothing could stabilize the convergence behavior of the strong form method for very simple test cases, such as the two-strip problem. However, we have also seen other test cases where the strong form method failed to provide accurate solutions, regardless of whether we smooth or not. Thus, the first consequence is that artificial smoothing should be avoided because it distorts the data from the mathematical model and does not yield advantages for more complex test cases.

The hybrid method produces good results and provides first-order convergence for elliptic interface problems. In the performance benchmark in \ref{app:performance}, we see that the computation time does not increase compared to the fully strong form method, due to a lower number of iterations in the linear solver. Additionally, our method guarantees the M-matrix property and thus positivity of numerical results and fulfills the discrete divergence theorem locally. However, one limitation is the need for computing a Voronoi cell. We observed that the presented node selection routine only works well with high discontinuities in the diffusion coefficient, as seen in \cref{ssec:two_strip,ssec:curved_interface}. It is a major drawback that the node selection depends on the jump magnitude and not the mere existence of a jump in the diffusion coefficient. On top of that, it is also susceptible to low point cloud quality, see \ref{app:node_selection}, but this issue can be managed by better point cloud management routines.

Finally, we saw that placing conservative Neumann boundary conditions is necessary for the {three-strip} test case.

    \section{Conclusion} \label{sec:conclusion}
    In this paper, we presented two formulations of discretizing the diffusion operator and a new way of combining these approaches in a locally conservative and positivity preserving method. Our approach combines the cell average-based conservative formulation and the strong form weighted least squares-based formulation into a new method. It puts the finite volume method idea inside the context of GFDM and uses its advantages without sacrificing the advantages that GFDM offers, namely the flexibility of the point cloud. This is a major benefit with respect to time-dependent problems in a Lagrangian formulation with a moving point cloud, which will be an extension of the present work.

We performed an in-depth parameter study with different test cases and methods. We considered jumps up to ten orders of magnitude whereas existing methods in literature were limited to three orders of magnitude. Our numerical results show that the classical weighted least squares approach fails for high jump magnitudes and that using a flux-conserving scheme overcomes these problems. Finally, we have seen that it can be necessary to incorporate conservative Neumann boundary conditions to guarantee the conservation property of the entire numerical scheme. Since we carried out the computations on non-conforming point clouds, we showed that the presented method does not depend on aligning the point cloud to the interface.

The discretization of the diffusion operator is subject to ongoing research and as such open questions remain. For the diffusion operator, it is possible to use the notion of shapeless volumes and higher-order GFDM methods. Additionally, the jump identification technique can be extended for it to be able to identify lower jumps. While this is not needed for our applications, it might be necessary for others.

    \section{Acknowledgements}
    Pratik Suchde would like to acknowledge support from the European Union’s Horizon 2020 research and innovation program under the Marie Skłodowska-Curie Actions grant agreement No. 892761 \enquote{SURFING}.

    \appendix
    \section{Selection of nodes for conservative method} \label{app:node_selection}

\begin{figure}
        \newcommand{\subfigurehspace}{.49\linewidth}
        \begin{subfigure}{\subfigurehspace}
            \centering
            \includegraphics[width=\linewidth]{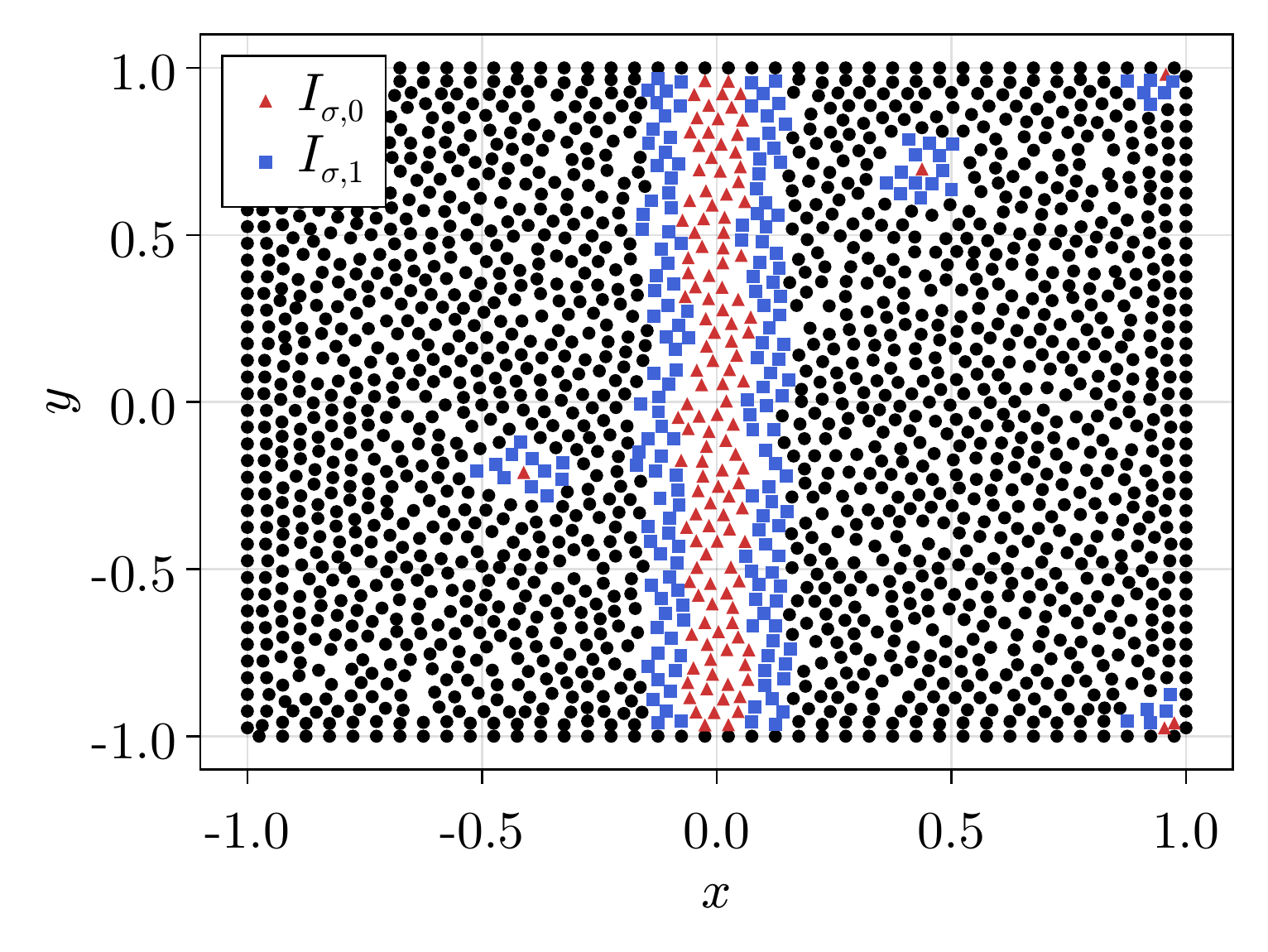}
            \caption{Global selection}
            \label{fig:sigma_global}
        \end{subfigure}
        \begin{subfigure}{\subfigurehspace}
            \centering
            \includegraphics[width=\linewidth]{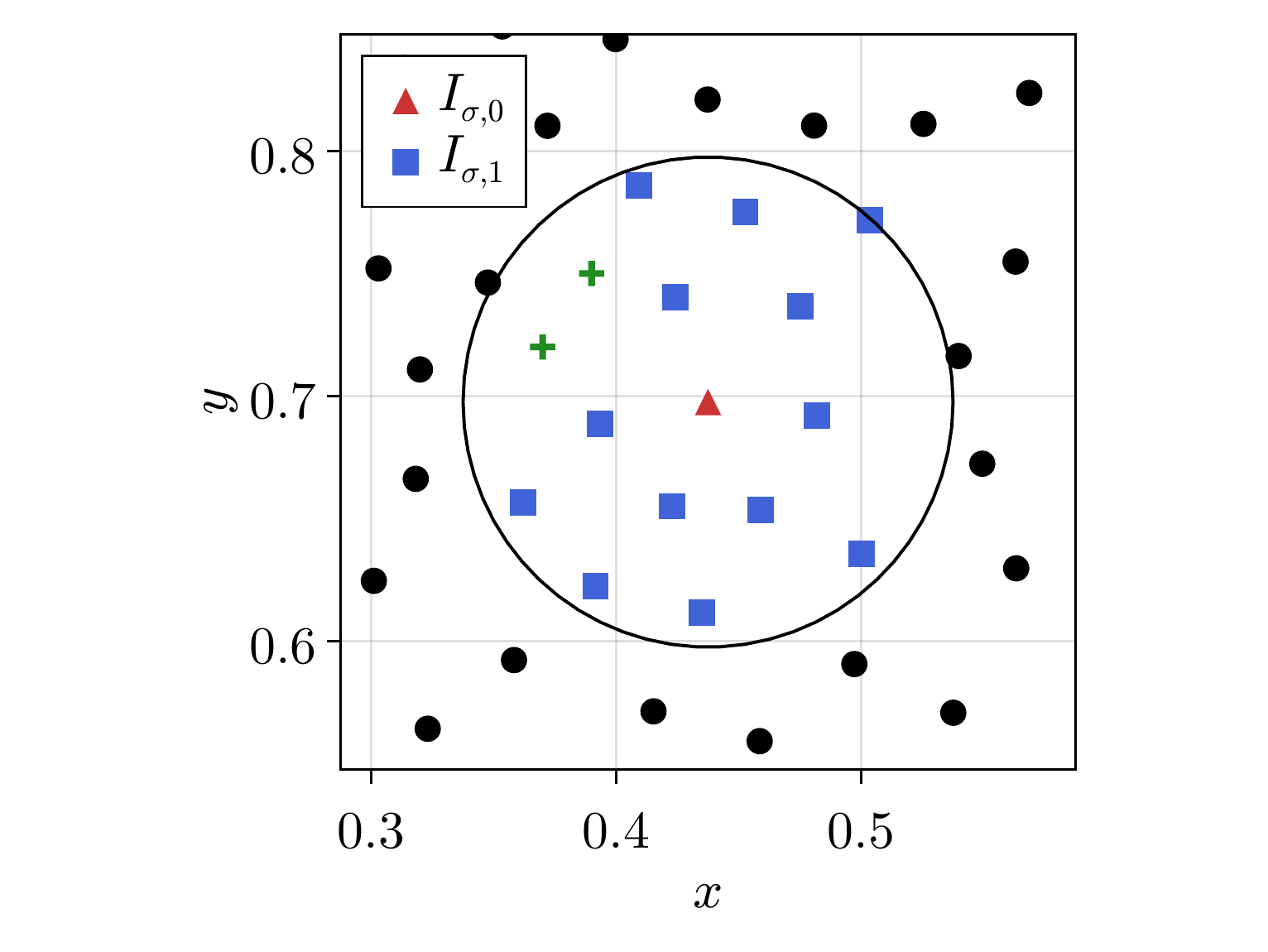}
            \caption{Local selection due to irregularity}
            \label{fig:sigma_irregular}
        \end{subfigure}
        \caption{Nodes selected for the hybrid regime by the diagonal dominance criterion for an irregular point cloud with $h = \num[]{5e-2}$. Nodes whose addition lead to a diagonally dominant operator are marked with a green plus sign.}
        \label{fig:sigma}
    \end{figure}

    In this section, we discuss the selection of nodes for the conservative method. For that, we define the index set $I$ of interior points. We denote the subset of points where the diagonal dominance criterion \eqref{eq:hybrid_condition} triggers by
    \[ I_{\sigma, 0} = \set{i \in I \mid \sigma_i > 0} \]
    and the index sets of their neighbors by
    \[ I_{\sigma, 1} = \set{i \in I\setminus I_{\sigma, 0} \mid S_i \cap I_{\sigma, 0} \ne \emptyset}. \]
    Note that for $I_{\sigma, 0}$, the criterion $\sigma_i > 0$ is unpractical since numerical errors are not factored in and thus the majority of points are selected for the conservative scheme. To overcome this drawback, we provide a tolerance $\epsilon > 0$, in our case $\epsilon = \num{1e-12}$, and use the criterion $\sigma_i > \epsilon$ instead. We tested other values $\epsilon \in (\num{1e-12}, \num{1e-4})$ but a remarkable difference in the selected nodes could not be observed. In \cref{fig:sigma_global} the index sets are depicted for a jump of $\delta\eta = \num{1e4}$ and we can see that $I_{\sigma, 0}$ generally conforms to the interface \[ \Gamma = \set{(x,y)^T \in \R^2 \mid x = 0} \] well. But we also see outliers that are in the region where the diffusivity $\eta$ is smooth, one of which is shown in \cref{fig:sigma_irregular}. These outliers are most likely due to the irregularity or locally low quality of the point cloud. The addition of two points at the locations $(0.37, 0.72)$ and $(0.39, 0.75)$ in the top-left empty region of the neighborhood, marked as a plus in \cref{fig:sigma_irregular}, leads to a diagonally dominant row.

    In the following, we will discuss the fraction of nodes that fall under the hybrid regime for the two-strip test case with varying jump magnitude and smoothing length.  For this, we define the index set of nodes that are chosen for the conservative scheme by \[ I_c = I_{\sigma, 0} \cup I_{\sigma, 1} \]. To define the set of nodes that are close to the interface, we exploit the fact that $\eta$ is piecewise constant and define the proximity to the interface locally through the neighborhoods
    \[ I_\Gamma = \set{ i \in I \mid \exists j\in S_i : \eta_i \ne \eta_j }. \]
    Finally, we define the set of nodes that consist of all the neighborhoods of points that are close to the interface
    \[ I_\Gamma^+ = \bigcup_{i \in I_\Gamma} S_i \cap I. \]

    \begin{table}
        \centering
        \begin{tabular}{c|cccc|cccc}
            & \multicolumn{4}{c|}{$h = \num{5e-2}$} & \multicolumn{4}{c}{$h = \num{1.25e-2}$} \\
            \cline{2-9}
            $\delta\eta$ & $\frac{\abs{I_{\sigma, 0}}}{\abs{I}}$ & $\frac{\abs{I_c}}{\abs{I}}$ & $\frac{\abs{I_c \cap I_\Gamma}}{\abs{I_c}}$ & $\frac{\abs{I_c \cap I_\Gamma^+}}{\abs{I_c}}$ & $\frac{\abs{I_{\sigma, 0}}}{\abs{I}}$ & $\frac{\abs{I_c}}{\abs{I}}$ & $\frac{\abs{I_c \cap I_\Gamma}}{\abs{I_c}}$ & $\frac{\abs{I_c \cap I_\Gamma^+}}{\abs{I_c}}$ \\
            \hline
            \num{1e+00} & \SI{0.28}{} & \SI{2.21}{}  & \SI{0.00}{}  & \SI{0.00}{}  & \SI{0.13}{} & \SI{1.57}{} & \SI{0.00}{}  & \SI{0.00}{} \\
            \num{1e+02} & \SI{5.13}{} & \SI{15.72}{} & \SI{51.58}{} & \SI{85.96}{} & \SI{1.35}{} & \SI{4.96}{} & \SI{41.27}{} & \SI{69.19}{} \\
            \num{1e+04} & \SI{7.00}{} & \SI{17.43}{} & \SI{46.52}{} & \SI{87.34}{} & \SI{1.87}{} & \SI{5.43}{} & \SI{37.64}{} & \SI{71.89}{} \\
            \num{1e+06} & \SI{7.78}{} & \SI{18.15}{} & \SI{44.68}{} & \SI{87.84}{} & \SI{2.04}{} & \SI{5.59}{} & \SI{36.56}{} & \SI{72.70}{} \\
            \num{1e+08} & \SI{8.11}{} & \SI{18.53}{} & \SI{43.75}{} & \SI{88.10}{} & \SI{2.12}{} & \SI{5.67}{} & \SI{36.08}{} & \SI{73.06}{} \\
            \num{1e+10} & \SI{8.33}{} & \SI{18.70}{} & \SI{43.36}{} & \SI{88.20}{} & \SI{2.16}{} & \SI{5.69}{} & \SI{35.93}{} & \SI{73.17}{} \\
        \end{tabular}
        \caption{Node fractions in \% for hybrid method with different smoothing lengths $h$ and jump magnitudes $\delta\eta$.}
        \label{tab:fractions}
    \end{table}

    We summarized different node fractions in \cref{tab:fractions}. First, let us take a look at the fraction of nodes that are selected for the hybrid regime by the $\sigma_i > 0$ criterion, given by $\frac{\abs{I_{\sigma, 0}}}{\abs{I}}$. In the table, we observe that even without a jump, some points are selected for the hybrid regime. By increasing the jump, we observe that the number of points that do not fulfill the diagonal dominance criterion rises but appears to be bounded for both smoothing lengths. As expected, the fraction of such nodes decreases with a decreasing smoothing length $h$ and thus increasing number of points. The fraction of points chosen for the conservative scheme, indicated by $\frac{\abs{I_c}}{\abs{I}}$, displays the same behavior.

    In the next step, we take a look at the fraction of nodes that are selected for the hybrid regime and those that are actually within the region of the jump, given by $\frac{\abs{I_c \cap I_\Gamma}}{\abs{I_c}}$. We observe that for an increasing jump, the number of points that use the hybrid method that are within the interface region decreases, meaning that for an increasing jump, we select points that are not within the interface region. This can be explained by points $\bvec{x}_i$ that have one or two points within their neighborhoods $\bvec{x}_j \in B_i$ that have the other $\eta$-value, meaning $\eta_i \ne \eta_j$. For low jumps, it is possible for such points to obtain a diagonally dominant discrete operator; with an increasing jump however, these few neighbors have an increasing influence on the numerical gradient which then leads to instabilities in the discrete operator. This assumption can be confirmed by taking a look at the ratio $\frac{\abs{I_c \cap I_\Gamma^+}}{\abs{I_c}}$ that indicates the fraction of points selected for the hybrid regime that also lie in the extended proximity of the interface $I_\Gamma^+$. This ratio is increasing with an increasing jump magnitude.

\section{Performance benchmarks} \label{app:performance}
It is reasonable to assume that the hybrid method presented in \cref{ssec:hybrid} would decrease the performance of the numerical solver as compared to using the strong form method from \cref{ssec:classical} exclusively. Not only do we have to compute a Voronoi cell for the selected points, but we are required to compute the conservative operator in addition to the strong form one. In our self-written Julia code, a global Delaunay triangulation of the point cloud is computed. Since the Voronoi tesselation is the dual of the Delaunay triangulation, we can easily reconstruct the Voronoi cells from the triangulation by connecting the circumcenters of each triangle. For large-scale 3D applications, it is generally not feasible to compute a global Delaunay triangulation. In the software suite MESHFREE \cite{MESHFREE}, a Delaunay triangulation is computed locally for point cloud management purposes. Computing Voronoi cells from an underlying (local) Delaunay triangulation is thus not that expensive compared to constructing Voronoi cells directly.

We measured the runtime of the self-written numerical solver in Julia \cite{Bezanson_Edelman_Karpinski_Shah_2017} on an Intel\textsuperscript{\circledR} Core\textsuperscript{\texttrademark} i5-8600K processor from the start until the end of the simulation. We used the interior interface test case from \cref{ssec:interior_interface} with a jump $\delta\eta = \num{1e10}$. The global Delaunay triangulation was only computed for the hybrid method using a Julia wrapper of Triangle \cite{Shewchuk_1996}. The global linear system \cref{eq:poisson_discrete} was solved using BiCGstab(2) \cite{Sleijpen_1993} and an incomplete LU preconditioning. For each smoothing length, we computed the mean elapsed time of a sample of ten simulations. The results are summarized in \cref{tab:performance_comparison}.

We can see that there is no performance drawback of the new method, the hybrid method is even faster in most cases. This can be explained due to the Delaunay triangulation being not very complex in 2D on a square domain. On the finest point cloud, the computation of the Delaunay triangulation took only \SI{1}{\second} which is roughly \SI{3}{\percent} of the total computation time. However, the hybrid method needed fewer BiCGstab iterations and thus can compensate for the lost time due to the Delaunay triangulation and additional computation of the conservative operators.

\begin{table}
    \centering
    \begin{tabular}{r|rr|rr}
        & \multicolumn{2}{|c|}{Strong Form} & \multicolumn{2}{|c}{Hybrid} \\
        \hline
        $h$ & $N_{\text{iter}}$ & Mean elapsed time & $N_{\text{iter}}$ & Mean elapsed time \\
        \num{2e-1} & 1 & \SI{28.2}{\milli\second} & 2 & \SI{28.4}{\milli\second} \\
        \num{1e-1} & 2 & \SI{92.7}{\milli\second} & 3 & \SI{88.7}{\milli\second} \\
        \num{5e-2} & 3 & \SI{136}{\milli\second} & 2 & \SI{119}{\milli\second} \\
        \num{2.5e-2} & 4 & \SI{376}{\milli\second} & 3 & \SI{355}{\milli\second} \\
        \num{1.25e-2} & 9 & \SI{1.76}{\second} & 6 & \SI{1.62}{\second} \\
        \num{6.25e-2} & 18 & \SI{7.27}{\second} & 11 & \SI{6.98}{\second} \\
        \num{3.125e-2} & 33 & \SI{36.21}{\second} & 21 & \SI{33.93}{\second}
    \end{tabular}
    \caption{Number of BiCGstab iterations $N_{\text{iter}}$ and mean elapsed time for the strong form and hybrid methods with different smoothing lengths $h$.}
    \label{tab:performance_comparison}
\end{table}

    \bibliographystyle{elsarticle-num-names}
    % \bibliography{bibliography_bibtex.bib}

\end{document}